\documentclass[a4,11pt]{article}
\usepackage{amssymb}
\usepackage{amsmath}
\usepackage{amsthm}
\usepackage{bbm}
\usepackage{graphicx}
\usepackage{tikz}
\usepackage{psfrag}
\usepackage[all]{xy}
\usepackage{supertabular,multicol}
\usepackage{authblk}

\usepackage[totalwidth=17cm,totalheight=24.2cm]{geometry}

\makeatletter
\let\mcnewpage=\newpage
\newcommand{\Trick}{
\renewcommand\newpage{%
        \if@firstcolumn
            \hrule width\linewidth height0pt
                \columnbreak
        \else
                \mcnewpage
        \fi
}
}
\makeatother

\theoremstyle{plain}

\newtheorem{theorem}{Theorem}[section]
\newtheorem{proposition}[theorem]{Proposition}
\newtheorem{corollary}[theorem]{Corollary}
\newtheorem{lemma}[theorem]{Lemma}

\newtheorem{maintheorem}{Theorem}

\newtheorem{mainlemma}[maintheorem]{Lemma}

\theoremstyle{definition}
\newtheorem{definition}[theorem]{Definition}
\newtheorem{example}[theorem]{Example}
\newtheorem{remark}[theorem]{Remark}
\newtheorem{problem}[theorem]{Problem}
\newtheorem{question}[theorem]{Question}

\newtheorem*{notation}{Notation}

\makeatletter

\newenvironment{figurehere}
  {\def\@captype{figure}}
  {}
\makeatother

\def\tr{{\rm tr}\,}

\def\th{\vartheta}
\def\ZZ{\mathbb{Z}}

\def\RR{\mathbb{R}}
\def\CC{\mathbb{C}}
\def\PP{\mathbb{P}}

\def\SU{\mathrm{SU}}
\def\SL{\mathrm{SL}}
\def\SO{\mathrm{SO}}
\def\PSL{\mathrm{PSL}}
\def\dev{\mathrm{dev}}

\def\Graph{{G}}

\def\Sph{\mathbb{S}^2}
\def\Spt{\mathbb{S}^3}

\def\ol#1{\overline{#1}}

\def\bm#1{\text{\boldmath$#1$}}

\def\vc#1{{\bm{#1}}}   

\def\rar{\rightarrow}
\def\lra{\longrightarrow}
\def\pa{\partial}
\def\1{\vc{1}}
\def\0{\vc{0}}

\newcommand{\sqdiamond}[1][fill=white]{\tikz [x=1.2ex,y=1.85ex,line width=.1ex,line join=round, yshift=-0.285ex] \draw  [#1]  (0,.5) -- (.8,1) -- (1.6,.5) -- (.8,0) -- (0,.5) -- cycle;}%

\def\cube{{\Box}}
\def\tcube{{\sqdiamond}}
\def\hcube{\triangle}

\def\Aang{\mathcal{A}}
\def\Hang{\mathcal{H}}
\def\Pang{\mathcal{P}}
\def\Nang{\mathcal{N}}
\def\Kang{\mathcal{K}}
\def\Fang{\mathcal{F}}
\def\inte{\mathrm{int}}

\begin{document}

\title{Spherical metrics with conical singularities\\
on a $2$-sphere: angle constraints}

\author[1]{Gabriele Mondello}
\affil[1]{{\small{``Sapienza'' Universit\`a di Roma, Department of Mathematics (mondello@mat.uniroma1.it)}}}

\author[2]{Dmitri Panov}
\affil[2]{{\small{King's College London, Department of Mathematics (dmitri.panov@kcl.ac.uk)}}}

\maketitle

\abstract{\noindent In this article we give a criterion for
the existence of a metric of curvature $1$ on a $2$-sphere
with $n$ conical singularities of prescribed angles $2\pi\th_1,\dots,2\pi\th_n$ and non-coaxial holonomy.
Such a necessary and sufficient condition is expressed in terms
of linear inequalities in $\th_1,\dots,\th_n$.}

\tableofcontents

\setlength{\parskip}{0.2\baselineskip}

\section{Introduction}


\subsection{Formulation of the problem}

The aim of this paper is to prove the existence of Riemannian metrics of curvature $1$ with $n\geq 2$ conical singularities of assigned angles on a compact connected orientable surface of genus $0$.

In order to state the problem more formally, we recall some terminology.

\begin{notation}
By {\it{surface}} we always mean a smooth $2$-dimensional real manifold, possibly with boundary, and we will call {\it{a sphere}} just a compact connected
orientable surface diffeomorphic to $S^2$.
By {\it{metric}} we always mean a Riemannian metric and
by {\it{spherical metric}} just a Riemannian metric of constant curvature $1$. 
We keep the notation $\Sph$ for the unit $2$-sphere endowed with the standard metric.
\end{notation}

From the local point of view, spherical metrics are easy to describe: it is a classical result (see Killing \cite{killing:constant} and Hopf \cite{hopf:constant}) that a surface endowed with a spherical metric
is locally isometric to a portion of $\Sph$.

From a global point of view spherical metrics on compact connected orientable surfaces only exist in genus $0$ by Gauss-Bonnet; moreover, in this case they are all isometric to each other.

The situation becomes more interesting if we allow our metrics to admit conical singularities. 

\begin{definition}
A Riemannian metric $g$ of curvature $1$ on a surface $S$ has a
{\it{conical singularity
of angle $2\pi\alpha>0$ at $y\in S$}} if it can be locally written as $g=dr^2+\alpha\sin(r)d\theta^2$, where $(r,\theta)$ are local polar coordinates on $S$ centered at $y$. We will say that the angle is {\it{integral}} if
$\alpha\in\ZZ$.
\end{definition}

Our goal is to answer the following question.

\begin{question}[Existence of metrics]\label{question:intro} 
For which $\vc{\th}=(\th_1,\dots,\th_n)\in\RR_+^n$
there exists a spherical metric $g$ on a sphere $S$ with $n$ conical singularities 
of angles $2\pi\cdot\vc{\th}=(2\pi \th_1,\dots,2 \pi \th_n$)?
\end{question}


\subsection{Context and known results}

We remark that Question \ref{question:intro} is actually different from the following more classical problem.

\begin{question}[Existence of conformal metrics]\label{question:conformal}
Fix a connected Riemann surface $S$ with distinct points $x_1,\dots,x_n$.
For which $\vc{\th}=(\th_1,\dots,\th_n)\in\RR_+^n$ there exists
a {\it{conformal}} metric of constant curvature on $S$ with
conical singularity at $x_i$ of angle $2\pi\th_i$?
\end{question}

Notice that the curvature
of the desired metric in Question \ref{question:conformal}
must have the same sign as $\chi(S,\vc{\th}):=\chi(S)+\sum_i (\th_i-1)$.

\begin{remark}
Both problems described above can be formulated in terms of moduli spaces of metrics
$\mathfrak{Met}(S,x,\vc{\th})$ of constant curvature on the surface $S$ with conical singularities $x_1,\dots,x_n$ of assigned angles $2\pi\cdot\vc{\th}$ (up to isotopies that fix the singularities). Thus, Question \ref{question:intro} can be rephrased in terms of non-emptiness of such $\mathfrak{Met}(S,x,\vc{\th})$; on the other hand, Question \ref{question:conformal} asks whether the map $\mathfrak{Met}(S,x,\vc{\th})\rar\mathfrak{T}(S,x)$ to Teichm\"uller space that remembers the underlying conformal structure is onto.
In this paper we will not push this point of view farther.
\end{remark}

For $n=0$ and $\chi(S)\geq 0$ it is a standard fact 
that in every given conformal class there exists a
metric of constant curvature and that such a metric is unique up to rescaling and conformal automorphisms of $S$;
whereas for $\chi(S)<0$ such an existence and uniqueness statement is provided by the classical uniformization theorem proven by Koebe \cite{koebe:I}-\cite{koebe:II} and Poincar\'e \cite{poincare:uniformization}.

Assume now $n>0$.
Existence and uniqueness results were proven
by Thurston \cite{thurston:shapes} 
and Troyanov \cite{troyanov:euclidean}
for $\chi(S,\vc{\th})=0$
and by McOwen \cite{mcowen:conical} 
and Troyanov \cite{troyanov:conical}
for
$\chi(S,\vc{\th})<0$.

As for the case $\chi(S,\vc{\th})>0$,
existence and uniqueness still holds in
the subcritical case (and so in particular
when all angles are smaller than $2\pi$)
again by Troyanov \cite{troyanov:conical}.
On the other hand, 
it is known that such uniqueness does not hold any more
in the supercritical case.
For instance,
an existence theorem was proven by Bartolucci-De Marchis-Malchiodi \cite{BDMM:supercritical} for $\chi(S)\leq 0$ and a lower bound for the number of such metrics is also provided.
Notice that the general case of $\chi(S,\vc{\th})>0$ and
$\vc{\th}\notin (0,1]^n$ is not covered by the above works.

Another manifestation of the non-uniqueness of the solution is provided by Scherbak \cite{scherbak:critical}, who counted
the exact number of such metrics 
for almost all configurations of $x_1,\dots,x_n\in S$, when
$S$ is a sphere, $\chi(S,\vc{\th})>0$ and 
all $\th_i$ are integers.

\bigskip

Back to Question \ref{question:intro}, it is easy to see that the only
possibility for $n=1$ is a surface isometric to the standard $\Sph$.
An answer to this question is also known for 
$n=2$ by work of Troyanov \cite{troyanov:bigons} and for
$n=3$ by work of Eremenko \cite{eremenko:three}. 
A detailed analysis of spherical polygons with two non-integral angles
is done in \cite{EGT:polygons}, and more extensively for $n=4$ in
\cite{EGT:four}; spherical quadrilaterals with three non-integral angles are studied in \cite{EGT:polygons3}.

In this paper we will give an {\it almost complete} answer to this question for $n\ge 4$.

\subsection{Main results}

%
Our first main result shows that the existence of a spherical metric
on a sphere $S$ with conical singularities of angles
$2\pi\cdot\vc{\th}$ imposes
some restrictions on the vector $\vc{\th}=(\th_1,\dots,\th_n)$.

%

\begin{notation}
Denote by $\|\cdot\|_1$ the standard $\ell^1$ norm on $\RR^n$ and by $d_1$ the associated $\ell^1$ distance, and let
$\ZZ^n_o$ denote the subset of {\it{odd}} points of $\RR^n$,
namely of points $\vc{m}=(m_1,\dots,m_n)$ in $\ZZ^n$ such that
$\|\vc{m}\|_1$ is odd.
\end{notation}

\begin{maintheorem}[Holonomy constraints]\label{nonexist}
Suppose there exists a sphere $S$ with a spherical metric
with conical singularities of angles
$2\pi\th_1,\dots,2\pi\th_n$. Then the 
following inequalities hold:
\begin{align}
\left.
\begin{array}{r@{\hspace{0.5em}}l}
\vc{\th} & >\phantom{-}0 \\
\sum_{i=1}^n (\th_i-1)
& >-2
\end{array} \right\} \tag{P}\label{gauss-bonnet}\\
\begin{array}{r@{\hspace{0.5em}}l}
d_1(\vc{\th}-\vc{\1}, \ZZ^n_o)  & \ge \phantom{-}1 
\end{array}\phantom{\Big\}}
\tag{H}\label{theinquality}
\end{align}
where $\vc{\1}=(1,1,\dots,1)\in\RR^n$.\\
Moreover, if equality in \eqref{theinquality} is attained, then the holonomy of the metric is coaxial.
\end{maintheorem}

\begin{remark}
The {\it{positivity constraints}}
\eqref{gauss-bonnet} follow from the positivity
of the angles and the positivity of the area,
via the Gauss-Bonnet theorem.
\end{remark}

As the set of points for which
the {\it{holonomy constraints}} \eqref{theinquality} do not hold is the union of disjoint octahedrons, we also have the following.

\begin{mainlemma}[Connectedness]\label{lemma:connected}
The set of points in $\RR^n$ that strictly satisfy
the holonomy constraints \eqref{theinquality} is non-empty for $n\ge 3$ and connected for $n\ge 4$. The same holds for the subset of
points that satisfy the positivity constraints \eqref{gauss-bonnet} and
the holonomy constraints \eqref{theinquality} strictly.
\end{mainlemma}

The proof of Theorem \ref{nonexist} consists of a few steps.
We first
associate to each spherical 
metric on $S$ with conical singularities $x_1,\dots,x_n$
the holonomy representation $\rho$ of the free group $\pi_1(S\setminus \{x_1,\cdots,x_n\})$ in $\SO(3,\RR)$. Then we
show that, since $S$ is a sphere,
such holonomy representation admits a canonical lift $\hat{\rho}$ to $\SU(2)$. Thus, we
relate representations $\pi_1(S\setminus\{x_1,\dots,x_n\})\rar\SU(2)$ to closed broken geodesics on $\Spt$ and we
verify that the wished closed broken geodesic on $\Spt$ exists if and only if Inequalities \eqref{theinquality} are satisfied. This explains
the name ``holonomy constraints''.

A special role will be played by ``generic'' holonomy representations, namely whose image does not belong to a $1$-parameter subgroup of $\SO(3,\RR)$.

\begin{definition}
A representation $\rho$ in $\SO(3,\RR)$ is {\it{coaxial}} if
its image consists of rotations about the same axis.
\end{definition}

%

The second main result of this paper is the
following partial converse to Theorem \ref{nonexist}.

\begin{maintheorem}[Existence of spherical metrics]\label{existence}
Let $\th_1,\dots,\th_n$ be real numbers that satisfy both
the positivity constraints \eqref{gauss-bonnet} and
the holonomy constraints \eqref{theinquality} strictly.
Then there exists a sphere $S$ with a spherical metric
with conical points of angles $2\pi\th_1,\dots,2\pi\th_n$
and non-coaxial holonomy.
\end{maintheorem}

\begin{remark}
The cases that are not covered by this theorem is when (\ref{theinquality}) becomes an equality,
when 
the holonomy of such a spherical metric is necessarily coaxial
(provided such a metric exists!).
\end{remark}

In order to prove Theorem \ref{existence}, we proceed as follows.

First we construct such metrics for $n=2,3,4$ (the cases $n=2$ and $n=3$ were previously treated by Troyanov and Eremenko respectively).

The idea is then to inductively produce metrics with $n\geq 5$ conical points close to degenerate ones by picking a spherical metric with fewer
singularities and splitting a conical point.
%
%
More precisely, given $\th_1,\dots,\th_{n}$ as in Theorem \ref{existence}, we show that the wished spherical metric on $S$
can be obtained starting from a spherical metric on $S'$ with $n-1$
conical singularities by operating a surgery in a neighbourhood
of a conical point. Typically, the surgery will produce two points of angles
$2\pi\th_i$, $2\pi\th_j$ very close to each other on $S$ out of a single conical point of angle $2\pi(\th_1\pm\th_j-1+\eta)$ on $S'$, where $\eta$ is a small fee that we have to pay for the performed cut-and-paste operation.
In order to take care of this little $\eta$, we use a deformation result of Luo \cite{luo:monodromy}.

\begin{remark}
The presence of such a possibly nonzero $\eta$ is what
forces us to require that the holonomy constraints \eqref{theinquality} are satisfied {\it{strictly}}.
\end{remark}

Finally, the combinatorial result that tells us which conical point to split is the following.

\begin{maintheorem}[Algebraic merging]\label{linearalgebra} 
Assume $n\geq 5$.
Let $\th_1,\dots,\th_n$ be real numbers 
that satisfy both
the positivity constraints \eqref{gauss-bonnet}
and
the holonomy constraints
\eqref{theinquality} strictly.
Then there is a choice of distinct $i,j\in\{1,\dots,n\}$
such that at least
one of the following two $(n-1)$-tuples
\[
\left(\th_1,\dots,\hat{\th}_i,\dots,\hat{\th}_j,\dots,\th_n, (\th_{i}+ \th_{j}-1)\right),
\qquad
\left(\th_1,\dots,\hat{\th}_i,\dots,\hat{\th}_j,\dots,\th_n, (\th_{i}- \th_{j}-1)\right)
\]
satisfies positivity constraints and holonomy constraints strictly.
\end{maintheorem}


\subsection{Structure of the paper}

The paper is divided into two parts: the former deals with the holonomy constraints and the latter provides the actual geometric constructions of the wished metrics.

\medskip

In Section \ref{sec:rep} Theorem \ref{nonexist} is proven.
More precisely,
we first recall some well-known facts about the developing map and the holonomy representation associated to a spherical metric.
Then we prove that such a representation admits a natural lift to $\SU(2)$ and that such a lift is canonical, if we are working with a sphere. In the remaining subsections we show that representations in $\SU(2)$ carry (almost) the same information as closed broken geodesics on $\Spt$, whose existence is equivalent to the holonomy constraints \eqref{theinquality}. The case of Abelian and coaxial holonomy is briefly discussed.

In Section \ref{sec:algebraic-merging} the holonomy and positivity constraints are studied from an algebraic point of view
and Lemma \ref{lemma:connected} and Theorem \ref{linearalgebra} are proven. 

Section \ref{sec:big-tri} is devoted to the study of spherical bigons and triangles.
As our final goal is to prove an existence theorem, we do not try to state a full characterization of them (which can be found in \cite{troyanov:bigons} and \cite{eremenko:three}), but we only provide constructions. The last subsection deals with triangles which are close to a bigon or a union of two bigons: these will be the key ingredients for operating the surgery that splits a conical point.

Section \ref{sec:cut-paste} is rather short and presents three typical cases of surgery.
The first one takes place near a conical point and will be used to split
a conical singularity.
The second and the third one are performed along a path and will be used to produce spheres with angles $2\pi(\vc{\th}+\vc{e_i}+\vc{e_j})$ starting
from spheres with angles $2\pi\cdot\vc{\th}$.

Spheres with four conical points are constructed in Section \ref{sec:four-points}. Most of them can be obtained by doubling
a spherical (convex and non-convex) quadrilateral. Two sporadic one-parameter families of metrics 
requires an ad-hoc treatment.

Finally, in Section \ref{sec:induction} we show how to apply the previously developed tools to inductively construct all desired metrics
with $n\geq 5$ conical points and so to prove Theorem \ref{existence}.

\subsection{Acknowledgements}
This paper represents the first output of our joint effort
throughout a few years.
G.M. would like to thank Rafe Mazzeo for interesting conversations on this topic.\\
The research of G.M. was partially funded by
FIRB 2010 national grant {\it{``Low-dimensional geometry and topology''}} {{(code: RBFR10GHHH\_003)}}. D.P. is a Royal Society University Research Fellow.

\section{Algebraic constraints}

\subsection{Holonomy constraints}\label{sec:rep}

\subsubsection{Holonomy representation}

We recall the following well-known fact.

\begin{lemma}[Developing simply connected surfaces]\label{lemma:dev}
Let $\Omega$ be a connected surface endowed with a Riemannian metric of curvature $1$.
Then the following hold.
\begin{itemize}
\item[(i)]
$\Omega$ is locally isometric to $\Sph$.
\item[(ii)]
If $\Omega$ is simply connected, then these local isometries patch together to define a global {\it{developing map}} $\dev:\Omega\rar \Sph$, which is a local isometry.
\item[(iii)]
Let $\tilde{p}\in\Omega$ and $\tilde{v}\in T^1_{\tilde{p}} \Omega$ be a unit tangent vector.
If two developing maps $\dev,\dev':\Omega\rar \Sph$ agree on $(\tilde{p},\tilde{v})$, then they coincide on the whole $\Omega$.
\end{itemize}
\end{lemma}

Even if the surface is not simply connected, we can still develop
paths on spherical surfaces.

\begin{lemma}[Developing paths]\label{lemma:dev-path}
Let $\Sigma$ be a surface with a metric of curvature $1$ and let
$\gamma:[0,1]\rar\Sigma$ be a continuous path.
\begin{itemize}
\item[(i)]
There exists a simply connected surface $\Omega$ with a metric of curvature $1$ such that
the path $\gamma$ factorizes as $\gamma=j\circ\tilde{\gamma}$, where
$\tilde{\gamma}$ is a map $\tilde{\gamma}:[0,1]\rar \Omega$ and $j:\Omega\rar \Sigma$ is a local isometry onto its image.
\item[(ii)]
The composition of $\dev:\Omega\rar \Sph$ with $\tilde{\gamma}$
defines a developing map 
$\dev_\gamma=\dev\circ\tilde{\gamma}:[0,1]\rar \Sph$
for $\gamma$. Similarly, the composition of $d(\dev)$ and $d\tilde{\gamma}$
induce a $d(\dev_\gamma):\gamma^* T_{\Sigma}\rar  T_{\Sph}$.
%
\end{itemize}
\end{lemma}

The surface $\Omega$ should be thought of as a thickening of $\gamma$: for instance, if $\gamma$ is an embedding, we can choose $j$ to be the inclusion of a tubular neighbourhood of $\gamma([0,1])$.

\begin{proof}[Proof of Lemma \ref{lemma:dev-path}]
Take for instance $j:\Omega\rar\Sigma$ to be a universal cover and put
on $\Omega$ the pull-back metric from $\Sigma$. Then clearly $\gamma$ factorizes as desired and so (i) follows. Assertion (ii) is a consequence of Lemma \ref{lemma:dev}.
\end{proof}

In light of the previous lemma, the following is very natural.

\begin{definition}
Two paths $\gamma$ on $S$ and $\gamma'$ on $S'$ are {\it{isometric}} if their developing maps $\dev_{\gamma}$ and
$\dev_{\gamma'}$ agree up to an isometry of $\Sph$.
\end{definition}

Now we want to attach an element of $\SO(3,\RR)$ to every based loop on $\Sigma$.
Fix a basepoint $p\in \Sigma$ and a unit tangent vector $v\in T^1_p\Sigma$. Choose also a point $\bar{p}\in \Sph$ and a $\bar{v}\in T^1_{\bar{p}}\Sph$.

For every $\gamma$ loop on $\Sigma$ based at $p$,
let $\gamma=j\circ\tilde{\gamma}$ and $\Omega$ be as in Lemma \ref{lemma:dev-path}.
For $t=0,1$ there is a unique unit vector $\tilde{v}_t\in T^1_{\tilde{\gamma}(t)}\Omega$ such that $dj_{\tilde{\gamma}(t)}(\tilde{v}_t)=v$. 
Moreover, there exists a unique choice of $\dev:\Omega\rar \Sph$ that takes
$(\tilde{\gamma}(0),\tilde{v}_0)$ to $(\bar{p},\bar{v})\in T^1 \Sph$.

We will call $\rho(\gamma)$
the unique element of $\SO(3,\RR)$ that acts on $T^1 \Sph$ by taking
$d(\dev)_{\tilde{\gamma}(0)}(\tilde{v}_0)=(\bar{p},\bar{v})$ to $d(\dev)_{\tilde{\gamma}(1)}(\tilde{v}_1)$.

The conclusion is the following well-known fact.

\begin{corollary}[Holonomy representation]\label{cor:hol-SO(3,R)}
The association $\gamma\mapsto\rho(\gamma)$ descends to a well-defined 
homomorphism $\rho:\pi_1(\Sigma,p)\rar\SO(3,\RR)$, called
holonomy representation.
\end{corollary}

\begin{remark}
$\rho$ is unique up to {\it{global conjugation}}, namely
a different choice of $v$ and of $(\bar{p},\bar{v})$ produce the representation
$B\rho B^{-1}$ for some $B\in\SO(3,\RR)$.
\end{remark}

\begin{remark}\label{rmk:conj}
Given a free loop in $\Sigma$, the same construction
determines a conjugacy class of elements in $\SO(3,\RR)$.
\end{remark}

We will be particularly interested in the following application
of Corollary \ref{cor:hol-SO(3,R)}.

Let $S$ be a surface homeomorphic to a sphere. Let $x_1,\dots,x_n$ be distinct points of $S$ and let $\th_1,\dots,\th_n>0$. Denote by $\dot{S}$ the punctured surface $S\setminus\{x_1,\dots,x_n\}$.

\begin{corollary}[Holonomy representation for cone surfaces]\label{cor:SO(3,R)}
For every spherical metric $g$ on $\dot{S}$
with conical singularities of angles
$2\pi\th_i$ at $x_i$ and for every $p\in\dot{S}$, 
the induced holonomy representation
$\rho:\pi_1(\dot{S},p)\rar\SO(3,\RR)$ is well-defined up to global conjugation. 
Moreover, if $\gamma_j$ is a loop that simply winds around $x_j$,
then $\rho(\gamma_j)$ is a rotation of angle $2\pi\th_j$.
\end{corollary}

In order to perform cut-and-paste constructions, we will need to establish a certain technical deformability property of spherical metrics.


\begin{definition}\label{def:deformable} 
A spherical metric $g$ on $\dot{S}$ with conical singularities of angles
$2\pi\cdot\vc{\th}=(2\pi\th_1,\dots,2\pi\th_n)$ is {\it{angle-deformable}}
if there exists a {{neighbourhood $N$}} of $\vc{\th}\in\RR^n$ such that
the following property holds:
\begin{itemize}
\item There exists a continuous family of spherical metrics $(g_{\vc{\nu}})$ on $\dot{S}$ parametrised by $\vc{\nu} \in N$  with conical singularities of angles $2\pi \cdot\vc{\nu}$,
such that $g_{\vc{\th}}=g$.
\end{itemize}

\end{definition}

\begin{notation}
We say that the angle vector $\vc{\th}$ or the associated
defect vector $\vc{\delta}$ (see Definition \ref{def:defect}) is
{\it{deformable}} if there exists an angle-deformable metric with conical singularities
of angles $2\pi\cdot\vc{\th}$.
\end{notation}

A corollary of a theorem by Luo \cite{luo:monodromy} on $\CC\PP^1$-structures with moderate singularities can be specialized
to the case 
of spherical metrics with non-integral angles and non-coaxial holonomy.
Here we formulate it in a simplified form, well-suited to our needs.

\begin{theorem}[Deformability]\label{luo}
Let $\vc{\th}=(\th_1,\dots,\th_n)$ with each $\th_i$ positive and non-integral.
Suppose that there exists a spherical metric $g$ on $\dot{S}$
with conical singularities of angles $2\pi\th_i$ at $x_i$
and with non-coaxial holonomy. Then $g$ is
angle-deformable.
\end{theorem}

On the other hand, non-coaxiality of the holonomy
can be easily checked as follows.

\begin{lemma}[Non-coaxiality criterion]\label{noncoaxhol}
Consider a spherical metric $g$ on $\dot{S}$ with conical singularities
of angles $2\pi\th_i$ at $x_i$.
Suppose that there exists a smooth geodesic path $\gamma$
of length $\ell\notin \pi\ZZ$
between two distinct points $x_j$ and $x_k$ such that
$\th_j,\th_k\not\in\ZZ$.
Then the holonomy of $g$ is non-coaxial.

Moreover, if
$\ell$ is not a multiple of $\pi/2$ or if
$\th_j$ is not half-integral,
then the holonomy of $g$ is non-Abelian.
\end{lemma}
\begin{proof}
Fix a basepoint $p\in\dot{S}$ and let
$[\gamma_i]\in\pi_1(\dot{S},p)$ be the class of a loop
that simply winds about $x_i$. Let $\rho:\pi_1(\dot{S},p)\rar \SO(3,\RR)$ be the holonomy representation
associated to the metric $g$.
Since $\th_j$ and $\th_k$ are not integers, the transformations
$\rho(\gamma_j)$ and $\rho(\gamma_k)$
have well-defined axes $A_j,A_k\subset\RR^3$, which are in fact
distinct because $\ell$ is not a multiple of $\pi$. Hence, the holonomy
is not coaxial.

If (a) or (b) is satisfied, then $\rho(\gamma_j)$ does not fix $A_k$.
This shows that $\rho(\gamma_j)$ and $\rho(\gamma_k)$ do not commute and so the holonomy is not Abelian.
\end{proof}

\begin{remark}
It is easy to see that the only Abelian but non-coaxial holonomy
representation takes values (up to conjugation) in the non-cyclic
subgroup of order $4$ of diagonal matrices in $\SO(3,\RR)$.
Such a holonomy is indeed realized, for instance by a spherical surface
of genus $0$ with three conical points of angles $\pi$ and by suitable
branched covers of it.
\end{remark}

\subsubsection{Canonical lift to $\SU(2)$}

The statement of Corollary \ref{cor:SO(3,R)}
can be slightly improved as follows.

\begin{proposition}[Lift of the holonomy to $\SU(2)$]\label{prop:SU(2)}
Let $S$ be a surface homeomorphic to a sphere and
$x_1,\dots,x_n$ be distinct points of $S$ 
and let $p\in\dot{S}=S\setminus\{x_1,\dots,x_n\}$ be a basepoint.
Suppose that $\dot{S}$ is endowed with a Riemannian
metric of curvature $1$ with conical singularities of angles
$2\pi\th_j>0$ at $x_j$.
Then its holonomy representation
$\rho$ admits a canonical lift $\hat{\rho}:\pi_1(\dot{S},p)\rar\SU(2)$.
Moreover, if $\gamma_j$ is a loop that winds simply around $x_j$,
then $\hat{\rho}(\gamma_j)$ has eigenvalues $e^{\pm i\pi (\th_j-1)}$.
\end{proposition}

\begin{remark}
Analogously as before, $\hat{\rho}$ is well-defined up to
global conjugation. Moreover, a free loop in $\dot{S}$
determines a conjugacy class of elements in $\SU(2)$.
\end{remark}

\begin{remark}\label{rmk:SU(2)}
It is well-known that the $\PSL(2,\CC)$-valued holonomy representation
associated to any $\CC\PP^1$-structure on a compact Riemann surface $\Sigma$
can be lifted to $\SL(2,\CC)$ and that such lifts correspond to complex line bundles $L$ on $\Sigma$ such that $L^{\otimes 2}\cong T_{\Sigma}$ (see Lemma 1.3.1 in \cite{GKM}, for instance).

In our case, we are considering the punctured surface $\dot{S}$. 
Since $S$ has genus zero,
requiring that $\hat{\rho}(\gamma_j)$ has eigenvalues $e^{\pm i\pi(\th_j-1)}$ already guarantees the uniqueness of the lift, so that we only have to check the existence.
If we were dealing with a punctured surface $\dot{\Sigma}$ of positive genus, then the requirement on the eigenvalues of $\hat{\rho}(\gamma_j)$ would restore the correspondence between lifts of the holonomy representations and line bundles $L$ on $\Sigma$ such that
$L^{\otimes 2}\cong T_{\Sigma}$.
Since it is not needed here, we will not analyze this case.
%
%
%
\end{remark}

\begin{definition}
A {\it{standard set of matrices for $\vc{\th}\in\RR_+^n$}} is
a $n$-uple $(U_1,\dots,U_n)$ of elements of $\SU(2)$ such that
$U_1\cdot U_2\cdots U_n=I$ and the eigenvalues of $U_j$
are $e^{\pm i\pi(\th_j-1)}$ for $j=1,\dots,n$.
\end{definition}

An immediate consequence of Proposition \ref{prop:SU(2)} is the following.

\begin{corollary}[From metrics to standard matrices]\label{cor:matrices}
Let $\th_1,\dots,\th_n$ be positive real numbers.
Suppose that there exists a metric of curvature $1$ on a sphere $S$
with conical singularities of angles $2\pi\th_1,\dots,2\pi\th_n$.
Then there exists a standard set of matrices for $\vc{\th}$.
%
\end{corollary}

Fix a standard set of generators $\gamma_1,\dots,\gamma_n$ of
$\pi_1(\dot{S},p)$, namely
\begin{itemize}
\item
$\gamma_j:[0,1]\rar \dot{S}$ is a smooth simple loop that winds
counterclockwise around $x_j$;
\item
the images of $\gamma_j$ intersect only at $p$;
\item
$\gamma_1\ast\dots\ast\gamma_n\simeq c_p$ the constant path at $p$.
\end{itemize}

\begin{proof}[Proof of Corollary \ref{cor:matrices}]
Just let $U_j:=\hat{\rho}(\gamma_j)$, where $\hat{\rho}$ is
the canonical lift provided by Proposition \ref{prop:SU(2)}.
\end{proof}

Though Proposition \ref{prop:SU(2)} may be phrased
in the more general context of $\CC\PP^1$-structure with conical singularities, 
%
we wish to provide a complete proof
tailored to our needs in the setting of spherical metrics.

\begin{proof}[Proof of Proposition \ref{prop:SU(2)}]
We break the proof into three main steps.\\

{\it{Step 1: construction of the lift $\hat{\rho}$.}}\\
Since $S$ is homeomorphic to a sphere, there exists
an open disk $D\subset S$ that contains all
$x_1,\dots,x_n$ and the images of $\gamma_1,\dots,\gamma_n$.
We can for instance assume that $S\setminus D$ consists of a single point $q$.

Choose a nowhere zero smooth vector field $V$ on $D$ 
and let $\hat{V}:=\frac{V}{\|V\|}$ the induced unit vector field on $D\cap\dot{S}$.
Let $\eta$ be a path contained in a coordinate chart near $q$ and that simply winds around $q$.
Clearly, in such a coordinate chart $V|_{\eta}$ and so $\hat{V}|_{\eta}$ have winding number $\pm 2$.
Finally, choose a point $\bar{p}\in \Sph$ and a $\bar{v}\in T^1_{\bar{p}} \Sph$.

Represent an element in $\pi_1(\dot{S},p)$ as a path $\gamma:[0,1]\rar \dot{S}\cap D$. By Lemma \ref{lemma:dev-path}(ii) and Lemma \ref{lemma:dev}(iii), there exists a unique 
developing map $\dev_\gamma$
that takes $(p,\hat{V}(p))\in T^1_S$ to $(\bar{p},\bar{v})\in T^1 \Sph$. 
For every $t\in[0,1]$ let $R_\gamma(t)\in\SO(3,\RR)$ be
the unique transformation that takes
$d(\dev_\gamma)_0 (V(\gamma(0)))=(\bar{p},\bar{v})$ to $d(\dev_\gamma)_t(V(\gamma(t)))$.
%
%
%
The path $R_\gamma:[0,1]\rar\SO(3,\RR)$ is clearly continuous
and satisfies $R_\gamma(0)=I$ and $R_\gamma(1)=\rho(\gamma)$.

Let $\hat{R}_\gamma:[0,1]\rar \SU(2)$ be the unique continuous
lift of $R_\gamma$ via the standard double cover $\SU(2)\rar \SO(3,\RR)$
such that $\hat{R}_\gamma(0)=I$. Define $\hat{\rho}(\gamma):=\hat{R}_\gamma(1)$.

If $s\mapsto \gamma_s$ is a continuous family of loops in $\dot{S}\cap D$ based at $p$, then
$\rho(\gamma_s)=\rho(\gamma_0)$ and so
$\hat{\rho}(\gamma_s)=\hat{\rho}(\gamma_0)$ by continuity.
Thus, two loops based at $p$ that are homotopic in $\dot{S}\cap D=\dot{S}\setminus\{q\}$ have the same 
$\SU(2)$-holonomy: this defines a representation $\pi_1(\dot{S}\setminus\{q\},p)\rar \SU(2)$.

In order to see that the constructed $\SU(2)$-representation
descends to $\pi_1(\dot{S},p)\cong \pi_1(\dot{S}\setminus\{q\},p)/\langle \eta\rangle$,
it is enough to check that the $\SU(2)$-holonomy along $\eta$ is trivial.
As $\hat{V}|_{\eta}$ winds twice, $q$ is a smooth point for the metric of $S$
and $\eta$ is freely homotopic to $\gamma_1\ast\dots\ast\gamma_n$,
we obtain $\hat{\rho}(\gamma_1)\cdots\hat{\rho}(\gamma_n)=\hat{\rho}(\eta)=I$. 

Hence, we conclude that
$\hat{\rho}:\pi_1(\dot{S},p)\rar\SU(2)$ is a well-defined representation that lifts $\rho$.\\

{\it{Step 2: eigenvalues of $\hat{\rho}(\gamma_j)$.}}\\
Let $p'$ be a point very close to $x_j$ and let $\beta$
be a loop based at $p'$ that keeps at constant distance from $x_j$
and that simply winds around $x_j$ at constant speed.
Clearly, the path $\gamma_j$ is homotopic to $\gamma'_j=\alpha^{-1}\ast\beta\ast\alpha$, where
$\alpha$ is a suitable simple path from $p$ to a point $p'$.

Thus, $R_{\gamma'_j}(1)$ can be written as
$A^{-1}BA$, where $A=R_{\gamma'_j}(\frac{1}{3})$ and
$B=R_{\gamma'_j}(\frac{2}{3})R_{\gamma'_j}(\frac{1}{3})^{-1}$,
and so $\hat{R}_{\gamma'_j}(1)=\hat{A}^{-1}\hat{B}\hat{A}$
for $\hat{A}=\hat{R}_{\gamma'_j}(\frac{1}{3})$ and
$\hat{B}=\hat{R}_{\gamma'_j}(\frac{2}{3})\hat{R}_{\gamma'_j}(\frac{1}{3})^{-1}$.
%
%
%
%
%
By our choice of $\beta$, the path $[0,1]\ni t\mapsto R_{\gamma'_j}(\frac{1+t}{3})R_{\gamma'_j}(\frac{1}{3})^{-1}$ is very close to be a rotation about a fixed axis of constant speed $2\pi\th_j$ and $B$ is a rotation of angle $2\pi\th_j$.
As a consequence, $[0,1]\ni t\mapsto\hat{R}_{\gamma'_j}(\frac{1+t}{3})\hat{R}_{\gamma'_j}(\frac{1}{3})^{-1}$ is very close
to be conjugate to the path $[0,1]\ni t\mapsto \mathrm{diag}(e^{it\pi(\th_j-1)},e^{-it\pi(\th_j-1)})$. Thus, $\hat{B}$
has eigenvalues $e^{\pm i\pi(\th_j-1)}$ and so the same holds
for $\hat{\rho}(\gamma_j)$.\\
%
%
%

{\it{Step 3: the lift $\hat{\rho}$ is canonical.}}\\
Consider first another nowhere vanishing vector field $W$ on $D$
and let $\hat{W}=\frac{W}{\|W\|}$.
There exists a continuous function 
$\bar{a}:D\rar \RR/\ZZ$  such that
$\hat{W}(x)$ is obtained from $\hat{V}(x)$ by a counterclockwise rotation of an angle $2\pi \bar{a}(x)$.
Since $D$ is simply connected, the function $\bar{a}$ lifts to
a continuous $a:D\rar\RR$. If we call $\hat{V}_s(x)$ the vector obtained by rotating $\hat{V}(x)$ counteclockwise by $s\cdot 2\pi a(x)$ for all $s\in[0,1]$, then $s\mapsto \hat{V}_s$ is a continuous family of unit vector fields on $D$ with $\hat{V}_0=\hat{V}$ and $\hat{V}_1=\hat{W}$. This induces a continuous family $s\mapsto\hat{\rho}_s$ of lifts of $\rho$, which must thus be the constant family. Hence, $\hat{\rho}_0=\hat{\rho}_1$ and so $\hat{\rho}$ does not depend on the choice of the vector field.

Finally, if $D'=S\setminus\{q'\}$ is another disk, then there is an isotopy that moves $q$ to $q'$ fixing $\{x_1,\dots,x_n\}$ and so
it moves $\dot{S}\cap D$ to $\dot{S}\cap D'$. Again, this determines a continuous family of lifts of $\rho$, which must then be constantly equal to $\hat{\rho}$.
\end{proof}

We remark that the coaxiality condition for the $\SO(3,\RR)$-valued holonomy representation can be rephrased in a more familiar way in terms of its lift.

\begin{lemma}[Non-coaxial subgroups]\label{lemma:non-coaxial}
Let $\hat{G}$ be a subgroup of $\SU(2)$ and let
$G$ be its image via the natural projection $\SU(2)\rar\SO(3,\RR)$.
\begin{itemize}
\item[(a)]
The group $\hat{G}$ is commutative if and only if $G$ belongs to a $1$-parameter subgroup of $\SO(3,\RR)$.
Hence, the canonical lift $\hat{\rho}$ is Abelian $\iff$ the representation $\rho$ is coaxial.
\item[(b)]
If $G$ is non-coaxial and $\tau\in\PSL(2,\CC)$ commutes with all elements in $G$, then $\tau\in\SO(3,\RR)$.
Hence, if $\tau\rho\tau^{-1}=\rho$ and $\rho$ is non-coaxial, then
$\tau\in\SO(3,\RR)$.
\end{itemize}
%
\end{lemma}
\begin{proof}
Since unitary matrices are diagonalizable, then $\hat{G}$ is Abelian
if and only if all matrices in $\hat{G}$ are simultaneously diagonalizable. This occurs if and only if $\hat{G}$ is contained in a $1$-parameter subgroup of $\SU(2)$, which is equivalent to
asking that $G$ is contained in a $1$-parameter subgroup of $\SO(3,\RR)$. This proves (a).

As for (b), let
$\hat{\tau}\in\mathrm{SL}(2,\CC)$ be a lift of $\tau$, so that
%
$\hat{\tau}\hat{\gamma}=\pm\hat{\gamma}\hat{\tau}$ for every $\hat{\gamma}\in\hat{G}$.
Up to conjugation by a matrix in $\SU(2)$, we can assume that
\[
\hat{\tau}=\left(
\begin{array}{cc}
\lambda & z\\
0 & \lambda^{-1}
\end{array}
\right)
\]
with $\lambda,z\in\CC$ and $|\lambda|\geq 1$, so that
\[
h:=\ol{\hat{\tau}}^T\hat{\tau}=
\left(
\begin{array}{cc}
|\lambda|^2+|z|^2 & z\ol{\lambda}^{-1}\\
\ol{z}\lambda^{-1} & |\lambda|^{-2}
\end{array}
\right)
\]
has $\mathrm{det}(h)=1$ and $t:=\frac{1}{2}\tr(h)=\frac{1}{2}\big(|\lambda|^2+|\lambda|^{-2}+|z|^2\big)\geq 1$. Thus, $h$ is diagonalizable and it has eigenvalues $\mu^2$ and $\mu^{-2}$, with $\mu=\sqrt{t+\sqrt{t^2-1}}\geq 1$.
It is easy to see that $\|\tau v\|\leq \mu \|v\|$ for every
$v\in\CC^2$ and equality holds
if and only if $v$ belongs to the $\mu^2$-eigenspace $E_{\mu^2}\subseteq\CC^2$ of $h$.
Since $\|\tau(\hat{\gamma}(v))\|=\|\hat{\gamma}(\tau (v))\|=
\|\tau (v)\|\leq \mu\|v\|=\mu\|\hat{\gamma}(v)\|$,
every $\hat{\gamma}\in\hat{G}$ preserves $E_{\mu^2}$.

By (a), the group $\hat{G}$ is not Abelian and so $E_{\mu^2}$ cannot be $1$-dimensional. This implies that $t=1$
and so $|\lambda|=1$ and $z=0$, which shows that $\hat{\tau}\in\SU(2)$ and finally $\tau\in\SO(3,\RR)$.
\end{proof}

\subsubsection{Matrices in $\SU(2)$ and broken geodesics on $\Spt$}

In view of Corollary \ref{cor:matrices}, it is natural first
to discuss the following.

\begin{problem}\label{problem:matrices}
Find a criterion for the existence of
a standard set of matrices $U_1,\dots,U_n\in\SU(2)$
for $\vc{\th}\in\RR^n_+$ that do not simultaneusly
commute.
\end{problem}

This problem was addressed in many papers and explicit
inequalities are known (see \cite{biswas:parabolic}).
In order to motivate these inequalities,
we recall how this question is equivalent to a different
question about broken geodesics on the standard
$3$-sphere $\Spt$.
%

\begin{notation}
By {\it{broken geodesic}} on $\Spt$ we will mean
a piecewise geodesic path with endpoints $v_0$ and $v_n$ that passes through an ordered collection of points $v_0,\dots,v_n$ of $\Spt$ in such a way that each {\it{side}} $s_j$ going from
the {\it{vertex}} $v_{j-1}$ to the vertex $v_j$ is of minimal length (and so at most $\pi$).
\end{notation}

Given a broken geodesic on $\Spt$ with vertices $v_0,\dots,v_n$,
we define $U_j\in \SU(2)$ as the unique transformation that
takes $v_{j-1}$ to $v_j$ for $j=1,\dots,n$.

Vice versa, given matrices $U_1,\dots,U_n$ in $\SU(2)$
and fixed a basepoint $v_0:=(1,0)$ on the unit sphere $\Spt\subset\CC^2$, 
we define $v_j:=U_j(v_{j-1})=U_j U_{j-1}\cdots U_1(v_0)$ for $j=1,\dots,n$.
A broken geodesic $\Gamma$ is then obtained by drawing a {\it{shortest}} geodesic $s_j$ from $v_{j-1}$ to $v_{j}$ for all $j=1,\dots,n$. Notice that, given $v_{j-1}$, the segment $s_j$ is uniquely determined unless $U_j=-I$.


%
%

Clearly, the matrices $U_j$ satisfy $U_1\cdots U_n=I$ if and only
if $v_n=v_0$, i.e. if and only if the broken geodesic is {\it{closed}}.
%


\begin{definition}\label{def:defect}
Let $\vc{\th}\in\RR^n$ be an {\it{angle vector}}. Its associated {\it{defect vector}}
is $\vc{\delta}:=\vc{\th}-\vc{\1}\in\RR^n$, where $\vc{\1}=(1,1,\dots,1)$.
%
The associated {\it{reduced angle vector}} $\vc{\ol{\th}}\in\RR^n$ is defined in such a way that $\ol{\th}_j\in[0,2)$ 
and $\th_j-\ol{\th}_j\in 2\ZZ$. Finally, the {\it{reduced defect vector}} is $\ol{\vc{\delta}}:=\ol{\vc{\th}}-\vc{\1}\in [-1,1)^n$.
\end{definition}

\begin{remark}
The definition of $\ol{\vc{\th}}$ is motivated by the fact that
the edge $s_j$ of the broken geodesic on $\Spt$ associated to $U_1,\dots,U_n$ has length $\ell_j=\pi |1-\ol{\th}_j|=\pi|\ol{\delta}_j|$.
\end{remark}

We summarize the content of the above discussion into the following.

\begin{lemma}[Broken geodesics and standard set of matrices]
Let $\th_1,\dots,\th_n>0$.
Then the following are equivalent:
\begin{itemize}
\item[(a)]
there exists a closed broken geodesic on $\Spt$ with $n$ edges 
of length $\ell_j=\pi|\ol{\delta}_j|$ for $j=1,\dots,n$;
\item[(b)]
there exists a standard sets of matrices $U_1,\dots,U_n\in\SU(2)$
for $\vc{\th}$.
\end{itemize}
\end{lemma}

Notice that, through the identification $\SU(2)\ni U\mapsto
U(v_0)\in \Spt$,
there is a correspondence between
$1$-parameters subgroups of $\SU(2)$ and
maximal circles on $\Spt$ through $v_0$.
Thus, the matrices $U_1,\dots,U_n$ simultaneously commute
if and only if $v_0,\dots,v_n$ all belong to the same maximal circle.

The following result
is essentially contained in \cite{biswas:parabolic}.

\begin{theorem}[Constraints for broken geodesics]\label{thm:polygon}
There exists a closed broken geodesic on $\Spt$ with $n$ edges of length
$\ell_1,\dots,\ell_n\in [0,\pi]$ if and only if
\begin{equation}\tag{Pol}\label{eq:polygon}
\sum_{j\in X}(\pi-\ell_j)+\sum_{k\in X^c}\ell_k\geq \pi 
\end{equation}
for all $X\subseteq\{1,\dots,n\}$ with $|X|$ odd.
\end{theorem}

\begin{remark}
These inequalities are  generalizations of the
following simple statement:   there cannot be a closed broken
geodesic on $\Spt$
with odd number of edges of length $\pi$. Moreover,
even if we replace each length  $\ell_j=\pi$ by $\ell_j=\pi\pm \varepsilon_j$
with $\sum_j |\varepsilon_j|< \pi$, then
a closed broken geodesic cannot exist.
\end{remark}

\begin{lemma}[Broken geodesics on a maximal circle]\label{lemma:equality}
Equality in \eqref{eq:polygon} has the following geometric counterpart.
\begin{itemize}
\item[(i)]
If equality is attained in \eqref{eq:polygon} for a certain
$X$, then every closed broken geodesic on $\Spt$ with edges of lengths $\ell_j$ sits on a maximal circle.
\item[(ii)]
If a closed broken geodesic on $\Spt$ with edges of lengths $\ell_j$ sits on a maximal circle, then
\[
\sum_{j\in Y}\ell_j -\sum_{k\in Y^c}\ell_k \equiv 0 \quad
\text{(mod $2\pi$)}
\]
for some subset $Y\subseteq\{1,\dots,n\}$.
\end{itemize}
\end{lemma}
\begin{proof}
Indeed, (ii) is immediate: once fixed an orientation on the maximal circle, just let $Y$ be the collection of positively oriented edges of the broken geodesic.
About (i), it is enough to notice that, if a broken geodesic $\Gamma$ does not sit on a maximal circle, then it can be deformed in such a way that the quantity on the left-hand side of
\eqref{eq:polygon} decreases. 
\end{proof}

The above characterization of
the lengths of the edges of 
closed broken geodesics on $\Spt$
gives a criterion for the existence of
matrices 
that satisfy
the conditions of Problem \ref{problem:matrices}.
Now we want to rewrite such criterion in a more compact way.

\begin{corollary}[Angle constraints for representations in $\SU(2)$]\label{cor:holonomy}
Given $\th_1,\dots,\th_n>0$, the following facts are equivalent.
\begin{itemize}
\item[(1)]
There exists a standard set of matrices $U_1,\dots,U_n\in\SU(2)$ for $\vc{\th}$.
\item[(2)]
The following inequalities hold
\begin{equation}\tag{Pol'}\label{eq:hol1}
\sum_{j\in X}\left(1-\left|\ol{\delta}_j\right|\right)
+\sum_{k\in X^c}\left|\ol{\delta}_k\right|\geq 1
\end{equation}
for all $X\subseteq\{1,\dots,n\}$ with $|X|$ odd.
\item[(3)]
The following inequality holds
\begin{equation}\tag{H}\label{eq:holonomy}
d_1(\vc{\delta},\ZZ^n_o)\geq 1
\end{equation}
\end{itemize}
Also, equality holds in (3) if and only if it holds in (2) for a certain $X$.
Moreover:
\begin{itemize}
\item[(I)]
if $d_1(\vc{\delta},\ZZ^n_o)=1$, then 
all standard $n$-uples of matrices $U_1,\dots,U_n$
for $\vc{\th}$
%
simultaneously commute and in fact they belong to the same
$1$-parameter subgroup of $\SU(2)$;
\item[(II)]
if there exists
a standard $n$-uple of matrices
$U_1,\dots,U_n$ that belong to
the same $1$-parameter subgroup of $\SU(2)$,
then 
\[
\sum_{j\in Y}\th_j -\sum_{k\in Y^c}\th_k\equiv 0
\quad\text{(mod $2$)}
\]
for a certain subset $Y\subseteq \{1,\dots,n\}$.
\end{itemize}
\end{corollary}
\begin{proof}
The equivalence (1)$\iff$(2) is just a rephrasing of Theorem \ref{thm:polygon}. Moreover, it is clear that (I) and (II)
are rephrasings of (i) and (ii) in Lemma \ref{lemma:equality}.
So it is enough to show that (2)$\iff$(3), which is a consequence 
of the following equality
\begin{equation}\tag{H=Pol'}\label{eq:HP'}
d_1(\vc{\delta},\ZZ^n_o)=
\inf_{\text{$|X|$ odd}}
\left(
\sum_{j\in X}\left(1-\left|\ol{\delta}_j\right|\right)
+\sum_{k\in X^c}\left|\ol{\delta}_k\right|\right)
\end{equation}

Let $\vc{m}\in\ZZ^n_o$ and call $X$ the subset of indices in $\{1,\dots,n\}$ for which
$m_i$ is odd. Clearly, $|X|$ is odd because $\|\vc{m}\|_1$ is.
It is easy to see that  $|\delta_j-m_j|
\geq 1-\left|\ol{\delta}_j\right|$
for $j \in X$,
and $|\delta_k-m_k|\geq 
\left|\ol{\delta}_k\right|$ for $k\in X^c$.
Thus, $d_1(\vc{\delta},\vc{m})\geq \sum_{j\in X}\left(1-\left|\ol{\delta}_j\right|\right)
+\sum_{k\in X^c}\left|\ol{\delta}_k\right|$.
Moreover, the equality is attained for those $\vc{m}\in\ZZ^n_o$ for which $|\delta_j-m_j|\leq 1$ for all $j=1,\dots,n$.
Thus, equation \eqref{eq:HP'} holds.
\end{proof}

As a consequence, we can determine a necessary condition
for the existence of a metric of curvature $1$ on $\Sph$
with cone points of angles $\th_1,\dots,\th_n$.

\begin{proof}[Proof of Theorem \ref{nonexist}]
It follows by combining Corollary \ref{cor:matrices} and
implication
(1)$\implies$(3) in
Corollary \ref{cor:holonomy}.
\end{proof}

\subsection{Algebraic merging}\label{sec:algebraic-merging}

The main goal of this section is to prove Theorem 
\ref{linearalgebra}.
In order to do that, we first need to set some notation.

\medskip

Given $\vc{\th}=(\th_1,\dots,\th_n)\in\RR_+^n$,
we recall that the defect vector is $\vc{\delta}=\vc{\th}-\vc{\1}$. 
Throughout this section it will turn more practical to directly work with $\vc{\delta}$ instead of $\vc{\th}$.

\begin{notation}
Denote by $\Hang^n$ the subset of $\vc{\delta}\in\RR^n$ 
that satisfy the {\it{holonomy constraints}}, namely such that
$d_1(\vc{\delta},\ZZ^n_o)\geq 1$. This is the complement in $\RR^n$ of a union of octahedrons centred at points 
of $\ZZ^n_o$.
Denote by $\Pang^n$ the subset of $\vc{\delta}\in\RR^n$ that
satisfy the {\it{positivity constraints}}, namely such that
$\delta_1+\dots+\delta_n>-2$ and $\delta_1,\dots,\delta_n>-1$.
The locus $\Pang^n\cap \Hang^n$ of {\it{admissible}} defect vectors will be denoted by $\Aang^n$.
\end{notation}

Let $n\geq 4$ and let $i,j\in\{1,\dots,n\}$ be distinct.
We define the {\it{positive/negative algebraic merging operation}}
\[
M_{(i\pm j)}:
\RR^n\lra \RR^{n-1}
\]
as $M_{(i+j)}(\delta_1,\dots,\delta_n):=
(\delta_1,\dots,\widehat{\delta}_i,\dots,\widehat{\delta}_j,\dots,\delta_n,
\delta_i+\delta_j)$
and $M_{(i-j)}(\delta_1,\dots,\delta_n):=
(\delta_1,\dots,\widehat{\delta}_i,\dots,\widehat{\delta}_j,\dots,\delta_n,
\delta_i-\delta_j-2)$.

\begin{lemma}[Basic properties of $M_{(i\pm j)}$]\label{lemma:alg-merging-properties}
Every algebraic merging operation $M:\RR^n\rar\RR^{n-1}$
satisfies the following properties:
\begin{itemize}
\item[(a)]
$M(\ZZ^n_o)=\ZZ^{n-1}_o$;
\item[(b)]
$M$ is contracting for the $\ell^1$ metrics;
\item[(c)]
$M(\vc{\delta})\in \inte(\Hang^{n-1})
\implies
\vc{\delta}\in\inte(\Hang^n)$.
\end{itemize}
\end{lemma}
\begin{proof}
Property (a) is obvious and claim (c) follows from (a) and (b).
As for (b), up to reordering the coordinates we can assume that 
$M=M_{(1+2)}$ or $M=M_{(1-2)}$.

Let $\vc{m},\vc{m'}\in\RR^n$. Then
\begin{align*}
d_1(M_{(1+2)}(\vc{m}),M_{(1+2)}(\vc{m'}))
& = |(m_1+m_2)-(m'_1+m'_2)|+\sum_{j=3}^n |m_j-m'_j| \leq\\
& \leq |m_1-m'_1|+|m_2-m'_2|+\sum_{j=3}^n |m_j-m'_j| =
 d_1(\vc{m},\vc{m'})\, .
\end{align*}
The proof for $M=M_{(1-2)}$ is completely analogous.
\end{proof}

The main result of this section is the following more precise version of
Theorem \ref{linearalgebra}.

\begin{theorem}[Algebraic merging]\label{5points}
Let $n\geq 5$ and
suppose that $\vc{\delta}\in\inte(\Aang^n)$.
Then there exist distinct indices $i,j\in\{1,\dots,n\}$
such that at least one of the following holds:
\begin{itemize}
\item[(a)]
$M_{(i+j)}(\vc{\delta})\in\inte(\Aang^{n-1})$;
\item[(b)]
$M_{(i-j)}(\vc{\delta})\in\inte(\Aang^{n-1})$
and $\delta_i,\delta_j,\delta_i-\delta_j\notin\ZZ$.
\end{itemize}
\end{theorem}

In order to prove the above result, we will
separately analyze three different cases.
We will see that in most situations it is possible
to find indices $i,j$ such that (a) holds.
%


\subsubsection{Intersection of $\Aang^n$ with a unit integer cube}\label{sec:analysis}

Throughout this section, we assume $n\geq 3$.

\begin{notation}
We will use the symbol
$\cube^n$ to denote any closed unit cube with integer vertices in
$\RR^n$ and 
the symbol $\tcube^n$ to denote the truncated cube
obtained by intersecting $\cube^n$ with $\Hang^n$.
Sometimes, we will use the notation $\cube_{\vc{c}}$
to indicate the unit cube with center in $\vc{c}=(c_1,\dots,c_n)\in\RR^n$.
\end{notation}


\begin{center}
\begin{figurehere}
\psfrag{odd}{\parbox[t]{3cm}{Odd-integral\\vertices of the cube}}
\includegraphics[width=0.3\textwidth]{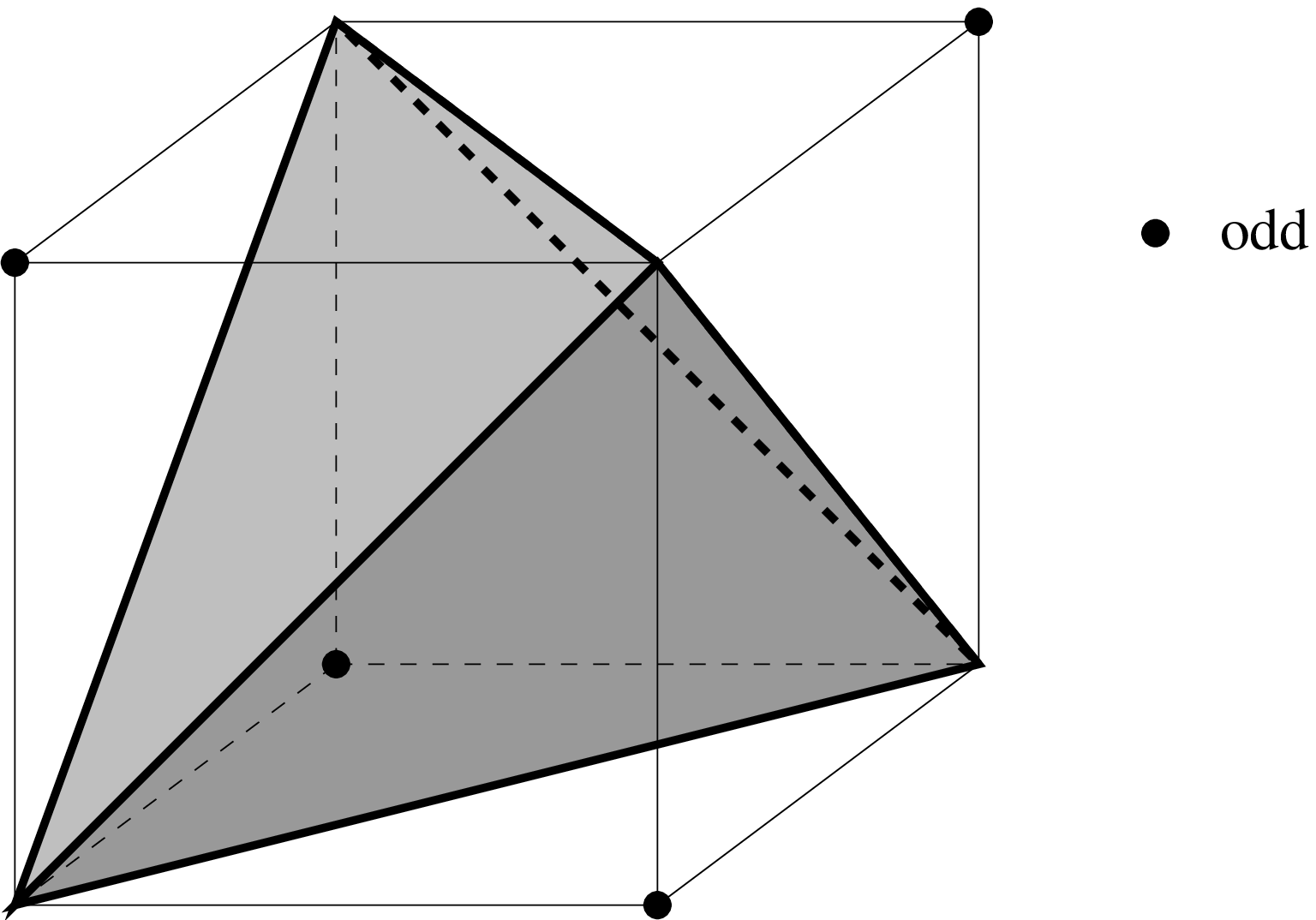}
\caption{{\small A $3$-dimensional truncated cube $\tcube^3$.}}\label{fig:tcube}
\end{figurehere}
\end{center}

\begin{lemma}[Truncated cubes]
Let $\cube^n$ be  a unit cube in $\RR^n$ with integer vertices. The intersection $\tcube^n=\cube^n\cap \Hang^n$ is the 
convex hull of all even vertices of $\cube^n$. 
\end{lemma}

\begin{proof}
Note that $\cube^n\cap \Hang^n$ consist of all points of $\cube^n$ that 
are on $\ell^1$ distance at least one from all odd vertices of $\cube^n$. If $\vc{m}$ is 
such an odd vertex, then the points in $\cube^n$ at distance at most $1$ from 
$\vc{m}$ are formed by the simplex spanned by $\vc{m}$ and the $n$ even vertices
of $\cube^n$ at distance $1$ from $\vc{m}$. Hence, the set  $\cube^n\cap \Hang^n$ is 
obtained from $\cube^n$ by cutting away $2^{n-1}$ simplices 
corresponding to the odd vertices of $\cube^n$. 
It follows that $\tcube^n$ is a convex
polytope and it is easy to see that its vertices are
the even vertices of $\cube^n$.
\end{proof}

\begin{remark} 
Since $\Pang^n$ is convex and $\Aang^n=\Pang^n\cap \Hang^n$ it follows from 
the lemma that the intersection of any integral unit cube $\cube^n$ with $\Aang^n$ is convex.
\end{remark}

As a consequence, we deduce the connectedness of $\Hang^n$ and of $\Aang^n$ for $n\geq 4$.

\begin{proof}[Proof of Lemma \ref{lemma:connected}]
Each $n$-dimensional simplex corresponding to an even vertex of $\cube^n$ has volume $1/n!$ and so
$\tcube^n$ has volume $1-\frac{2^{n-1}}{n!}>0$, because $n\geq 3$. Hence, all $\tcube^n$ have non-empty interior and so, in particular, $\Hang^n$ and $\Aang^n$ are non-empty.

Let now $n\geq 4$. We claim that, if $\vc{c}=(c_1,\dots,c_n)$ is the center of a unit cube $\cube_{\vc{c}}$ and $\vc{e_i}$ is a vector in the standard basis of $\RR^n$, then the interior of $\tcube_{\vc{c}}\cup\tcube_{\vc{c}+\vc{e_i}}$ is connected.
In fact, the two adjacent truncated cubes $\tcube_{\vc{c}}$ and $\tcube_{\vc{c}+\vc{e_i}}$ share a face $\Fang$ isometric to a lower dimensional truncated cube $\tcube_{\vc{c'}}$, where
$\vc{c'}=(c_1,\dots,\widehat{c}_i,\dots,c_n)$.
As $n-1\geq 3$, the interior $\inte(\Fang)\cong \inte(\tcube_\vc{c'})$ (as a subset of $\RR^{n-1}$) is non-empty: let $\vc{\delta}$ be a point in $\inte(\Fang)$.
As $\inte(\tcube_{\vc{c}}\cup\tcube_{\vc{c}+\vc{e_i}})$ 
is star-shaped with respect to $\vc{\delta}$, the claim follows.

One then easily concludes that $\inte(\Hang^n)$ and $\inte(\Aang^n)$ are connected for $n\geq 4$.
\end{proof}

Note that the boundary $\partial\tcube^n$ of a truncated $n$-cube is made of $2^{n-1}$ faces isometric to $(n-1)$-simplices (one for each odd vertex of $\cube^n$) and $2n$ faces isometric to truncated $(n-1)$-cubes (one for each face of $\cube^n$). 
All faces have an interior part for $n\geq 4$ (see Figure \ref{fig:tcube2}), whereas the non-simplicial faces of $\partial\tcube^3$ are degenerate, as it appears clearly
in Figure \ref{fig:tcube}.

\begin{center}
\begin{figurehere}
\psfrag{simplicialface}{Simplicial face}
\psfrag{nonsimplicialface}{Non-simplicial face}
\psfrag{p}{$\vc{\delta}$}
\psfrag{c}{$\vc{c}$}
\psfrag{dpp}{$\vc{\delta^\pi}$}
\includegraphics[width=0.4\textwidth]{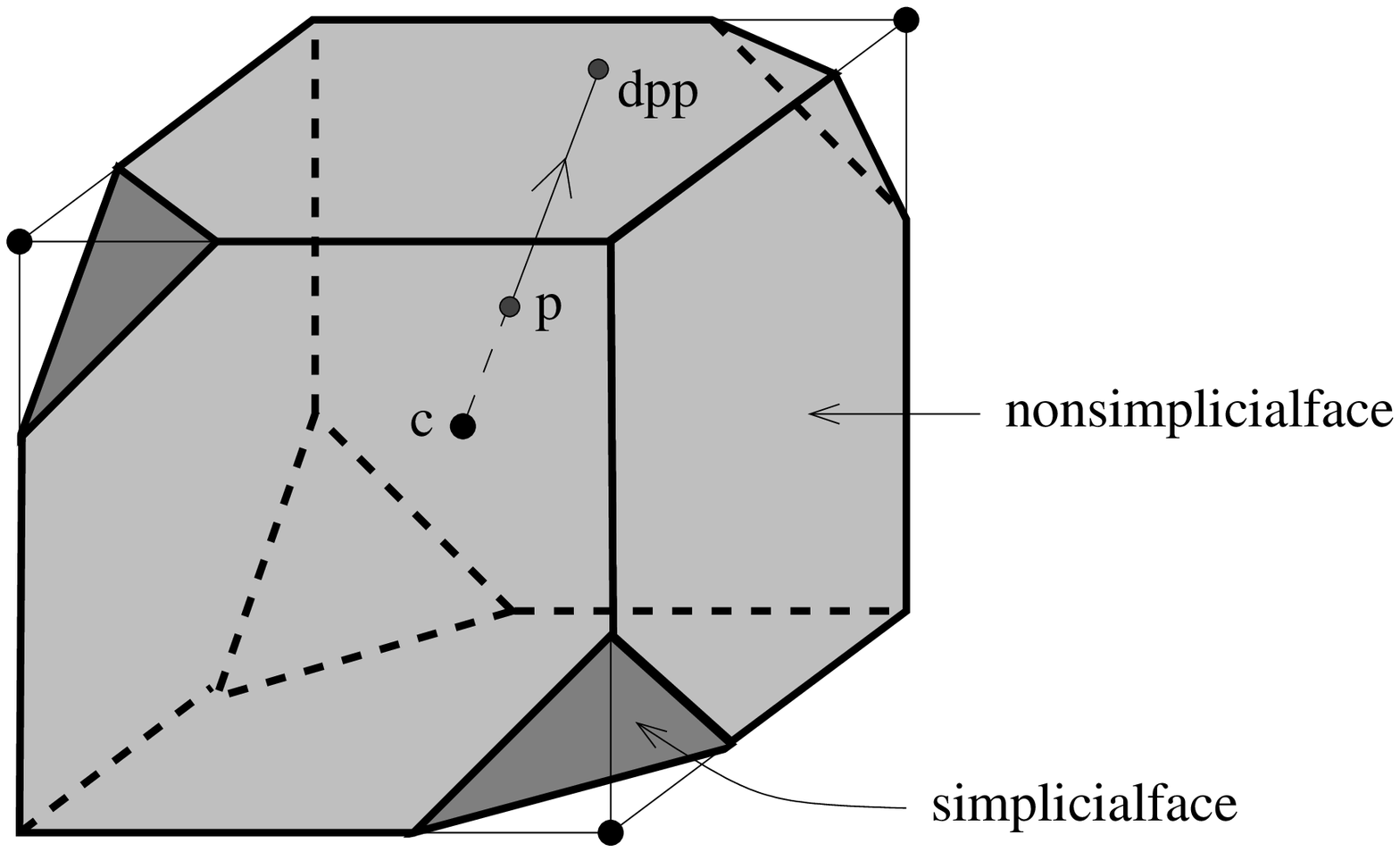}
\caption{{\small Symbolic picture of a truncated cube $\tcube^n$ with $n\geq 4$.}}\label{fig:tcube2}
\end{figurehere}
\end{center}

\begin{notation}
Denote by $\vc{c}$ the center of $\cube^n$, which is a point with half-integral coefficients. For any $\vc{\delta}\in \tcube^n$ different from $\vc{c}$ denote by $\vc{\delta}^{\pi}$ the
projection of $\vc{\delta}$ to the boundary of $\tcube^n$ from the center $\vc{c}$, i.e. the unique
point on $\partial\tcube^n$ such that 
$\vc{\delta}$ belongs to the segment that joins $\vc{\delta}^{\pi}$ and $\vc{c}$.
\end{notation}

In what follows we will distinguish two types of points in $\tcube^n$.

\begin{definition} 
A point $\vc{\delta}\in \tcube^n$ different from its center is called {\it simplicial} if 
$\vc{\delta}^{\pi}$ belongs to a simplicial face of $\partial\tcube^n$;
otherwise, $\vc{\delta}$ is called
{\it non-simplicial}.
\end{definition}

The following lemma summarizes some simple useful properties of simplicial and
non-simplicial points.

\begin{lemma}[Boundary of truncated cubes]\label{simpleproperties} 
Let $\vc{\delta}$ be a point in a truncated cube $\tcube^n$ different from its center $\vc{c}$.
Then the following hold.
\begin{itemize}
\item[(a)]
If $\vc{\delta}$ is non-simplicial, then 
$d_1(\vc{\delta}^{\pi},\ZZ_o^n)>1$ and
there exists an $i$ such that $\delta^{\pi}_i\in \ZZ$.
\item[(b)]
$\vc{\delta}$ is simplicial if and only if $d_1(\vc{\delta}^{\pi},\ZZ_o^n)=1$.
\item[(c)]
Suppose that $\vc{\delta}$ is a simplicial point and let $\vc{m}$ be a corresponding odd-integral vertex of $\cube^n$.
Then $\vc{m}$ is a closest point to $\vc{\delta}$ among 
all odd-integral points.
\end{itemize} 
\end{lemma}

\begin{proof}
The first two statements directly follow from the definitions, so we only prove (c). Let $\vc{m^{(1)}},\dots,\vc{m^{(n)}}$ be the even vertices of $\cube^n$ sitting at $\ell^1$ distance $1$
from $\vc{m}$.
Since $\vc{\delta}$ is simplicial it lies in the convex hull $\Kang$ of 
points $\vc{c},\vc{m^{(1)}},\dots,\vc{m^{(n)}}$.

Let $\vc{m'}$ be any other point in $\ZZ^n_o$ and denote $\vc{\delta'}$ a point in $\Kang$. Note that both functions 
$f=d_1(\bullet,\vc{m})$ and $f'=d_1(\bullet,\vc{m'})$ are affine on $\Kang$. So in order to prove that $d_1(\vc{\delta},\vc{m})\leq d_1(\vc{\delta},\vc{m'})$ it is enough to show that $f(\vc{p})\leq f'(\vc{p})$ for every vertex $\vc{p}$ of $\Kang$.
%
This is indeed so, because $d_1(\vc{m'},\vc{c})\ge d_1(\vc{m},\vc{c})=\frac{n}{2}$ and $d_1(\vc{m'},\vc{m^{(i)}})\ge d_1(\vc{m},\vc{m^{(i)}})=1$.
\end{proof}


\subsubsection{Simplicial and non-simplicial merging}

Even if we begin with an admissible defect vector $\vc{\delta}$,
the output of a positive algebraic merging operation might be no longer admissible.
As positivity issues are generally easier to keep under control,
here we focus on the problem of determining 
whether $M_{(i+j)}(\vc{\delta})$ belongs to $\inte(\Hang^{n-1})$
for given $i\neq j$ and $\vc{\delta}\in\inte(\Hang^n)$.

First, a simple observation about merging integral angles.

\begin{lemma}[Integral merging]\label{lemma:int-merging}
Let $\vc{\delta}\in \RR^n$ be a vector such that $\delta_i\in\ZZ$ for some $i$.
Let $M$ be a merging operation of type $M_{(i+j)},M_{(i-j)},M_{(j-i)}$
for some $j\neq i$.
Then $d_1(M(\vc{\delta}),\ZZ_o^{n-1})=d_1(\vc{\delta},\ZZ^n_o)$.
Hence, $M(\vc{\delta})\in \inte(\Hang^{n-1})$ if and only if $\vc{\delta}\in \inte(\Hang^{n})$.
%
\end{lemma}
\begin{proof}
Since $M(\ZZ^n_o)\subseteq \ZZ^{n-1}_o$ and $M$ is contracting for the $\ell^1$ distances, we have $d_1(M(\vc{\delta}),\ZZ^{n-1}_o)\leq d_1(\vc{\delta},\ZZ^n_o)$.
It is then enough to show that $d_1(\vc{\delta},\ZZ^n_o)\leq d_1(M(\vc{\delta}),\ZZ^{n-1}_o)$.
We will prove it for $M=M_{(i+j)}$, the other cases being analogous.
Moreover, up to reordering the coordinates, we can assume that $i=1$ and $j=2$.

For every $\vc{m}\in\ZZ^{n-1}_o$, we define $\vc{m'}:=(\delta_1,m_{n-1}-\delta_1,m_1,m_2,\dots,m_{n-2})
\in\ZZ^n_o$
so that $M_{(1+2)}(\vc{m'})=\vc{m}$.
Now, $d_1(\vc{\delta},\vc{m'})=|\delta_1-\delta_1|+|\delta_2-(m_{n-1}-\delta_1)|+\sum_{j=3}^n|\delta_j-m_{j-2}|=
d_1(M_{(1+2)}(\vc{\delta}),\vc{m})$,
since $M_{(1+2)}(\vc{\delta})=(\delta_3,\dots,\delta_n,\delta_1+\delta_2)$. This shows that $d_1(\vc{\delta},\ZZ^n_o)\leq d_1(M_{(1+2)}(\vc{\delta}),\ZZ^{n-1}_o)$.
\end{proof}

Now we will state two sufficient conditions for
a merging operation to satisfy the holonomy constraints:
one for simplicial points
and one for non-simplicial points.

\begin{lemma}[Non-simplicial merging]\label{cubmerging}
Let $n\geq 4$.
Let $\vc{\delta}\in \tcube^n$ be non-simplicial and suppose that  $\delta^{\pi}_i$ is an integer. Let $M$ be a merging operation
of type $M_{(i+j)},M_{(i-j)},M_{(j-i)}$ with $j\neq i$.
Then the point $M(\vc{\delta})$ lies is $\inte(\Hang^{n-1})$.
%
\end{lemma}
 
\begin{proof}
Let $\vc{c}$ be the center of $\tcube^n$.
We will prove that the image of the segment $[\vc{\delta^{\pi}},\vc{c}]$ 
under the map $M:\RR^n\rar\RR^{n-1}$ lies in $\tcube^{n-1}$ and at worst one of its points,
namely $M(\vc{c})$, sits at distance $1$ from $\ZZ_o^{n-1}$. 

Note first that $M([\vc{\delta^{\pi}},\vc{c}])$ belongs entirely
to some unit integer cube $\cube^{n-1}$. Indeed, $c_i+c_j$ and $c_i-c_j-2$ are integers, but $|(\delta^{\pi}_i+\delta^{\pi}_j)-(c_i+c_j)|\le 1$
and $|(\delta^{\pi}_i-\delta^{\pi}_j-2)-(c_i-c_j-2)|\le 1$.

Moreover, both $M(\vc{\delta^{\pi}})$ and $M(\vc{c})$
satisfy holonomy constraints, i.e. they belong to some $\tcube^{n-1}$. Indeed, using the fact that $\delta^{\pi}_i$ is integer and using  Lemma \ref{simpleproperties}(a) and Lemma \ref{lemma:int-merging}, we have
\[
d_1(M(\vc{\delta^{\pi}}), \ZZ_o^{n-1})=d_1(\vc{\delta^{\pi}}, \ZZ_o^n)>1
\]
and so $M(\vc{\delta^{\pi}})$ does not belong to a simplicial face of $\tcube^{n-1}$.
At the same time, $d_1(M(\vc{c}),\ZZ_o^{n-1})=\frac{n-2}{2}\ge 1$, since each coordinate of $\vc{c}$ is half-integral.

Since $\tcube^{n-1}$ is convex, and both ends of 
$M([\vc{\delta^{\pi}},\vc{c}])$ belong to the same $\tcube^{n-1}$, the whole segment belongs 
to it as well. 
On the other hand, $M(\vc{\delta^{\pi}})$ does not belong to a simplicial face of $\tcube^{n-1}$ and so
at worst one point of the segment 
$M([\vc{\delta^{\pi}},\vc{c}])$
belongs to a simplicial face,
namely $M(\vc{c})$. Clearly, this happens only if $n=4$.
\end{proof}

\begin{lemma}[Simplicial merging]\label{simplexmerging} 
Let $n\geq 5$.
Let $\vc{\delta}\in \tcube^n$ be simplicial and let $\vc{m}$ be a corresponding odd vertex. Consider a merging operation $M$
of positive type $M_{(i+j)}$ or of negative type $M_{(i-j)}$
and
suppose that $d_1(M(\vc{m}), M(\vc{\delta}))>1$. Then
$M(\vc{\delta})$ belongs to $\inte(\Hang^{n-1})$.
\end{lemma}

\begin{proof}
As in the proof of Lemma \ref{cubmerging}, note that the segment $M([\vc{\delta^{\pi}}, \vc{c}])$ belongs to some unit integer cube $\cube^{n-1}$. Let $\tcube^{n-1}$ be the truncated cube associated to such a $\cube^{n-1}$.

Since $n\ge 5$, the point $M(\vc{c})$ does not belong to any simplicial face
of $\tcube^{n-1}$, because $d_1(M(\vc{m'}),M(\vc{c}))=\frac{n-2}{2}>1$ for every vertex $\vc{m'}$ of $\cube^{n-1}$.
From
$
d_1(M(\vc{m}), M(\vc{\delta^{\pi}}))\le 
d_1(\vc{m},\vc{\delta^{\pi}})=1
$
we deduce that
the segment $M([\vc{\delta^{\pi}}, \vc{c}])$ intersects
the simplicial face of $\tcube^{n-1}$ corresponding to
the odd vertex $M(\vc{m})$.
Let us denote such a point of intersection by $\vc{z}$
and observe that the point $M(\vc{c})$ is the center of a non-simplicial face of $\tcube^{n-1}$. Thus, the segment $[M(\vc{c}),\vc{z}]$ lies inside $\tcube^{n-1}$ and in fact
it is not contained in any simplicial face of $\tcube^{n-1}$. 
As $M(\vc{\delta})$ lies in the interior of the segment
$[M(\vc{c}),\vc{z}]$, the conclusion follows.
\end{proof}

The following example shows why the restriction $n>4$ is important.

\begin{example}[Case $n=4$]
Let $a\in (0,\frac{1}{2})$ and let 
$\vc{\delta}=(a,-a,-1+a,-1+a)\in \RR^4$.
Observe that $\vc{\delta}$ is a defect vector that satisfies
positivity and holonomy constraints strictly; nevertheless,
only the mergings $M_{(1+2)}$, $M_{(1+3)}$ and $M_{(1+4)}$
preserve the positivity constraints. At the same time neither of these three positive mergings strictly preserves the holonomy constraints.
\end{example}

\subsubsection{Case (a): $\delta_1\leq 0$}

The following observation is elementary and so we omit the proof.

\begin{lemma}
The domain in $\RR^n$ obtained by
intersecting $\inte(\Aang^{n})$
with the cube  $(-1,0]^n$
is described by the following system of $2n+1$ inequalities:
\[
\begin{cases}
\delta_i\leq 0 & \text{for all $i=1,\dots,n$}\\
\sum_{j=1}^n \delta_j < 2\delta_i & \text{for all $i=1,\dots,n$}\\
\sum_{j=1}^n \delta_j>-2
\end{cases}
\]
\end{lemma}

In order to simplify the notation, up to rearranging the indices
we will assume
now on that $\delta_1\ge \delta_2\ge \dots\ge \delta_n$.

\begin{proposition}\label{allneg} 
Assume $n\geq 5$.
Suppose that $\vc{\delta}\in \inte(\Aang^n)\cap (-1,0)^n$.
Then $M_{(1+2)}(\vc{\delta})\in \inte(\Aang^{n-1})$.
\end{proposition}

\begin{proof}
Since $\sum _j\delta_j>-2$ and $n\ge 4$, we have $\delta_1+\delta_2>-1$.
Hence $M_{(1+2)}(\vc{\delta})$ satisfies the positivity constraints. 

Suppose now that $\vc{\delta}$ is a simplicial point. It is easy to see that the point  $-\vc{e_n}$ is an odd-integer point closest to $\vc{\delta}$.  
We have $d_1(M_{(1+2)}(\vc{\delta}), M_{(1+2)}(-\vc{e_n}))=d_1(\vc{\delta}, \vc{e_n})>1$. So by Lemma \ref{simplexmerging}, the vector $M_{(1+2)}(\vc{\delta})$
strictly satisfies the holonomy constraints.

Suppose now that $\vc{\delta}$ is not a simplicial point,
and so $\vc{\delta^{\pi}}$ strictly satisfies the holonomy constraints
and there exists an index $i$ such that $\delta_i^{\pi}$ is an integer.
As $-\vc{e_i}$ is an odd-integer vector, we have
$1+2\delta_i^{\pi}-\sum_{j}\delta^\pi_j=d_1(-\vc{e_i},\vc{\delta^{\pi}})>1$.
On the other hand,
$\sum_j \delta_j^{\pi}>-2$, because
$\sum_j \delta_j>-2$ and $n\geq 4$.
As a consequence,
$2\delta^{\pi}_i>\sum_{j}\delta^{\pi}_j>-2$
and so $\delta^{\pi}_i=0$.
The assumption $\delta_1\geq\delta_i$ necessarily implies
$\delta^\pi_1\geq\delta^\pi_i=0$ and so $\delta^\pi_1=0$,
because $\vc{\delta^\pi}\in[-1,0]^n$.
By Lemma \ref{cubmerging} we then conclude that $M_{(1+2)}(\vc{\delta})$ satisfies the holonomy constraints strictly. 
%
\end{proof}

\subsubsection{Case (b): $\delta_{1}>0$ and $\delta_{2}+\delta_{3}>-1$}

Up to rearranging the indices, assume
that $\delta_1\ge \delta_2\ge \dots\ge \delta_n$.

\begin{proposition}\label{partial} Suppose that $n\ge 5$ and  
$\vc{\delta}\in \inte(\Aang^n)$. Suppose moreover
that $\delta_{1}>0$ and $\delta_{2}+\delta_{3}>-1$. Then
there exist indices $i,j$ such that
$M_{(i+j)}(\delta)\in \inte(\Aang^{n-1})$.
\end{proposition}

\begin{proof}
Note that a positive merging $M_{(i+j)}$ preserves the sum
of the entries of the defect vector.
%
Thus, in order to prove that $M_{(i+j)}(\vc{\delta})$ satisfies the
positivity constraints, we only need to check
that $\delta_i+\delta_j>-1$.

Suppose first that $\vc{\delta}$ is non-simplicial, and
let $i$ be an index such that $\delta^{\pi}_i$ is integer. 
If $i\ne 1$, consider the positive merging $M_{(1+i)}$. According to Lemma \ref{cubmerging}, the vector
$M_{(1+i)}(\vc{\delta})$ satisfies the holonomy constraints strictly. At the same time 
$\delta_1+\delta_i>-1$, because $\delta_1>0$. 
Similarly, if $i=1$, we can consider the merging $M_{(1+2)}$.

Suppose now that $\vc{\delta}$ is simplicial and let
$\vc{m}$ be a closest odd-integer point.
Chose $i$ and $j$ distinct elements of 
$\{1,2,3\}$ so that $m_i-\delta_i$ and $m_j-\delta_j$ are of the same sign.
Then $d_1(M_{(i+j)}(\vc{m}), M_{(i+j)}(\vc{\delta}))=d_1(\vc{m},\vc{\delta})>1$ and so, according to 
Lemma \ref{simplexmerging}, the vector $M_{(i+j)}(\vc{\delta})$ strictly satisfies the holonomy constraints. At the same time $\delta_i+\delta_j>-1$ by hypothesis and so positivity constraints are also satisfied.
\end{proof}

\subsubsection{Case (c): $\delta_1>0$ and $\delta_2+\delta_3\le -1$}

Notice in particular that, in such a case, $\delta_j<0$ for all $j>1$.

The following technical definition is useful to clarify
when to apply a positive merging.

\begin{definition}
Let $n\geq 5$.
A vector $\vc{\delta}\in\inte(\Aang^n)$ with
$\delta_1\geq\dots\geq\delta_n$
is {\it{positively mergeable}} if the following three conditions are
not simultaneously satisfied:
\begin{itemize}
\item[(a)]
the vector $\vc{\delta}$ is simplicial;
\item[(b)]
there exists an integer $l\geq 1$ such that
\begin{itemize}
\item[(b1)]
$\vc{m}=(l,-1,-1,\dots,-1)$ is a vector 
in $\ZZ^n_o$ closest
to $\vc{\delta}$;
\item[(b2)]
$l>\delta_1$;
\item[(b3)]
$d_1(M_{(1+n)}(\vc{\delta}), M_{(1+n)}(\vc{m}))\leq 1$.
\end{itemize}
\end{itemize}
\end{definition}


\begin{proposition}[Positive merging]\label{plusmerging} 
Let $n\geq 5$ and let $\vc{\delta}\in \inte(\Aang^n)$ be a 
positively mergeable defect vector.
Then there exists indices $i,j$
such that $M_{(i+j)}(\vc{\delta})\in \inte(\Aang^{n-1})$.
\end{proposition}

\begin{proof}
Suppose first that $\vc{\delta}$ is the center of a unit cube.
Then all its entries are half-integers. Since $n\geq 5$, the vector
$M_{(i+j)}(\vc{\delta})$ sits at distance $\geq \frac{3}{2}$ from $\ZZ^{n-1}_o$ for any distinct $i$ and $j$, and so
$M_{(i+j)}(\vc{\delta})\in\inte(\Aang^{n-1})$.

Thus, now on we can assume that $\vc{\delta}$ is not the center of a
unit cube. Thanks to Propositions   \ref{allneg} and \ref{partial}, it is enough to treat the case $\delta_1\geq 0$ and $\delta_2+\delta_3\le -1$. \\

{\it{Case (a) violated.}}\\
The vector $\vc{\delta}$ is non-simplicial. 
Consider the positive merging $M_{(1+j)}$, where either $\delta^{\pi}_1$ or $\delta^{\pi}_j$ is integer. Since $\delta_1\ge 0$, the vector $M_{(1+j)}(\vc{\delta})$
satisfies the positivity constraints; moreover, by Lemma \ref{cubmerging} it also satisfies holonomy constraints strictly.\\

{\it{Assume now on that (a) is satisfied.}}\\
The vector $\vc{\delta}$ is simplicial: let $\vc{m}$ be a point in $\ZZ^n_o$ closest to $\vc{\delta}$, which is necessarily of the following type
$\vc{m}=l\vc{e_1}-\sum_{j\in J}\vc{e}_j$, for some integer $l\geq 0$
and some $J\subset\{2,3,\dots,n\}$.

Suppose $3\notin J$, and so $2\notin J$ either.
Then $d_1(\vc{m},\vc{\delta})\geq d_1(\vc{m}-\vc{e_2}-\vc{e_3},\vc{\delta})$. Replacing $\vc{m}$ by $\vc{m}-\vc{e_2}-\vc{e_3}$, we can thus assume $3\in J$ and then $\{3,4,\dots,n\}\subset J$.\\

{\it{Assume now on that
either $\vc{m}=l\vc{e_1}-(\vc{e_2}+\dots+\vc{e_n})$ or
$\vc{m}=l\vc{e_1}-(\vc{e_3}+\dots+\vc{e_n})$.}}\\
%
%
%
If $\delta_1-m_1=\delta_1-l\geq 0$, 
then $d_1(M_{1+3}(\vc{\delta}),M_{1+3}(\vc{m}))=d_1(\vc{\delta},\vc{m})>1$. Thus,
$M_{(1+3)}(\vc{\delta})$ strictly satisfies the holonomy constraints by Lemma \ref{simplexmerging}.

If $\delta_1-l<0$ and $m_2=0$, then
$d_1(M_{1+2}(\vc{\delta}),M_{1+2}(\vc{m}))=d_1(\vc{\delta},\vc{m})>1$.
As above, $M_{(1+2)}(\vc{\delta})$ strictly satisfied the holonomy constraints by Lemma \ref{simplexmerging}.

Thus, we are left to deal with the case $\delta_1<l$ and $m_2=-1$.\\

{\it{Assume now on that (b1) and (b2) are satisfied.}}\\
%
%
%
Since $d_1(M_{(1+i)}(\vc{m}), M_{(1+i)}(\vc{\delta}))\leq
d_1(M_{(1+j)}(\vc{m}), M_{(1+j)}(\vc{\delta}))$
for all $2\leq i<j\leq n$,
we can again conclude by applying
Lemma \ref{simplexmerging} to the operation $M=M_{(1+n)}$, unless
$d_1(M_{(1+n)}(\vc{m}), M_{(1+n)}(\vc{\delta}))\leq 1$,
that is unless $\vc{\delta}$ is not positively mergeable.
%
%
\end{proof}

The remaining cases can be taken care of by negative merging.
%
%

\begin{proposition}[Negative merging]\label{no+merging} 
Suppose that $n\geq 5$ and $\vc{\delta}$ is not positively mergeable.
Then $M_{(1-n)}(\vc{\delta})$ satisfies positivity and strict holonomy constraints. Moreover, $\delta_1$, $\delta_n$ and $\delta_1-\delta_n$ are not integers.
\end{proposition}

\begin{proof}
By our assumptions,
\begin{align*}
1\ge d_1(M_{(1+n)}(\vc{\delta}), M_{(1+n)}(\vc{m})) & =
|1-l+\delta_1+\delta_n|+(n-2)+\sum_{1<j<n}\delta_j\geq \\
& \geq n-l-1+\sum_j \delta_j
\end{align*}
that is, $\sum_j \delta_j\leq 2+l-n$.
%
Since $\sum \delta_j>-2$, we conclude that $n-l<4$. Since $\vc{m}$ is odd, the integer $n-l$ is even and so $l\geq n-2$. 
Moreover $\delta_1\geq l-1$, because $\vc{m}$ is a vector in $\ZZ^n_o$ closest to $\vc{\delta}$.

Recall that $M_{(1-n)}(\vc{\delta})=(\delta_2,\dots,\delta_{n-1},\delta_1-\delta_n-2)$. By the above computations,
\[
\delta_1-\delta_n-2\geq (l-1)-\delta_n-2\geq n-5-\delta_n>-1
\]
and 
\begin{align*}
\delta_2+\dots+\delta_{n-1}+(\delta_1-\delta_n-2)
& \geq \delta_1+\dots+\delta_{n-2}-2\geq \delta_1+(n-3)\delta_{n-2}-2>\\
& > (l-1)-(n-3)-2\geq -2
\end{align*}
which shows that $M_{(1-n)}(\vc{\delta})$ satisfies the positivity constraints.

On the other hand,
\begin{align*}
d_1(M_{(1-n)}(\vc{\delta}),M_{(1-n)}(\vc{m})) & =
|(l+1)-(\delta_1-\delta_n)|+\sum_{1<j<n}(1+\delta_j)=\\
& =(l-\delta_1)+\sum_{j>1}(1+\delta_j)=d_1(\vc{\delta},\vc{m})>1
\end{align*}
and so $M_{(1-n)}(\vc{\delta})$ strictly satisfies the holonomy
constraints by Lemma \ref{simplexmerging}.
%
%
%

Since $\vc{\delta}$ is not positively mergeable, it is
easy to see that
$\delta_1$ and $\delta_n$ cannot be integers.
In order to show that $\delta_1-\delta_n$ is not an integer either,
it is enough to prove that $\delta_1-\delta_n>l$,
because $\delta_1<l$ and $\delta_n\in (-1,0)$.
This can be easily verified, since
\begin{align*}
1 & \geq |1-l+\delta_1+\delta_n|+(n-2)+\sum_{1<j<n}\delta_j
\geq \\
&\geq l-1-\delta_1-\delta_n+(n-2)(1+\delta_n)>
(l-\delta_1+\delta_n)+1
\end{align*}
and so $\delta_1-\delta_n>l$.
\end{proof}

So finally we have achieved our task.

\begin{proof}[Proof of Theorem \ref{5points}]
If $\vc{\delta}$ is positively mergeable or the center of a unit cube, then
a positive merging operation will work by
Proposition \ref{plusmerging} and so (a) holds.
On the other hand, if $\vc{\delta}$ is not positively mergeable,
then Proposition \ref{no+merging} ensures that a negative operation
of type $M_{(i-j)}$ works and that, in this case, the involved defects
$\delta_i$, $\delta_j$ and $\delta_i-\delta_j$ are not integers. Thus, (b) holds.
\end{proof}

\section{Geometric constructions}

\subsection{Spherical bigons and triangles}\label{sec:big-tri}

\begin{definition}
Let $n\geq 2$. A {\it{spherical $n$-gon}} is a bordered surface
homeomorphic to the closed unit disc, endowed with a Riemannian metric of constant curvature $1$, whose boundary consists of $n$ geodesic arcs (called {\it{edges}}) that form inner angles $\pi\th_1,\dots,\pi\th_n$.
We will say that such an $n$-gon is {\it{convex}} if all $\th_i\leq 1$.
\end{definition}

Mimicking Definition \ref{def:deformable}, we can consider 
angle-deformability of spherical $n$-gons.

\begin{definition}
A metric $g$ on a spherical $n$-gon with inner angles
$\pi\cdot\vc{\th}$ is {\it{angle-deformable}}
if there exists a neighbourhood $\Nang$ of $\vc{\th}\in\RR^n$
and a continuous family $\Nang\ni \vc{\nu}\mapsto g_{\vc{\nu}}$ of spherical metrics on the $n$-gon
such that $g_{\vc{\nu}}$ has conical singularities of angles 
$\pi \cdot\vc{\nu}$ for all $\vc{\nu}\in \Nang$ and
$g_{\vc{\th}}=g$.
\end{definition}

We will refer to a $2$-gon, $3$-gon and $4$-gon
respectively as a ``bigon'', ``triangle'' and ``quadrilateral''.
If $x_i,x_{i+1}$ are consecutive vertices of an $n$-gon, then
we will denote by $|x_i x_{i+1}|$ the length of the edge joining them.

\begin{notation}
Let $S$ be a compact spherical surface possibly with boundary and let $\gamma$ be a curve inside $S$. We say that $S'$ is {\it{obtained from $S$ by cutting along $\gamma$}} if $S'$ is the compact spherical
surface (possibly with boundary) obtained as a metric completion
of $S\setminus\gamma$.
\end{notation}

\begin{notation}
Let $S$ be a compact spherical surface with boundary.
The double $DS$ of $S$ is the spherical surface obtained
by gluing $S$ with $\bar{S}$ (another copy of $S$, with the opposite orientation) isometrically along their boundary.
We will say that $S$ is angle-deformable (resp. non-coaxial) if $DS$ is
(resp. if $DS$ has non-coaxial holonomy).
\end{notation}


\subsubsection{Bigons}

Pick spherical coordinates $\psi\in[0,2\pi)$ and $\phi\in[0,\pi]$ on $\Sph$.
Given $0<\alpha<1$ and $0<r\leq \pi$, we denote by
$B_\alpha(r)=\{ (\psi,\phi)\,|\,\text{$\psi\in[0,\pi\alpha]$ and $\phi\in[0,r)$}\}$ and by $\ol{B}_\alpha(r)$ its closure.
For $\alpha\geq 1$, we let $B_\alpha(r)$ be obtained from $k$ copies $B_1,\dots,B_k$ of $B_{\alpha/k}(r)$ 
by gluing one geodesic side of $B_i$ to one geodesic side of $B_{i+1}$ for $i=1,\dots,k-1$. Analogously for $\ol{B}_\alpha(r)$.

\begin{definition}
Let $\alpha>0$ and $r\in(0,\pi)$.
The {\it{standard (open) $r$-neighbourhood of a vertex of a spherical polygon}} of angle $\pi\alpha$ is the surface with boundary $B_\alpha(r)$. The {\it{standard (open) $r$-neighbourhood of a cone point}} of angle $2\pi\alpha$ is the spherical surface
$S_\alpha(r)$ obtained by doubling $B_{\alpha}(r)$.
In an analogous way we define the standard closed $r$-neighbourhoods.
\end{definition}
%

\begin{lemma}[Existence of bigons]\label{lemma:bigon}
\begin{itemize}
\item[(a)]
For every $\alpha>0$,
there exists a bigon $B_\alpha$
with both angles
$\pi\alpha$ and with cone points at distance $\pi$.
Such $B_\alpha$ is angle-deformable.
\item[(b)]
Let $d>0$ be an integer.
There exists a continuous family of bigons
$(0,2\pi)\ni \ell \mapsto B(d,\ell)$ with both angles $d\pi$ and the two sides of lengths
$(\ell,2\pi-\ell)$ for $d$ is odd, or $(\ell,\ell)$ for $d$ even.
%
\end{itemize}
\end{lemma}

\begin{notation}
In what follows, we will refer to a bigon to type (a) above as an {\it{(ordinary) bigon}} and to
a bigon of type (b) as an {\it{exceptional bigon}}.
\end{notation}

\begin{proof}[Proof of Lemma \ref{lemma:bigon}]
The bigon in (a) is clearly $B_\alpha=\ol{B}_\alpha(\pi)$ and deformability is obvious.

About claim (b), let $D\subset \Sph$ be a closed unit hemisphere.
The bigon $B(1,\ell)$ is easily obtained from $D$
by marking two points $x_1,x_2$ on $\partial D$ that break $\partial D$ into
two geodesic arcs of lengths $\ell$ and $2\pi-\ell$.
For $d>1$, consider cyclic cover $\tilde{S}\rar \Sph$ of degree $d$ branched at $x_1,x_2$.
If $B_{(1)},\dots,B_{(d)}$ are the lifts of $D$ and
$B'_{(1)},\dots,B'_{(d)}$ are the lifts of the other hemisphere $\overline{\Sph\setminus D}$, then $B(d,\ell)$ is obtained as
the union of all $B_{(2i+1)}$ and $B'_{(2i)}$.
%
\end{proof}

%
%

By doubling the bigons constructed above, we immediately have
the following.

\begin{corollary}[Spheres with two conical points]\label{cor:two-pointed}
\begin{itemize}
\item[(a)]
For every $\alpha>0$,
there exists a spherical surface $S_\alpha:=DB_\alpha$ homeomorphic to a sphere
with both angles
$2\pi\alpha$ and with cone points at distance $\pi$. Such $S_\alpha$ are angle-deformable.
\item[(b)]
Let $d>0$ be an integer.
There exists a continuous family
$(0,\pi)\ni \ell \mapsto DB(d,\ell)$ of
spherical surfaces
with cone points at distance $\ell$ and
both angles $2\pi d$.
\end{itemize}
\end{corollary}

\begin{center}
\begin{figurehere}
\psfrag{x1}{$x_1$}
\psfrag{x2}{$x_2$}
\psfrag{l}{$\ell$}
\psfrag{Ba}{$DB_\alpha$}
\psfrag{B(2,l)}{$DB(2,\ell)$}
\includegraphics[width=0.5\textwidth]{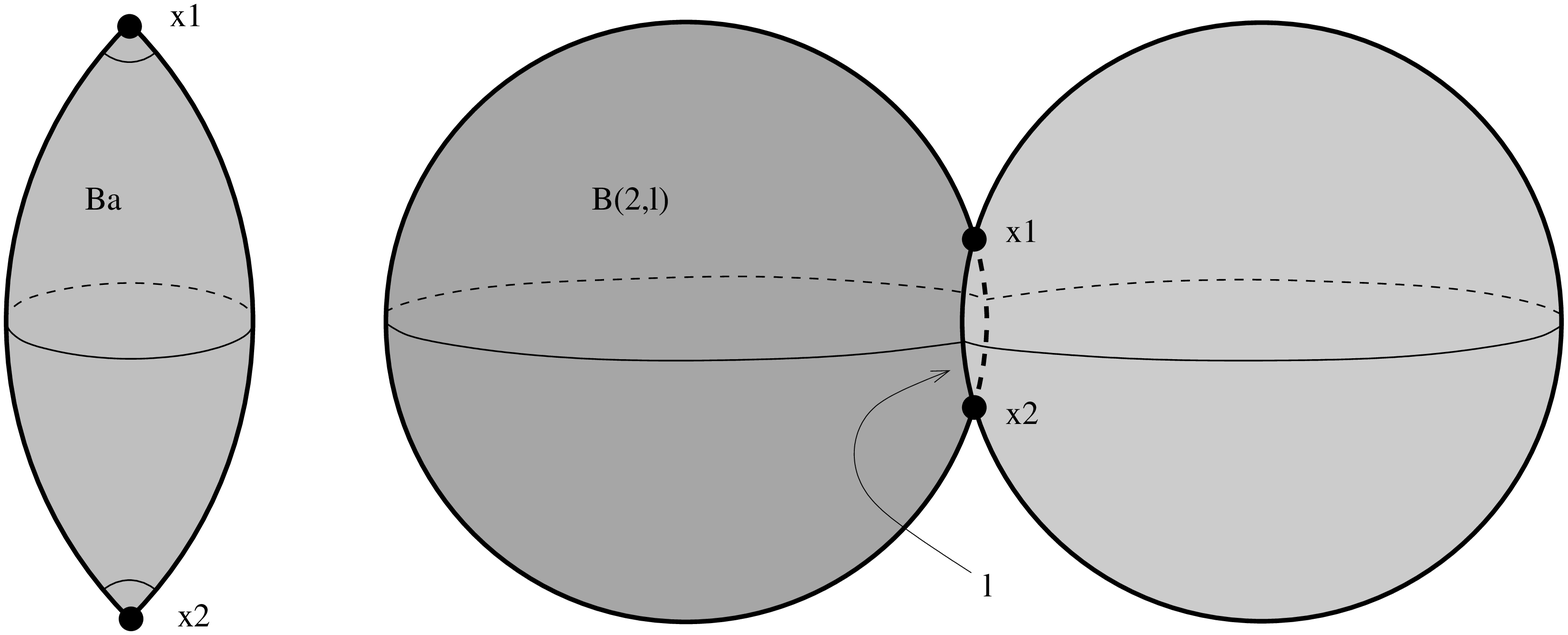}
\caption{{\small An ordinary and an exceptional sphere with two conical points.}}\label{fig:bigons}
\end{figurehere}
\end{center}

\begin{remark}
It can be easily seen that all bigons are of types (a) and (b) described in Lemma \ref{lemma:bigon}. Analogously, surfaces of curvature $1$ homeomorphic to a sphere with two conical points can be obtained by doubling such bigons as in Corollary \ref{cor:two-pointed} (see Troyanov \cite{troyanov:bigons}).
\end{remark}



As an application, here we characterize spherical surfaces with
non-integral angles and reducible holonomy.

\begin{lemma}[Metrics with reducible holonomy]\label{reducibledecsription} 
Let $S$ be a spherical surface
with conical singularities $x_1,\dots,x_n$ of angles $2\pi\th_1,\dots,2\pi\th_n$.
Suppose that the holonomy 
$\rho: \pi_1(\dot{S})\to \SO(3,\RR)$ 
is reducible
and that no $\th_i$ is integral.
Then there is a subset $J\subseteq\{1,\dots,n\}$ and
geodesic graph $\Graph\in S$ such that
\begin{itemize}
\item[(a)]
$S\setminus \Graph$ is the disjoint
union of the disks $S_{\th_i}(\pi/2)$ for $i\in  J^c$
and possibly of some hemispheres;
\item[(b)]
for every $j\in J$ the conical point $x_j$ is contained in $\Graph$
and has conical angle $2\pi(k_j+\frac{1}{2})$ for some $k_j\in\ZZ_{\geq 0}$.
\end{itemize}
\end{lemma} 
\begin{proof}
Let $\widetilde{\dot{S}}\rar \dot{S}$ 
be the universal cover.
Consider the developing map $\dev:\widetilde{\dot{S}}\rar \Sph$
and the holonomy representation $\rho:\pi_1(\dot{S})\rar\SO(3,\RR)$
associated to the given spherical metric.
Since $\rho$ is non-trivial and reducible, there is a plane $P\subset \RR^3$  and an orthogonal line $L$ invariant under
$\rho(\pi_1(\dot{S}))$. Clearly the map $\dev$ does not
ramify over $\Sph\setminus(P\cup L)$.
%

Define $\widetilde{\Graph}$ as the 
$\pi_1(\dot{S})$-invariant geodesic graph $\dev^{-1}(\Sph\cap P)\subset \widetilde{\dot{S}}$, which descends to a geodesic graph on $\dot{S}$. The closure $\Graph$ of such graph passes through the
conical points $\{x_j\,|\,j\in J\}$ for some $J\subseteq\{1,\dots,n\}$.

Let $R_j\in\SO(3,\RR)$ be the holonomy along a loop that simply
winds about $x_j$. Such an $R_j$ is not the identity and its axis lies in $P$; moreover, $R_j$ preserves $L$: the only possibility is that
$R_j$ is a rotation by an angle $\pi$ and so $\th_j=k_j+\frac{1}{2}$ for some $k_j\in\ZZ_{\geq 0}$.

Consider now an $i\in J^c$. Let $D_i$ be the component of $S\setminus\Graph$ that contains $x_i$ and let 
$\widetilde{\dot{D}}_i\rar\dot{D}_i$ be the universal cover of
$\dot{D}_i=D_i\setminus\{x_i\}$. The developing map
restricts to a cover of $\dot{D}_i$ over a component of
$\Sph\setminus (P\cup L)$ and so $D_i$
is isometric to $S_{\th_i}(\pi/2)$.

Let $D$ be a component of $S\setminus\Graph$ that does not
contain any $x_i$. Then $\dev$ induces an isomorphism between
$D$ a component of $\Sph\setminus P$, and so $D$ is a hemisphere.
\end{proof}

\subsubsection{Triangles}

The following theorem follows from \cite[Theorem 3]{eremenko:three}.

\begin{theorem}[Existence of triangles]\label{membrane} 
Let $\vc{\th}=(\th_1,\th_2,\th_3)$
be a triple of real numbers satisfying holonomy constraints (\ref{theinquality})
strictly and the positivity constraints (\ref{gauss-bonnet}). Then there exists an angle-deformable non-coaxial spherical triangle with 
inner angles $\pi\cdot\vc{\th}$.
\end{theorem}

Since we will need spherical triangles later, here we give
a short constructive proof of the above theorem.
The wished triangle
is assembled from pieces constructed 
in Lemma \ref{lemma:bigon}, Lemma \ref{twoconstructions} 
and Corollary \ref{twocases}.

\begin{lemma}[Existence of convex triangles]\label{triangle} 
Let $\vc{\th}=(\th_1,\th_2,\th_3)\in (0,1)^3$. 
A convex spherical triangle with angles $\pi\cdot\vc{\th}$ exists 
if and only if both conditions are satisfied:
\begin{itemize}
\item[(i)] 
the numbers $(1-\th_1,1-\th_2,1-\th_3)$ satisfy the triangular inequality;
\item[(ii)]
the following inequality holds: $\th_1+\th_2+\th_3-1>0$.
\end{itemize}
Moreover, such convex triangles are angle-deformable and non-coaxial.
\end{lemma}
\begin{proof}
Assume a convex spherical triangle with angles $\pi\cdot\vc{\th}$ exists. Then condition (i) holds since its dual spherical triangle has edges of lengths $\pi(1-\th_i)$. Moreover, (ii) holds too, since $\pi(\th_1+\th_2+\th_3-1)$ is the area of the triangle.
Vice versa, fix a hemisphere $D$ and a point $x_2$ on $\partial D$.
It is easy to see that, for every triple $(\th_1,\th_2,\th_3)$ satisfying
(i) and (ii), one can realize a triangle with such angles
as the convex hull in $D$ of $x_1,x_2,x_3$, for suitable
$x_3\in\partial D$ and $x_1\in D$, see Figure \ref{fig:triangles}(a).
Such triangles are angle-deformable by construction; moreover, since $\th_i\in(0,1)$ and the triangle is inscribed in a hemisphere, it is immediate to see that it is non-coaxial.
\end{proof}

This lemma settles Theorem \ref{membrane} for triangles with $\th_i<1$.
Indeed, inside the unit cube $[0,1]^3$
Inequalities (\ref{theinquality}) 
and conditions (i) and (ii) 
describe the tetrahedron with vertices 
$(1,0,0)$, $(0,1,0)$, $(0,0,1)$, $(1,1,1)$.

\begin{center}
\begin{figurehere}
\psfrag{x1}{$x_1$}
\psfrag{x2}{$x_2$}
\psfrag{x3}{$x_3$}
\psfrag{(a)}{(a)}
\psfrag{(b)}{(b)}
\psfrag{T}{$T$}
\psfrag{T'}{$T'$}
\includegraphics[width=0.6\textwidth]{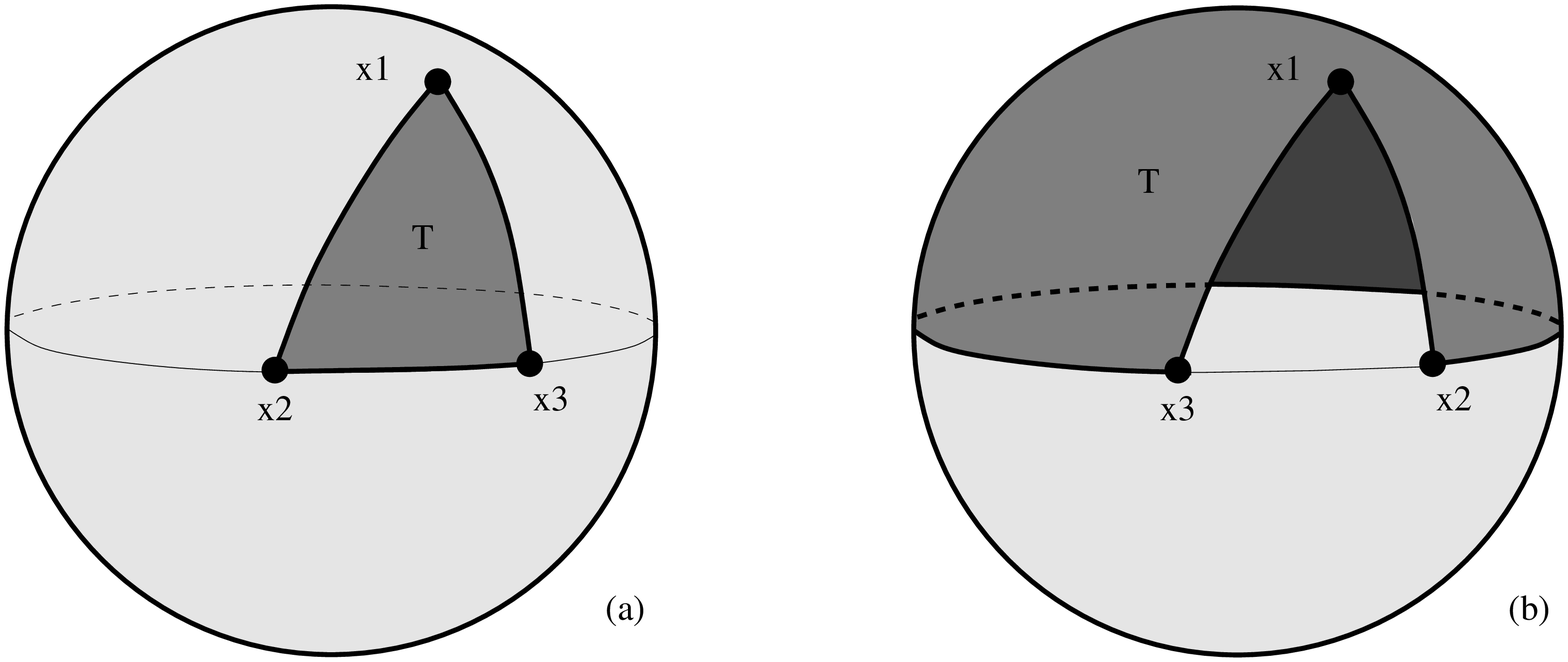}
\caption{{\small A convex triangle (a) and a complement (b) of a convex triangle in a hemisphere.}}\label{fig:triangles}
\end{figurehere}
\end{center}

\begin{corollary}[Triangles with small angles]\label{twocases} 
Suppose $\vc{\th}=(\th_1,\th_2,\th_3)$ 
satisfy Inequalities (\ref{theinquality}) strictly and
it belongs to the domain $\Pi^3:=[0,2]\times [0,1] \times [0,1]\subset \RR^3$.
Then there exists an angle-deformable non-coaxial spherical triangle with angles $\pi \cdot\vc{\th}$.
\end{corollary}
\begin{proof}
As we explained, the case when all $\th_i$ are less than
$1$ follows from Lemma \ref{triangle}. Suppose that $1<\th_1<2$. 
Then it is easy to check that the triple 
$\vc{\th'}=(2-\th_1, 1-\th_3, 1-\th_2)$ satisfies Inequalities (\ref{theinquality})
as well. Consider a convex spherical  triangle $T'\subset \Sph$ 
with angles $\pi(\th'_1,\th'_3,\th'_2)$ at the vertices $(x_1,x_3,x_2)$
and let $E\subset \Sph$ be the maximal circle that contains
the vertices $x_2$ and $x_3$.
Cut $\Sph$ along $E$ and let $D$ be the component that contains
$\mathrm{int}(T')$.
Then the spherical triangle obtained from
$D\setminus T'$ by metric completion
has angles $\pi\cdot\vc{\th}$, see Figure \ref{fig:triangles}(b).
Angle-deformability and non-coaxiality of $T'$ implies that the constructed triangles is angle-deformable and non-coaxial too.
\end{proof}

\begin{center}
\begin{figurehere}
\psfrag{y1}{$y_1$}
\psfrag{y2}{$y_2$}
\psfrag{y3}{$y_3$}
\psfrag{b1}{$b_1$}
\psfrag{b2}{$b_2$}
\psfrag{g}{$\gamma$}
\psfrag{a}{$\alpha$}
\psfrag{1-a}{$1-\alpha$}
\psfrag{1}{$2$}
\psfrag{d}{$2d$}
\psfrag{T}{$T(d,\ell,\alpha)$}
\psfrag{D}{$D$}
\includegraphics[width=0.6\textwidth]{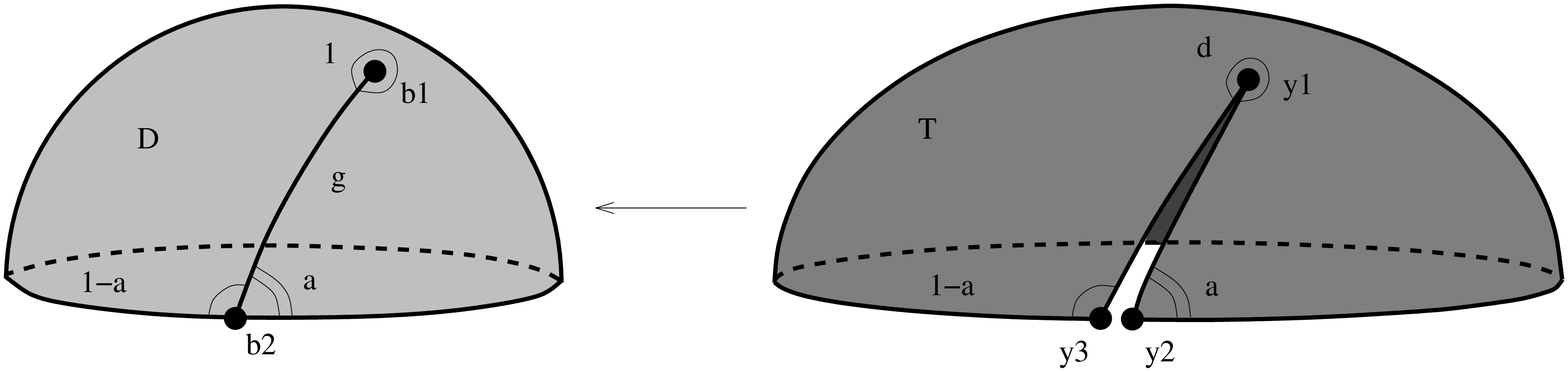}
\nopagebreak[4]
\caption{{\small The triangle $T(d,\ell,\alpha)$.}}\label{fig:triangle2}
\end{figurehere}
\end{center}

\begin{lemma}[The triangles $T(d,\ell,\alpha)$]\label{twoconstructions}
Let $d$ be a positive integer and let $0<\ell<\pi$ and $0<\alpha<1$.
There exists a spherical triangle $T(d,\ell,\alpha)$ of vertices $y_1,y_2,y_3$
with edges of lengths $|y_1y_2|=|y_1y_3|=\ell$ and $|y_2y_3|=2\pi d$
and angles $\pi (2d)$, $\pi\alpha$, and $\pi(1-\alpha)$
at the vertices $y_1$, $y_2$, $y_3$ correspondingly.
\end{lemma}

\begin{proof}
%
Let $D\subset \Sph$ be a closed hemisphere.
Choose a geodesic segment $\gamma$ on $D$ of length $\ell$, with one endpoint $b_1$ in the interior of $D$ and the other endpoint
$b_2$ on the boundary of $D$ and forming
angles $\pi\alpha$ and $\pi(1-\alpha)$ with $\partial D$.
Consider now a ramified degree $d$ cover $\tilde{D}\rar D$ 
that has an order $d$ branching at $b_1$. The wanted
spherical triangle is obtained by cutting $\tilde{D}$ along $\tilde{\gamma}$,
namely one of the $d$ geodesic preimages of $\gamma$,
as illustrated in Figure \ref{fig:triangle2}.
\end{proof}

Denote by $\Gamma^3\subset \ZZ_{\geq 0}^3$  the additive semigroup consisting of elements 
$\vc{m}=(m_1,m_2,m_3)$ such that $m_1\ge m_2\ge m_3$ and $m_1+m_2+m_3\in 2\ZZ$.

\begin{lemma}\label{Pi+} 
The subset $\{\vc{\th}=(\th_1,\th_2,\th_3)\in\RR^3\,|\,
\th_1\geq\th_2\geq\th_3\}\subset \RR_+^3$ 
is contained in the union $\bigcup_{\vc{m}\in \Gamma^3}(\Pi^3+\vc{m})$.
\end{lemma}
The previous lemma is completely elementary, and we omit the proof.

\begin{proof}[Proof of Theorem \ref{membrane}.]
Let $\vc{\th}=(\th_1,\th_2,\th_3)\in \RR^3_+$
be a triple satisfying Inequalities (\ref{theinquality}) strictly.
After reordering the coordinates,
we can assume $\th_1\ge \th_2\ge \th_3$.
We will construct now a spherical triangle with angles 
$\pi\cdot\vc{\th}$.

By Lemma \ref{Pi+},
there exists $\vc{m}\in\Gamma^3$ and $\vc{\th'}\in\Pi^3$
such that $\vc{\th}=\vc{\th'}+\vc{m}$.
Since Inequality ({\ref{theinquality}})
is invariant with respect to translations by integer vectors $\vc{m}$ with $m_1+m_2+m_3$ even, by Corollary \ref{twocases} there exists 
an angle-deformable non-coaxial spherical triangle $T'$ with 
angles $\pi\th'_1,\pi\th'_2,\pi\th'_3$ at its
vertices $x'_1,x'_2,x'_3$. Since no $\th'_i$ is an integer,
no edge of $T'$ has length multiple of $\pi$.

Consider now separately two cases.

{\it{Case (a): $m_1>m_2+m_3$.}}
The construction is illustrated in Figure \ref{fig:triangle-a}.
\begin{center}
\begin{figurehere}
\psfrag{x'1}{$x'_1$}
\psfrag{x'2}{$x'_2$}
\psfrag{x'3}{$x'_3$}
\psfrag{x1}{$x_1$}
\psfrag{x2}{$x_2$}
\psfrag{x3}{$x_3$}
\psfrag{B''}{$B''$}
\psfrag{B'''}{$B'''$}
\psfrag{x''1}{$x''_1$}
\psfrag{x''2}{$x''_2$}
\psfrag{x''3}{$x''_3$}
\psfrag{T}{$T$}
\psfrag{T'}{$T'$}
\psfrag{T''}{$T''$}
\psfrag{t'1}{$\th'_1$}
\psfrag{t'2}{$\th'_2$}
\psfrag{t'3}{$\th'_3$}
\psfrag{1-t'3}{$1-\th'_3$}
\psfrag{2d}{$2d$}
\psfrag{m2}{$m_2$}
\psfrag{m3}{$m_3$}
\includegraphics[width=0.6\textwidth]{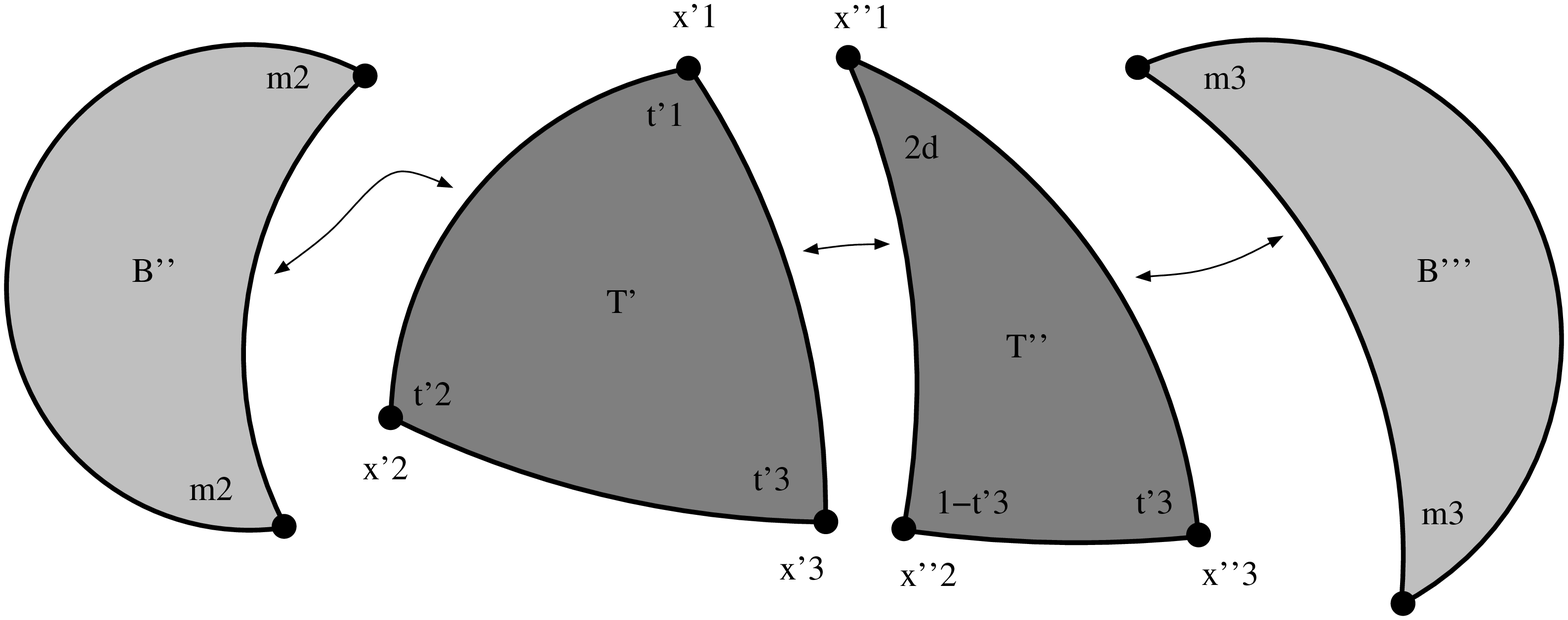}
\caption{{\small Building the triangle in case (a) of Theorem \ref{membrane}.}}\label{fig:triangle-a}
\end{figurehere}
\end{center}

Let $d=\frac{1}{2}(m_1-m_2-m_3)$. 
By Lemma \ref{twoconstructions}, there exists
a spherical triangle $T''=T(d,|x'_1x'_3|,1-\th'_3)$ 
with vertices $x''_1,x''_2,x''_3$, angles $\pi(2d,1-\th'_3,\th'_3)$
and $|x''_1 x''_2|=|x''_1 x''_3|=|x'_1x'_3|$,
Denote by $T$ the triangle with vertices $x_1,x_2,x_3$
obtained by identifying
the side $x''_1 x''_3$ of $T''$
with the side $x'_1x'_3$ of $T'$.
The angle at vertex $x_1$ of $T$ corresponding to $x'_1\sim x''_1$
is $\pi(2d+\th'_1)$, the angle at $x_2$ corresponding to $x'_2$
is $\pi\th'_2$ and the angle at $x_3$ corresponding to $x''_3$ is
$\pi\th'_3$.
Finally, take two exceptional
bigons $B''=B(m_2, |x_1x_2|)$ and $B'''=B(m_3, |x_1x_3|)$ and glue them with $T$ by isometrically identifying one side of $B''$ to $x_1x_2$ and 
one side of $B'''$ to $x_1x_3$. 
Angle-deformability of $T'$ implies angle-deformability of $T$ and so of the wished triangle. Since $|x_1 x_2|=|x'_1 x'_2|$ is not a multiple of $\pi$ and no $\th_i$ is an integer, the triangle $T$ is non-coaxial
and so is the constructed triangle.

{\it{Case (b): $m_1\le m_2+m_3$.}}
The construction is illustrated in Figure \ref{fig:triangle-b}.
\begin{center}
\begin{figurehere}
\psfrag{x'1}{$x'_1$}
\psfrag{x'2}{$x'_2$}
\psfrag{x'3}{$x'_3$}
\psfrag{x1}{$x_1$}
\psfrag{B'}{$B'$}
\psfrag{B''}{$B''$}
\psfrag{B'''}{$B'''$}
\psfrag{T'}{$T'$}
\psfrag{t'1}{$\th'_1$}
\psfrag{t'2}{$\th'_2$}
\psfrag{t'3}{$\th'_3$}
\psfrag{m1}{$m_1$}
\psfrag{m2}{$m_2$}
\psfrag{m3}{$m_3$}
\includegraphics[width=0.6\textwidth]{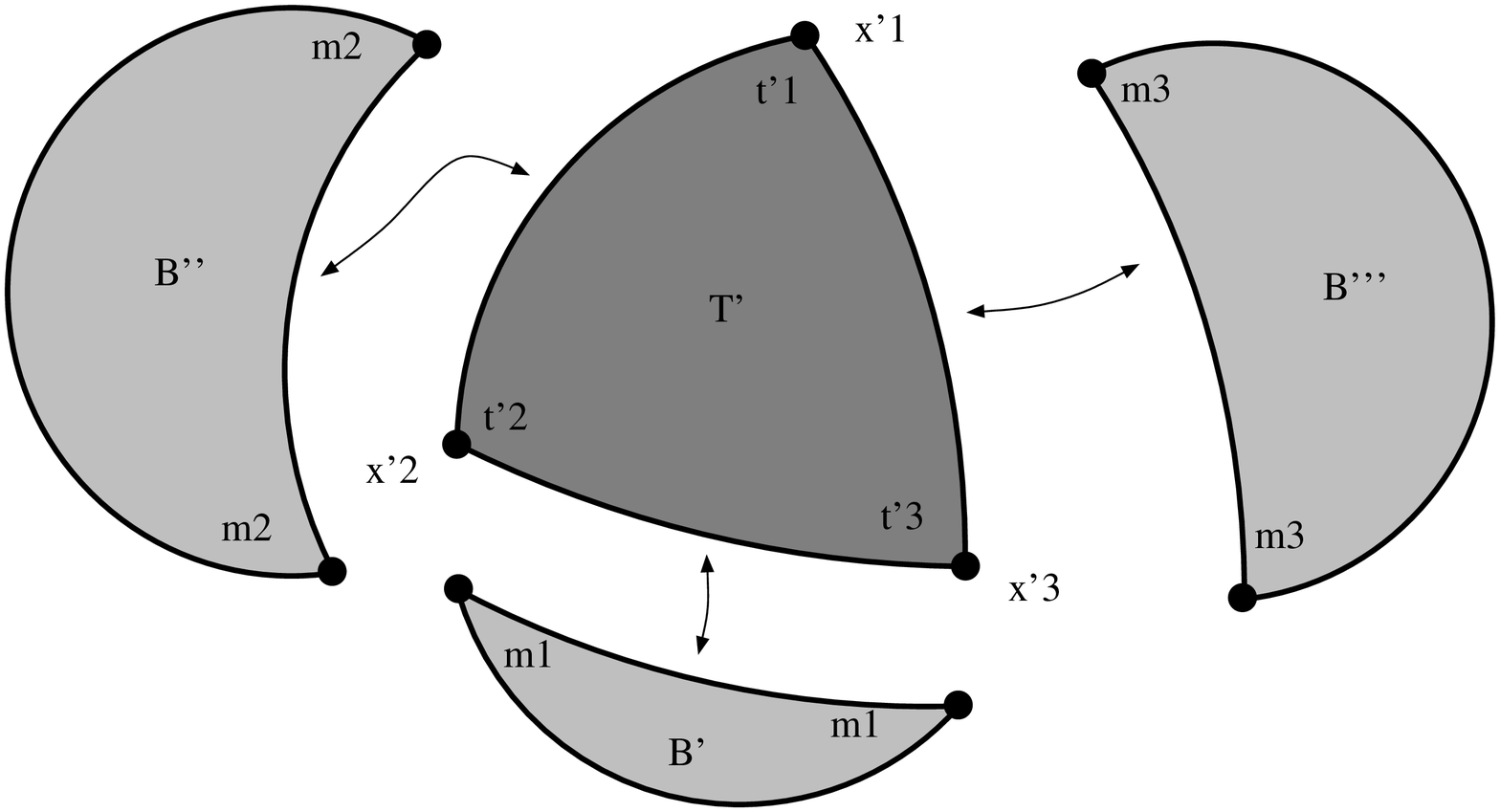}
\caption{{\small Building the triangle in case (b) of Theorem \ref{membrane}.}}\label{fig:triangle-b}
\end{figurehere}
\end{center}

In this case, set $d_1=\frac{1}{2}(m_2+m_3-m_1)$,
$d_2=\frac{1}{2}(m_3+m_1-m_2)$, $d_3=\frac{1}{2}(m_1+m_2-m_3)$.
Take the three exceptional bigons $B'=B(d_1, |x'_2x'_3|)$, $B''=B(d_2,|x'_3x'_1|)$, $B'''=B(d_3,|x'_1x'_2|)$
and glue them with $T'$
by isometrically identifying a side of $B'$ with $x'_2 x'_3$, a side of $B''$
with $x'_3x'_1$ and a side of $B'''$ with $x'_1x'_2$.
As before, angle-deformability and non-coaxiality of $T'$ implies angle-deformability and non-coaxiality of the constructed triangle.
\end{proof}

The above existence theorem for triangle allows to draw
the following conclusion about $3$-punctured spheres.

\begin{corollary}[Existence of $3$-punctured spherical metrics]
Let $\vc{\th}=(\th_1,\th_2,\th_3)$
be a triple of real numbers satisfying holonomy (\ref{theinquality})
strictly and the positivity constraints (\ref{gauss-bonnet}). Then there exists an angle-deformable non-coaxial spherical surface $S$ of genus $0$
with conical singularities of angle $2\pi\cdot\vc{\th}$.
\end{corollary}
\begin{proof}
The wished spherical surface is obtained by doubling
a spherical triangle with angles $\pi\cdot\vc{\th}$, whose existence
relies on Theorem \ref{membrane}.
\end{proof}

Indeed, a little more is true.

\begin{lemma}[Double of spherical triangles]\label{lemma:double-triangles}
Let $S$ be a sphere with distinct points $x_1,x_2,x_3$
endowed with a spherical metric $g$ with conical singularities of angle $2\pi\th_i$ at $x_i$ and non-coaxial holonomy. Then:
\begin{itemize}
\item[(a)]
$g$ is the unique spherical metric in its conformal class with
such conical singularities;
\item[(b)]
the spherical surface $(S,g)$ is obtained by doubling a spherical triangle.
\end{itemize}
\end{lemma}
\begin{proof}
The uniqueness of $g$ claimed in (a) was already noticed in \cite{eremenko:three}.
In fact, let $g'$ be a spherical metric conformal to $g$ and denote by
$J$ the underlying conformal structure.
Both
$g$ and $g'$ induce $\CC\PP^1$-structures $\Xi,\Xi'$ on the Riemann surface
$(\dot{S},J)$: their difference is thus encoded in a Schwarzian derivative
$\sigma(\Xi,\Xi')$, which is a holomorphic quadratic differential on $\dot{S}$.
A direct computation shows that $\sigma(\Xi,\Xi')$ has at most simple poles
at $x_1,x_2,x_3$, because $g$ and $g'$ have the same angles at the $x_i$,
and so $\sigma(\Xi,\Xi')\equiv 0$. This implies that the two $\CC\PP^1$-structures and so their holonomy representations agree.
Moreover, the developing maps of $g$ and $g'$ are conjugate through a M\"obius transformation $\tau\in\PSL(2,\CC)$ that commutes with the holonomy subgroup of $\SO(3,\RR)$. Since we assumed the $\SO(3,\RR)$-holonomy to be non-coaxial, 
Lemma \ref{lemma:non-coaxial}(b) ensures that
$\tau$ must lie in $\SO(3,\RR)$ and so $g=g'$.

As for (b), we remark that $(S,J)$ is biholomorphic to $\CC\PP^1$
through a map that takes $x_1,x_2,x_3\in S$ to $[1:0],[1:1],[0:1]\in\CC\PP^1$.
The conjugation is an anti-holomorphic (and so conformal) transformation of
$\CC\PP^1$ that fixes $[1:0],[1:1],[0:1]$ and so transports to a
conformal involution $\iota$ of $S$ that fixes $x_1,x_2,x_3$.
By (a), the metric $g$ must be fixed by $\iota$, which is thus an isometry
of $(S,g)$. It is then immediate to check that $S$ is isometric to the double $DT$, where $T$ is the spherical triangle $S/\!\iota$.
\end{proof}

\subsubsection{Almost degenerate triangles}

Spherical triangles can degenerate in several ways.
We are interested in describing two such degenerations:
in the first case, the triangle degenerates to an ordinary bigon;
in the second case, the triangle degenerates to a ``double bigon'',
that is the union of two ordinary bigons sharing a common vertex.

\begin{definition}
A spherical polygon is {\it{$r$-wide}} at a vertex $x_i$
of angle $\pi\alpha$ if the closed ball centered at $x_i$ of radius $r$ is isometric to $\ol{B}_\alpha(r)$ and does not contain any marked point other than $x_i$.
A spherical surface if {\it{$r$-wide}} at a cone point $x_i$ of angle
$2\pi\alpha$ if the closed ball centered at $x_i$ of radius
$r$ is isometric to $\ol{S}_\alpha(r)$ and does not contain any marked point other than $x_i$.
\end{definition}

\begin{notation}
If a spherical surface $S$ is $r$-wide at a conical point $y$ of angle $2\pi\alpha$, then
we denote by $U_y(r)$ the complement in $S$ of the open neighbourhood
of $y$ isometric to $B_\alpha(r)$.
\end{notation}


The triangles we are going to describe are needed in the surgery
operations that will split a conical point into a pair of conical singularities. In order to prove the angle-deformability of the so-constructed spherical surface, we need the following properties
from our triangles.

\begin{definition}
A spherical triangle $(T,g)$ with vertices $x_1,x_2,x_3$ and
angles $\pi\cdot\vc{\th}=\pi(\th_1,\th_2,\th_3)$ is {\it{$(x_1,x_2)$-angle-deformable}} if there exists a neighbourhood $\Nang'\subset\RR^2$ of $(\th_1,\th_2)$, a continuous map $\theta_3:\Nang'\rar\RR$
such that
$\theta_3(\th_1,\th_2)=\th_3$
and a continuous family of metrics $g_{\vc{\nu}}$ parametrized
by $\vc{\nu}\in \Nang'$ such that
$g_{(\th_1,\th_2)}=g$ and $g_{\vc{\nu}}$ has angles
$\pi(\nu_1,\nu_2,\theta_3(\nu_1,\nu_2))$.
\end{definition}

Notice that angle-deformability is clearly stronger than $(x_1,x_2)$-angle-deformability. On the other hand, the above definition is particularly meaningful for a $\vc{\th}$ that only weakly
satisfies the holonomy constraints, in which case absolute angle-deformability cannot hold.

Also, we recall that a spherical surface of genus $0$ with $3$ conical points is obtained by doubling a spherical triangle. Thus, angle-deformability of the surface is equivalent to angle-deformability of the triangle.


\begin{proposition}[Triangles close to an ordinary bigon]\label{sumtriangle}
Let $\th_1,\th_2,\th_3>0$ with $\th_3=\th_1+\th_2-1$.
For every $\varepsilon>0$ there exist $\eta\in(-\varepsilon,\varepsilon)$ and 
a spherical triangle $T$ 
with angles $\pi(\th_1, \th_2,\th_3+\eta)$ and vertices $x_1,x_2,x_3$, 
which is $\pi(1-\varepsilon)$-wide at $x_3$
and $(x_1,x_2)$-angle-deformable.
\end{proposition}

\begin{proof}
We divide the proof into four cases, illustrated in Figure \ref{fig:degenerate-tri-1}.

\begin{center}
\begin{figurehere}
\psfrag{1}{$1$}
\psfrag{(a)}{(a)}
\psfrag{(b)}{(b)}
\psfrag{(c)}{(c)}
\psfrag{(d)}{(d)}
\psfrag{x'1}{$x'_1$}
\psfrag{x'2}{$x'_2$}
\psfrag{x'3}{$x'_3$}
\psfrag{x1}{$x_1$}
\psfrag{x2}{$x_2$}
\psfrag{x3}{$x_3$}
\psfrag{y3}{$y_3$}
\psfrag{B'}{$B'$}
\psfrag{B''}{$B''$}
\psfrag{T'e}{$T'$}
\psfrag{T}{$T$}
\psfrag{Te}{$T$}
\psfrag{T''}{$B$}
\psfrag{t'1}{$\th'_1$}
\psfrag{t'2}{$\th'_2$}
\psfrag{t'3}{$\th'_3$}
\psfrag{t'3+e}{$\th'_3+\eta$}
\psfrag{d1}{$d_1$}
\psfrag{d2}{$d_2$}
\psfrag{t1}{$\th_1$}
\psfrag{t2}{$\th_2$}
\psfrag{t3}{$\th_3$}
\psfrag{t3+e}{$\th_3+\eta$}
\psfrag{1-t3-e}{$1-\th_3-\eta$}
\psfrag{e}{$e$}
\psfrag{e1}{$e_1$}
\psfrag{e2}{$e_2$}
\psfrag{a1}{$a_1$}
\psfrag{a2}{$a_2$}

\includegraphics[width=0.8\textwidth]{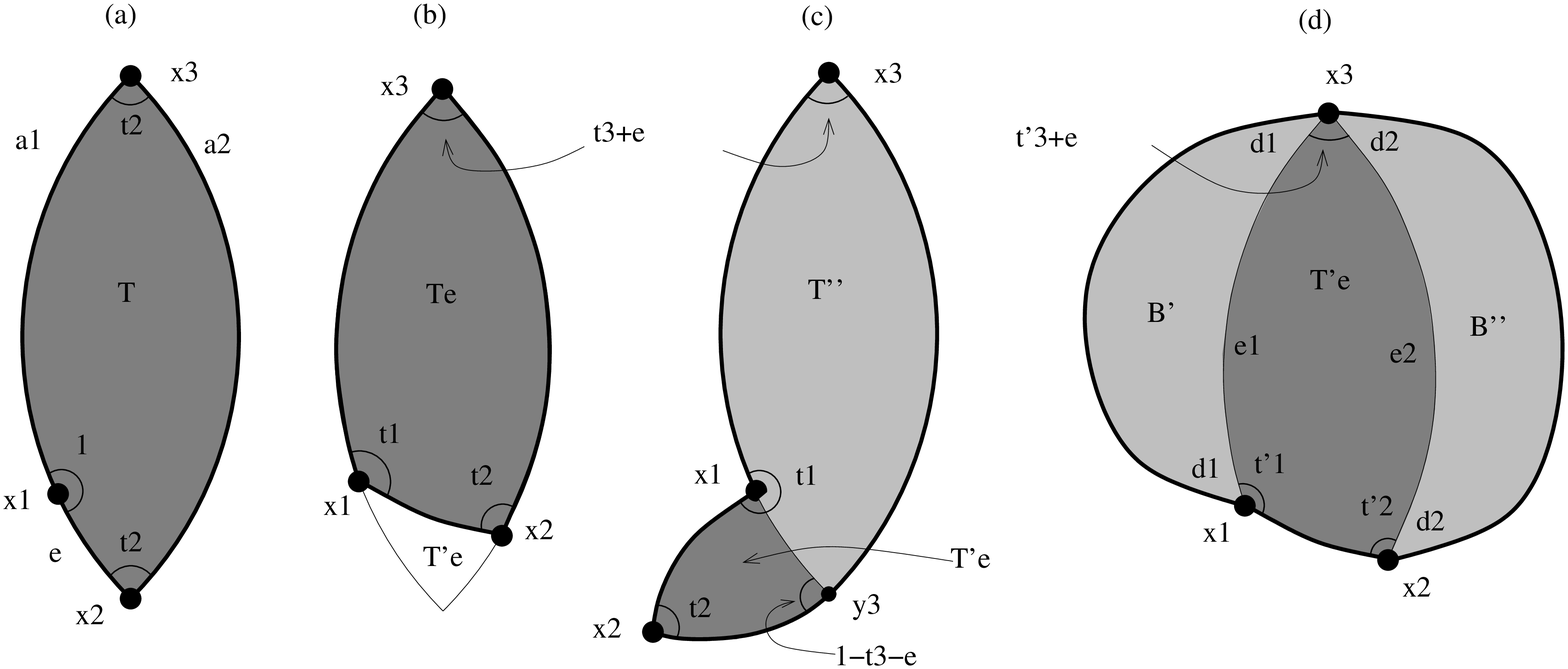}
\caption{{\small Triangles close to ordinary bigons.}}\label{fig:degenerate-tri-1}
\end{figurehere}
\end{center}

{\it{Case (a): $\th_1=1$.}}\\
Then we can take $T=B_{\th_2}$, mark the two vertices of the ordinary bigon by $x_2,x_3$ and place $x_1$ on $\partial T$ at distance $\pi-\frac{\varepsilon}{2}$ from $x_3$. For such a $T$, we have $\eta=0$ and both edges $x_1x_3$ and $x_2x_3$ have length
at most $\pi$.\\
Now, keep the edge $e$ from $x_1$ to $x_2$ fixed.
For $i=1,2$ and
for every $\nu_i$ close enough to $\th_i$
shoot the geodesic arc $a_i$ starting from $x_i$ and that forms
an angle $\pi\nu_i$ with $e$. Let $x_3$ be the first intersection of $a_1$ and $a_2$ and let $T_{(\nu_1,\nu_2)}$ the triangle bounded by $e,a_1,a_2$ and with internal angles $\pi(\nu_1,\nu_2,\theta_3)$, where $\theta_3$ is clearly a continuous function of $(\nu_1,\nu_2)$.
It is easy to see that $(\nu_1,\nu_2)\mapsto T_{(\nu_1,\nu_2)}$ is a continuous family of triangles
and that $T_{(\th_1,\th_2)}=T$. Hence, $T$ is $(x_1,x_2)$-deformable.

{\it{Case (b): $\th_1,\th_2<1$.}}\\
Let $\eta>0$ be smaller than
$2(1-\th_1),2(1-\th_2),1-\th_3$
so that the triple $(\th_1,\th_2,1-\th_3-\eta)\in (0,1)^3$ satisfies the triangular inequality.
By Lemma \ref{triangle}, there exists a convex triangle $T'$
with angles $\pi(1-\th_1,1-\th_2,\th_3+\eta)$ and vertices $(x_1,x_2,y_3)$. 
By construction, such a $T'$ is embedded inside an ordinary bigon $B_{\th_3+\eta}$ with vertices $y_3,x_3$.
The closure $T$ of the complement of $T'$ inside $B_{\th_3+\eta}$ is a triangle vertices $x_1,x_2,x_3$,
angles $\pi(\th_1,\th_2,\th_3+\eta)$ and $|x_1x_3|,|x_2x_3|<\pi$.
Notice that, as $\eta\rar 0$, the area of $T'$ (which depends on $\eta$) goes to zero
and so its diameter goes to zero too (because $\th_1,\th_2\in(0,1)$ are fixed).
Hence, for a sufficiently small $\eta$, the triangle $T$ is also $\pi(1-\varepsilon)$-wide at $x_3$.
Since the holonomy of $DT$ is clearly non-coaxial, it is angle-deformable and so in particular $T$ is $(x_1,x_2)$-angle-deformable.

{\it{Case (c): $\th_1\in(1,2)$, $\th_2\in(0,1)$ and $\th_3\in(0,1]$.}}\\
Let $\eta<0$ so that $|\eta|$ is smaller than $\varepsilon$,
$2(\th_1-1)$ and $2\th_2$.
Thus, the triple $(2-\th_1,1-\th_2,\th_3+\eta)\in (0,1)^3$ satisfies
the triangular inequality and by Lemma \ref{triangle} there exists a strictly convex triangle $T'$ with vertices $x_1,x_2,y_3$ and angles $\pi(\th_1-1,\th_2,1-\th_3-\eta)$. By construction, $|x_2 y_3|<\pi$. 
The triangle $T$ is the obtained by gluing $T'$ with
a standard bigon $B=B_{\th_3+\eta}$ with vertices $x_3,y'_3$ in such a way that $y'_3$ is identified to $y_3$ and $e_2$ is glued to a portion of an edge of $B_{\th_3+\eta}$.  Thus, $|x_1 x_3|<\pi$ and $|x_2 x_3|<2\pi$.
As before, it is clear that the length of $x_1 y_3$ goes to zero as
$\eta\rar 0$. Thus, $T$ is $\pi(1-\varepsilon)$-wide at $x_3$ for 
$|\eta|$ small enough.
As above, $DT$ is non-coaxial and so angle-deformable, hence $T$ is $(x_1,x_2)$-angle-deformable.

{\it{Case (d): $\th_3>1$.}}\\
Let $d_1,d_2$ be positive integers such that 
$\th'_1=\th_1-d_1\in (0,2)$, $\th'_2=\th_2-d_2\in (0,1]$ and
$\th'_3=\th_3-(d_1+d_2)\in (0,1]$.
By cases (a) or (b), there exists an $(x_1,x_2)$-angle deformable triangle $T'$
with angles $\pi(\th'_1,\th'_2,\th'_3+\eta)$ for some $|\eta|<\varepsilon$,
which is $\pi(1-\varepsilon)$-wide at $x_3$. Call $e_1,e_2$ the edges 
$x_1x_3$ and $x_2x_3$ of
$T'$, of lengths $\ell_1,\ell_2\in \pi(1-\varepsilon,2)$.
The triangle $T$ is then obtained by gluing an edge of the exceptional bigon $B(d_1,\ell_1)$ with $e_1$ and an edge of $B(d_2,\ell_2)$ with $e_2$.
Because $B(d_1,\ell_1)$ and $B(d_2,\ell_2)$ are $\pi(1-\varepsilon)$-wide at their vertices, such a $T$ is $\pi(1-\varepsilon)$-wide at $x_3$.
Since this gluing procedure can be performed in families, the obtained
triangle $T$ is $(x_1,x_2)$-angle-deformable.
%
%
%
%
%
%
%
%
%
%
%
\end{proof}

\begin{proposition}[Triangles close to a double bigon]\label{differencetriangle}
Let $\th_1,\th_2>0$ with $\th_3=\th_1-\th_2-1\geq 0$ and assume that $\th_2$ is not an integer.
Then for every $\varepsilon>0$ there exist $\eta\in(-\varepsilon,\varepsilon)$ and 
a spherical triangle $T$ 
with angles $\pi(\th_1, \th_2,\th_3+\eta)$ and vertices $x_1,x_2,x_3$,
which is $\pi(1-\varepsilon)$-wide at $x_3$
and $(x_1,x_2)$-angle-deformable.
\end{proposition}
\begin{proof}
Again we divide the proof in four cases, illustrated in Figure \ref{fig:degenerate-tri-2}.

\begin{center}
\begin{figurehere}
\psfrag{1}{$1$}
\psfrag{(a)}{(a)}
\psfrag{(b)}{(b)}
\psfrag{(c)}{(c)}
\psfrag{(d)}{(d)}
\psfrag{x1'}{$x'_1$}
\psfrag{x2'}{$x'_2$}
\psfrag{x3'}{$x'_3$}
\psfrag{x1}{$x_1$}
\psfrag{x2}{$x_2$}
\psfrag{x3}{$x_3$}
\psfrag{y3}{$y_3$}
\psfrag{x4}{$x_4$}
\psfrag{B'}{$B'$}
\psfrag{B'''}{$B'''$}
\psfrag{T}{$T$}
\psfrag{T'}{$T'$}
\psfrag{T''}{$T''$}
\psfrag{t'1}{$\th'_1$}
\psfrag{t'2}{$\th'_2$}
\psfrag{t'3}{$\th'_3$}
\psfrag{t'3+e}{$\th'_3+\eta$}
\psfrag{d1}{$d_1$}
\psfrag{d3}{$d_3$}
\psfrag{t1}{$\th_1$}
\psfrag{t2}{$\th_2$}
\psfrag{t3}{$\th_3$}
\psfrag{y2}{$y_2$}
\psfrag{t4}{$\th_4$}
\psfrag{1-t4}{$1-\th_4$}
\psfrag{t3+e}{$\th_3+\eta$}
\psfrag{1-t2}{$\small{1-\th_2}$}
\psfrag{t1-2}{$\th_1-2$}
\psfrag{2}{$2$}
\psfrag{e}{$e$}
\psfrag{e1}{$e_1$}
\psfrag{e2}{$e_2$}
\psfrag{a1}{$a_1$}
\psfrag{a2}{$a_2$}

\includegraphics[width=0.7\textwidth]{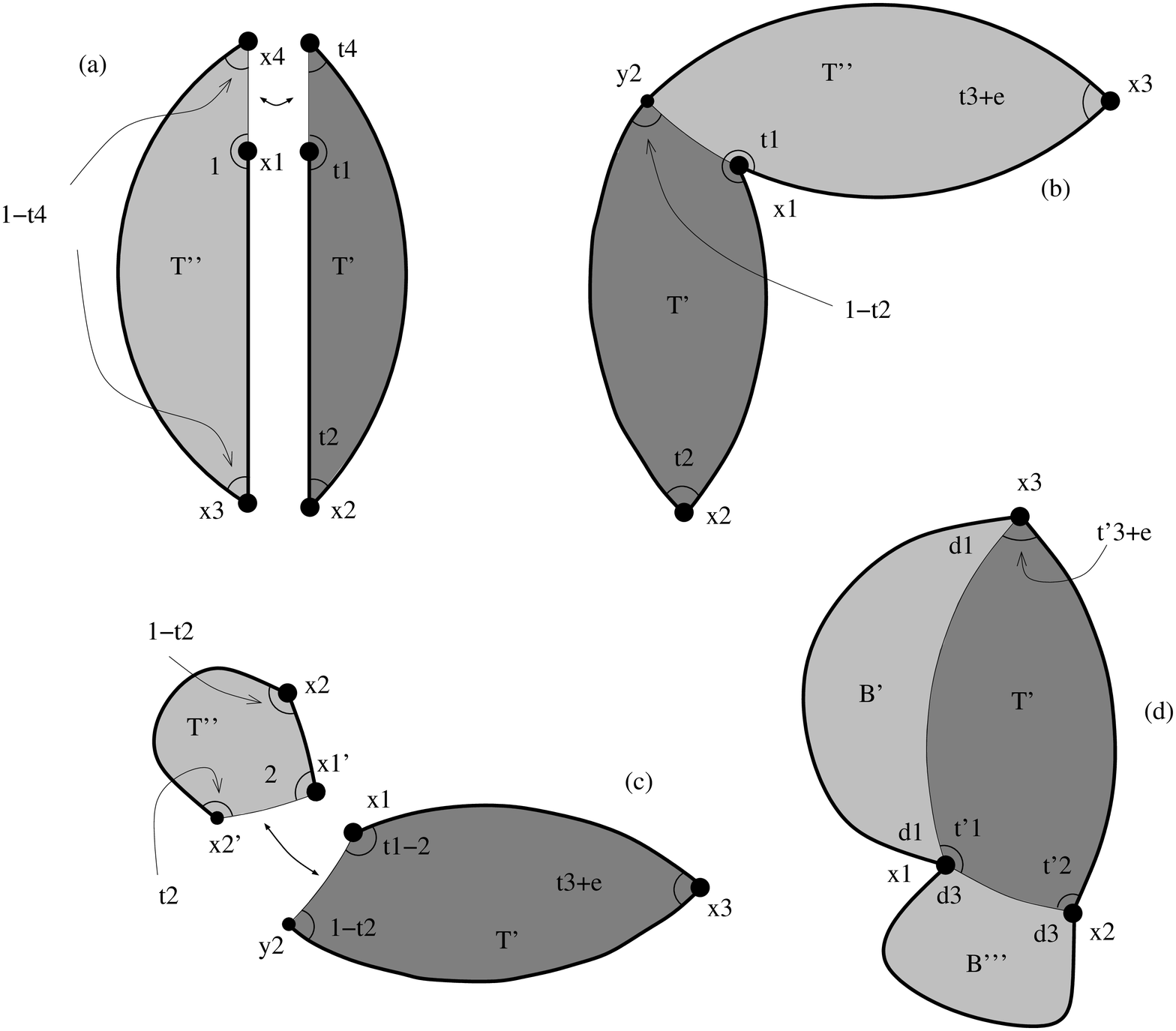}
\caption{{\small Triangles close to double bigons.}}\label{fig:degenerate-tri-2}
\end{figurehere}
\end{center}

{\it{Case (a): $\th_2\in(0,1)$, $\th_1=2$.}}\\
In this case $\th_3=1-\th_2>0$.
Then $T$ can be chosen to be the triangle $T(1,\pi(1-\varepsilon/2),\th_2)$
constructed in Lemma \ref{twoconstructions}.
Notice that, in this case, $\eta=0$ and that both edges $x_1 x_2$
and $x_1 x_3$ are shorter than $\pi$.
To see that such triangle is $(x_1,x_2)$-angle-deformable, label the midpoint of the edge $x_2x_3$ by $x_4$, so that the segment $x_1x_4$ splits $T$ into the triangles $T'$ with vertices $(x_1,x_2,x_4)$ and angles $\pi(1,\th_2,\th_2)$ and $T''$
with vertices $(x_1,x_4,x_3)$ and angles $\pi(1,\th_3,\th_3)$.
Also, the edge $x_1 x_4$ has length $\pi\varepsilon/2$.
By Proposition \ref{sumtriangle}, the triangle $T'$ is $(x_1,x_2)$-angle-deformable and so there exists a neighbourhood $\Nang'\subset\RR^2$ of
$(\th_1,\th_2)$, a continuous $\theta_4:\Nang'\rar\RR$ with
$\theta_4(\th_1,\th_2)=\th_2$ and a continuous family of triangles
$\Nang'\ni \vc{\nu}\rar T'_{\vc{\nu}}$ such that $T'_{(\th_1,\th_2)}=T'$
and $T_{\vc{\nu}}$ has angles $\pi(\nu_1,\nu_2,\theta_4(\nu_1,\nu_2))$.
Clearly, the length $\ell_{\vc{\nu}}$ of the edge $x_1x_4$ of $T'_{\vc{\nu}}$ depends continuously on $\vc{\nu}$.
Consider now the continuous family $\vc{\nu}\mapsto T''_{\vc{\nu}}$ of ordinary bigons with opposite vertices $(x_3,x'_4)$ and angles $\pi(\theta_4(\vc{\nu}),1-\theta_4(\vc{\nu}))$ and label by $x'_1$ a point
of the edge $x_3 x'_4$ that sits at distance $\ell_{\vc{\nu}}$ from $x'_4$. Gluing $T'_{\vc{\nu}}$ and $T''_{\vc{\nu}}$ along the
segments $x_1 x_4$ and $x'_1 x'_4$, we obtain the wished family of triangles parametrized by $\Nang'$.

{\it{Case (b): $\th_2\in(0,1)$ and $\th_1\in(1,2)$.}}\\
We proceed as in the proof of Proposition \ref{sumtriangle}, case (b).
Pick $\eta>0$ and smaller than $\varepsilon,2(2-\th_1),2\th_2$.
Then the triple $(2-\th_1,\th_2,1-(\th_3+\eta))$ satisfies the triangular
inequality and so there exists a convex triangle $T'$ with
vertices $x_1,y_2,x_3$ and angles $\pi(\th_1-1,1-\th_2,\th_3+\eta)$.
Moreover, $\eta$ can be chosen small enough so that such $T'$ is $\pi(1-\varepsilon)$-wide at $x_3$. Clearly,
the edge $y_2 x_1$ of $T'$ is shorter than $\pi$.
The desired triangle $T$ is then obtained by gluing
an ordinary bigon $B_{\th_2}$ with vertices $x_2$ and $y'_2$
to $T'$ by identifying $y_2$ to $y'_2$ and $y_2 x_1$ to a portion of an edge of $B_{\th_2}$. We underline that both $x_1 x_2$ and $x_1 x_3$ are shorter than $\pi$.
The double of such a triangle has non-coaxial holonomy and so the triangle is angle-deformable, and in particular $(x_1,x_2)$-angle-deformable.

{\it{Case (c): $\th_2\in(0,1)$ and $\th_1\in(2,3)$.}}\\
Pick $\eta<0$ such that $|\eta|$ is smaller than $\varepsilon,2(1-\th_2),\th_1-2$.
Then the triple $(3-\th_1,1-\th_2,1-(\th_3+\eta))$ satisfies the triangular
inequality and so there exists a convex triangle $T'$ with
vertices $x_1,y_2,x_3$ and angles $\pi(\th_1-2,\th_2,\th_3+\eta)$.
Moreover, $\eta$ can be chosen small enough so that such $T'$ is $\pi(1-\varepsilon)$-wide at $x_3$. Clearly,
the edge $x_1 y_2$ of $T'$ has length $\ell<\pi$.
Consider now a triangle $T''=T(1,\ell,\th_2)$ with vertices
$x_2$ of angle $\pi\th_2$,
$x'_2$ of angle $\pi(1-\th_2)$ and $x'_1$ of angle $2\pi$
and edges incident at $x'_1$ of length $l$.
The desired triangle $T$ is then obtained by gluing
$T''$ with $T'$ by identifying the edge
$x'_1 x'_2$ of the former to the edge $x_1 y_2$ of the latter.
As in the previous case, $x_1 x_2$ and $x_1x_3$ are shorter than $\pi$ and
the triangle is angle-deformable, and in particular $(x_1,x_2)$-angle-deformable.

{\it{Case (d): $\th_2$ not an integer.}}\\
Let $d_1,d_3$ be positive integers such that 
$\th'_2=\th_2-d_3\in (0,1)$ and
$\th'_3=\th_3-d_1\in (0,2)$ and so
$\th'_1=\th_1-d_1-d_3\in (1,3)$.
The previous cases ensure that there exists an $(x_1,x_2)$-angle-deformable triangle $T'$
with angles $\pi(\th'_1,\th'_2,\th'_3+\eta)$ for some $|\eta|<\varepsilon$,
which is $(\pi-\varepsilon)$-wide at $x_3$
and such that $|x_1 x_2|<\pi$ and $|x_1 x_3|<\pi$. 
The triangle $T$ is then obtained by gluing an edge of the exceptional bigon $B(d_1,|x_1 x_3|)$ with $x_1 x_3$ and an edge of $B(d_3,|x_1 x_2|)$ with $x_1 x_2$.
Such a $T$ is clearly $\pi(1-\varepsilon)$-wide at $x_3$.
Since this gluing procedure can be performed in families, the $(x_1,x_2)$-angle-deformability of $T$ follows from the analogous property of $T'$.
\end{proof}

\begin{remark}
The restriction $\th_2\notin\ZZ$ is not due to the chosen proof.
In fact, for $\th_2=1$ and $\th_1\notin\ZZ$, we would be looking for a bigon with different angles which are not multiples of $\pi$ and it is known that such bigons do not exists.
As another example, if $\th_2=d$, $\th_3=d'$ and $\th_1=d+d'+2$, 
the double $DT$ a triangle $T$ would be a (connected) ramified cover of $\Sph$ over $3$ points and this is clearly impossible, as
the product of a $d$-cycle and a $d'$-cycle in a group of permutations
cannot give a $(d+d'+2)$-cycle. 
\end{remark}

\subsection{Cut-and-paste operations}\label{sec:cut-paste}

We recall that, for every $\alpha>0$, we denoted by
$S_\alpha$ the spherical surface homeomorphic
to $\Sph$ with two cone points of angle $2\pi\alpha$ sitting
at distance $\pi$, as in Corollary \ref{cor:two-pointed}(a).
 

\subsubsection{Cut-and-paste at a conical point}

The goal of this section is to describe a cut-and-paste procedure
that permits to increase the number of conical points on a sphere
by modifying the metric in a neighbourhood of a conical point.

We remind that, if $S$ is a closed spherical surface which is $r$-wide at
a conical point $y$, then $U_y(r)$ denotes the closed subsurface of $S$
obtained by removing the closed ball of radius $r$ centered at $y$.
Thus, $\partial U_y(r)$ is a circle of length $2\pi\alpha \sin(r)$, where $2\pi\alpha$ is the angle at $y$.

The following lemma is obvious: the situation is illustrated in Figure \ref{fig:gluing-two-beans}.

\begin{lemma}[Surgery at conical points]\label{lemma:gluing-conical}
Let $S$ and $S'$ be two spherical surfaces with cone points
$y\in S$ and $y'\in S'$ of angles $2\pi\alpha$.
Suppose that $S$ is $r$-wide at $y$ and $S'$ is $(\pi-r)$-wide at $y'$
for some $r\in(0,\pi)$. Then the surface $S\#_r S'$ obtained
by gluing $U=U_y(r)$ and $U'=U_{y'}(\pi-r)$ through an isometry
$\partial U\cong\partial U'$ 
is a spherical surface with conical points.
Moreover, if $S$ or $S'$ has non-coaxial holonomy, the
same holds for $S\#_r S'$.
%
\end{lemma}

\begin{center}
\begin{figurehere}
\psfrag{S}{$S$}
\psfrag{S'}{$S'$}
\psfrag{y}{$y$}
\psfrag{y'}{$y'$}
\psfrag{r}{$r$}
\psfrag{pi-r}{$\pi-r$}
\psfrag{a}{$\alpha$}
\psfrag{U}{$U$}
\psfrag{U'}{$U'$}
\includegraphics[width=0.4\textwidth]{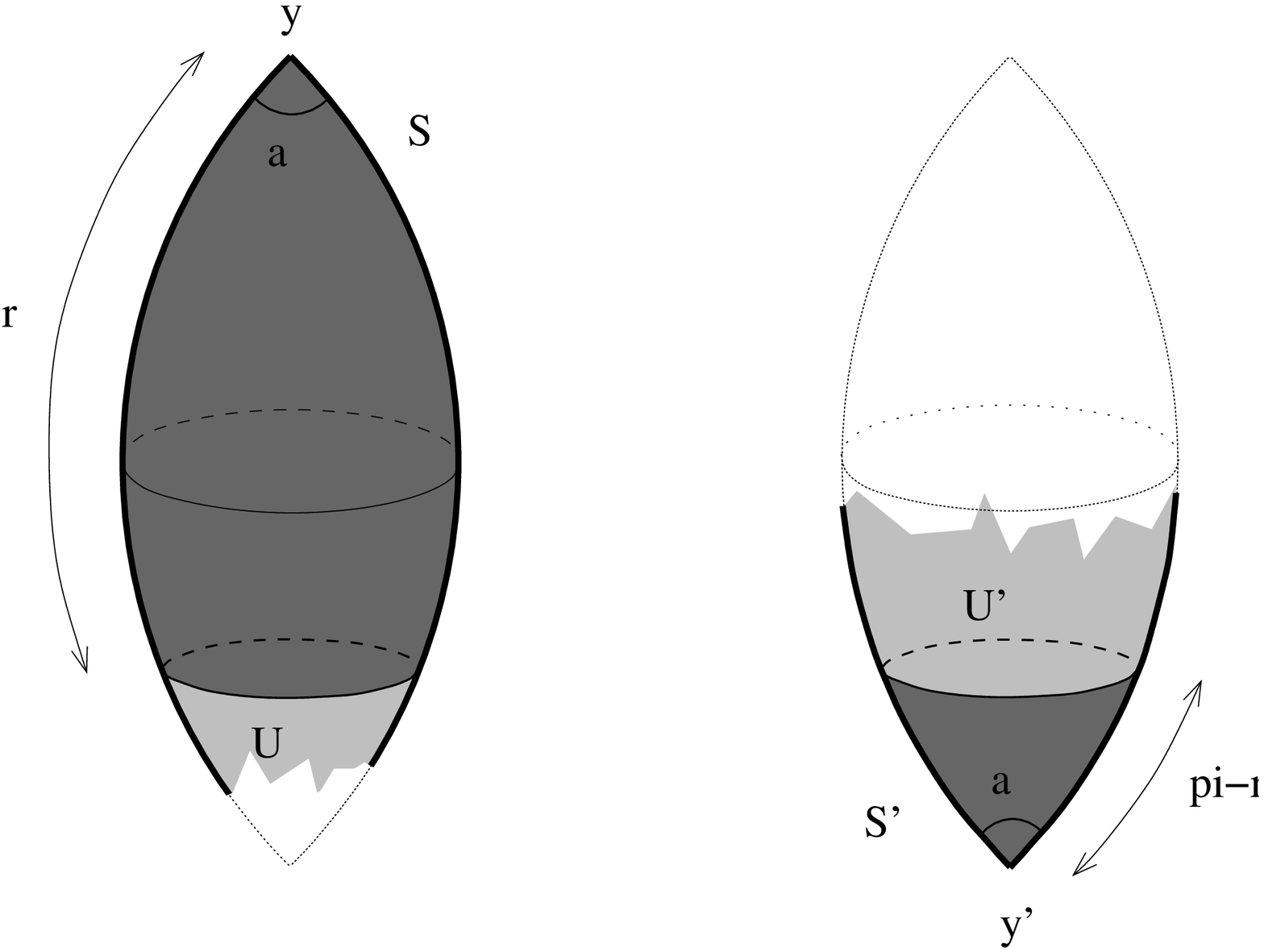}
\caption{{\small Chopping off the dark region and gluing $U$ and $U'$.}}\label{fig:gluing-two-beans}
\end{figurehere}
\end{center}

\subsubsection{Cut-and-paste along a path}


We recall that a path $\gamma$ on a surface $S$ is called {\it{simple}} if it is injective (i.e. if its image has no self-intersections).

\begin{definition}
A path $\gamma$ on a spherical surface $S$ is
{\it{simply developable}} if its developing map $\dev_\gamma$
is injective.
\end{definition}

The following lemma is obvious.


\begin{lemma}[Surgery along a path I]\label{gluing}
Let $S$ and $S'$ be two spherical surfaces.
Let $\gamma$ (resp. $\gamma'$) be a simple path on $S$ (resp. $S'$) running
from the conical point $y_1$ of angle $2\pi\alpha_1$ to the conical point $y_2$ of angle $2\pi\alpha_2$
(resp. from the conical point $y'_1$ of angle $2\pi\alpha'_1$ to the conical point $y'_2$ of angle $2\pi\alpha'_2$)
and intersecting the singularities nowhere else.
Suppose that $\gamma$ and $\gamma'$ are isometric.
Then the surface denoted by $S_\gamma\#_{\gamma'}S'$ 
and obtained by gluing $S\setminus\gamma$ and $S'\setminus\gamma'$ via the isometric identification of $\gamma$
with $\gamma'$ is a spherical surface;
moreover, 
the two points
$y_i$ and $y'_i$ are identified to a conical point 
of angle $2\pi(\alpha_i+\alpha'_i)$ on
$S_\gamma\#_{\gamma'}S'$
for $i=1,2$.
\end{lemma}

Using this lemma we get the following result. The situation is illustrated in Figure \ref{fig:gluing-path-1}.

\begin{proposition}[Surgery along a path II]\label{adding4piand2pi2pi}
Consider a spherical surface $S$ with conical points $y_1,\dots,y_k$
of angles $2\pi\beta_1,\dots,2\pi\beta_k$
and let $\gamma$ be a simple and simply developable path on $S$
that joins $y_1$ and $y_2$. Let also $d\in \ZZ_+$.
\begin{itemize}
\item[(a)]
The spherical surface obtained by gluing $S\setminus\gamma$
and $d$ copies of $\Sph\setminus\dev_\gamma$ via an isometric identification
of their boundaries has conical singularities $z_1,\dots,z_k$
of angles $2\pi(\beta_1+d,\beta_2+d,\beta_3,\dots,\beta_k)$.
%
%
\item[(b)]
Suppose  $\beta_2<1$ and that $\gamma$ is geodesic path
of length $\ell=|\gamma|<\pi$.
Then there exists a spherical surface $S'$ and a path $\gamma'$
on $S'$ isometric to $\gamma$ such that
$S_\gamma\#_{\gamma'}S'$
%
has conical singularities $z_1,\dots,z_k$ of angles
$2\pi(\beta_1+2d, \beta_2,\dots,\beta_k)$.
Moreover, 
the conical points $z_1$ and $z_2$
on $S_\gamma\#_{\gamma'}S'$
are joined by
a geodesic arc of length $\ell$.
%
\end{itemize}
Moreover, if $S$ is deformable (resp. non-coaxial), so are the constructed surfaces.
\end{proposition}

\medskip

\begin{center}
\begin{figurehere}
\psfrag{S-g}{$S\setminus\gamma$}
\psfrag{S'-g'}{$S'\setminus\gamma'$}
\psfrag{S'-g3}{$S'\setminus\gamma'_3$}
\psfrag{y1}{\parbox{0.2cm}{\[\begin{array}{c} {[y_1]}\\ || \\ z_1\end{array}\]}}
\psfrag{y2}{\parbox{0.2cm}{\[\begin{array}{c}z_2\\ || \\ {[y_2]}\end{array}\]}}
\psfrag{y2n}{$[y_2]$}
\psfrag{y'2}{$z_2=[y'_2]$}
\psfrag{b1}{$\beta_1$}
\psfrag{b2}{$\beta_2$}
\psfrag{1-b2}{$1-\beta_2$}
\psfrag{2d}{$2d$}
\psfrag{d}{$d$}
\psfrag{(a)}{(a)}
\psfrag{(b)}{(b)}
\includegraphics[width=0.7\textwidth]{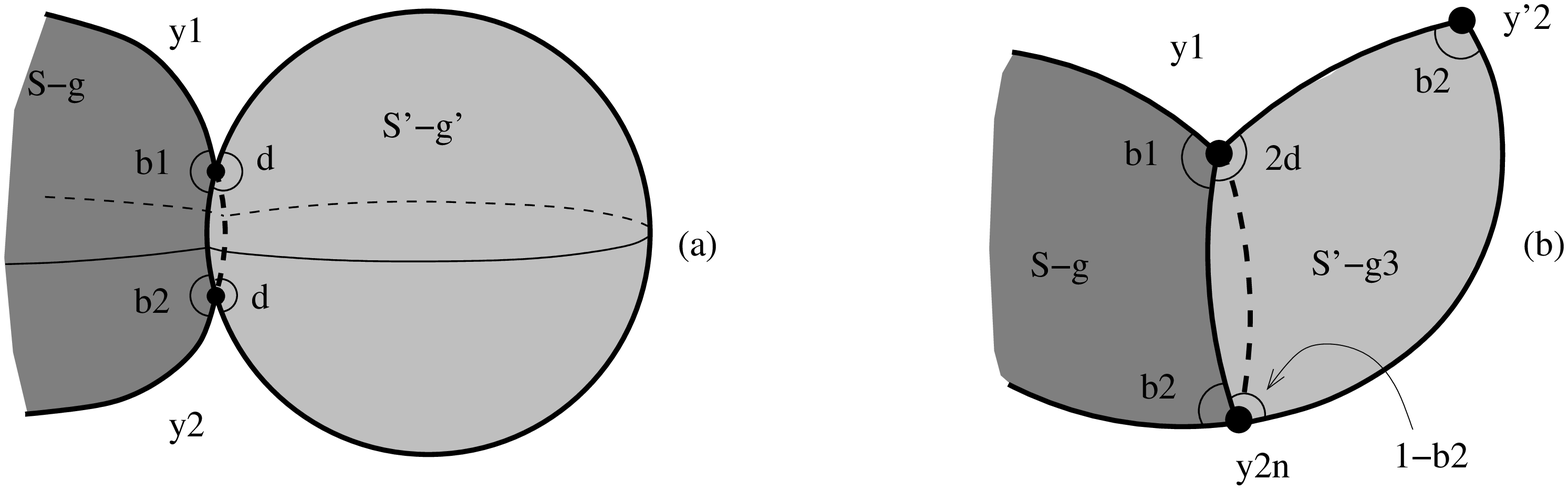}
\caption{{\small Increasing the angles by $d(\vc{e}_1+\vc{e}_2)$ and by $2d\vc{e}_1$.}}\label{fig:gluing-path-1}
\end{figurehere}
\end{center}

\begin{proof}[Proof of Proposition \ref{adding4piand2pi2pi}]
About (a),  according to Lemma \ref{gluing} it is
sufficient to construct a spherical surface $S'$ with two conical points of angles $2\pi d$
that are joined by a simple curve $\gamma'$ isometric to $\gamma$.
The conical points on the new spherical surface $S_\gamma\#_{\gamma'}S'$ will be the classes of the $y_i$, which will be denoted by $z_i$.

To do this, consider the developing map $\dev_\gamma$ to $\Sph$,
which is injective, and call $y'_1$ and $y'_2$ its starting and end points.
Then let $S'$ be the ramified cover of $\Sph$ of degree $d$ branched at $y'_1,y'_2$ and let $\gamma'$ be one of the $d$ lifts of $\dev_\gamma$. Indeed, such an $S'$ is obtained by gluing $d$ copies
of $\Sph\setminus\dev_\gamma$.

Observe that the path $\gamma$ is still simply developable for any
other spherical metric on $S$ sufficiently close to the given one.
Thus, if $S$ is deformable, so is the constructed $S_\gamma\#_{\gamma'}S'$.
%
%

Concerning (b), let $S'$ be the double of the spherical triangle $T(d,\ell,\beta_2)$
constructed in Lemma \ref{twoconstructions}. This is a sphere
with conical points $y'_1,y'_2,y'_3$ of angles $2\pi (2d)$, $2\pi \beta_2$ and $2\pi (1-\beta_2)$ and $y'_1$ is joined to $y'_2$ and $y'_3$ by
two geodesics of length $\ell=|\gamma|$: call them $\gamma'_2$ and $\gamma'_3$.
The new spherical surface is obtained by gluing $S\setminus \gamma$
to $S'\setminus \gamma'_3$ by identifying $\gamma$ to $\gamma'_3$,
$y_1$ to $y'_1$ and $y_2$ to $y'_3$.
Observe that $y_1$ and $y'_1$ merge to a conical point $z_1$ of angle $2\pi(\beta_1+2d)$, but $y_2$ and $y'_3$ merge to a regular point
(i.e. a point of angle $2\pi$). On the other hand,
the points $y'_2,y_3,\dots,y_k$ will give rise to conical points $z_2,\dots,z_k$
of angles $2\pi\beta_2,\dots,2\pi\beta_k$.
Finally, notice that $\gamma'_2$ descends on $S_\gamma\#_{\gamma'_3}S'$
to a geodesic path of length $\ell$ between $z_1$ and $z_2$.

Observe that, for any spherical metric on $S$ sufficiently close to the given one, the path $\gamma$ continuously deforms to a
geodesic between $y_1$ and $y_2$ of length less than $\pi$.
Thus, if $S$ is deformable, so is the constructed $S_\gamma\#_{\gamma'}S'$.

Moreover, considering $S\setminus\gamma$ inside $S_\gamma\#_{\gamma'}S'$, it is easy to see that in both cases (a) and (b) the spherical surface $S_\gamma\#_{\gamma'}S'$ has the same holonomy as $S$, and so it is non-coaxial if and only if $S$ is.
%
\end{proof}

\subsection{Spheres with four conical points}\label{sec:four-points}

In this section we will prove Theorem \ref{existence} for 
spheres with four conical points of angles not divisible by $2\pi$.

\begin{theorem}[Existence of $4$-punctured spherical metrics with non-integral angles]\label{thm:four}
Let $\th_1,\th_2,\th_3,\th_4$ be real non-integer numbers that satisfy
both the positivity constraints (\ref{gauss-bonnet}) and
the holonomy constraints (\ref{theinquality}) strictly.
Then there exists a sphere $S$ endowed with a spherical metric 
with $4$ conical singularities of angles $2\pi\th_1,\dots,2\pi\th_4$
and non-coaxial holonomy.
\end{theorem}

\begin{remark} By Luo's Theorem \ref{luo}, all metrics from Theorem \ref{thm:four} are deformable.
\end{remark}

The proof proceeds in two steps. First we study several types of spherical quadrilaterals embedded and immersed in $\Sph$. We construct 
embedded  quadrilaterals 
that have at most two angles larger than $\pi$ and immersed quadrilaterals with 
three angles less that $\pi$ and one angle in the interval $(2\pi, 3\pi)$.
By doubling such quadrilaterals we obtain all spherical metrics with non-integral
angles $2\pi\cdot (\th_1,\th_2,\th_3,\th_4)$ 
with $\vc{\th}\in (0,2)^2\times(0,1)^2$ or $\vc{\th}\in (2,3)\times(0,1)^3$,
apart from metrics with two exceptional one parameter families of angles.
Spherical surfaces with four conical points in the
exceptional classes are obtained by an alternative construction.
Finally, using cut-and-paste operations along paths we get all the remaining metrics.\\


Let $\mathfrak{S}_4$ be the group of permutations of $\{1,2,3,4\}$ and view
it as acting on $\RR^4$ in the obvious way $\mathfrak{S}_4\times\RR^4 \ni
(\sigma,\vc{\th})\mapsto \vc{\th}_\sigma\in\RR^4$, where $\vc{\th}_\sigma:=(\th_{\sigma(1)},\th_{\sigma(2)},\th_{\sigma(3)},\th_{\sigma(4)})$.
Let $\mathfrak{D}_8$ be the subgroup of $\mathfrak{S}_4$ generated by $(1234)$ and $(13)$, which is isomorphic to a dihedral group of order $8$.
The following simple observation will be useful.

\begin{notation}
The four vertices of a quadrilateral $Q$ are always cyclically labelled respecting
an orientation on $\pa Q$.
\end{notation}

\begin{lemma}[Allowed permutations]\label{lemma:permutations}
Suppose that there exists a spherical quadrilateral $Q$ with vertices $x_1,\dots,x_4$ and conical points of angles $\pi\cdot\vc{\th}$. 
Then for every $\sigma\in\mathfrak{D}_8$ there exists a spherical quadrilateral
$Q'$ with vertices $x'_1,\dots,x'_4$ and angles $\pi\cdot\vc{\th}_\sigma$.
\end{lemma}
\begin{proof}
We can produce $Q'$ out of $Q$ just cyclically permuting the labels or switching the orientation. In the former case
we easily see that this corresponds to the permutation $\sigma_1=(1234)$ or $\sigma_1=(4321)$;
in the latter case, this corresponds to one of the following
$\sigma_2=(12)(34)$, $\sigma_2=(13)(24)$ or $\sigma_2=(14)(23)$.
Since $\{\sigma_1,\sigma_2\}$ generates the $\mathfrak{D}_8$,
the conclusion follows.
\end{proof}

\begin{remark}
Given a surface $S$ of genus $0$ with $4$ conical points of angles $2\pi\cdot\vc{\th}$,
we can clearly produce an $S'$ with angles $2\pi\cdot\vc{\th}_\sigma$ for every
$\sigma\in\mathfrak{S}_4$: indeed, it is enough to relabel the conical points accordingly to $\sigma$. On the other hand, given a spherical quadrilateral $Q$ with
angles $\pi\cdot\vc{\th}$, it is not always possible to produce a quadrilateral $Q'$
with angles $\pi\cdot\vc{\th}_\sigma$ with $\sigma\in\mathfrak{S}_4$ but $\sigma\notin\mathfrak{D}_8$.
\end{remark}

\subsubsection{Convex quadrilaterals}
Let $\vc{c}=(c_1,c_2,c_3,c_4)\in\RR^4$ be a vector with strictly half-integral coordinates.
Recall that we denote by $\cube_{\vc{c}}$ the unit cube in $\RR^4$
with centre $\vc{c}$ and by $\tcube_{\vc{c}}$ the corresponding truncated cube. 
Notice that, since $n=4$ is even, $\vc{m}\in\ZZ_o^4$ if and only if $\vc{m}-\1\in\ZZ^4_o$:
thus, $\vc{\th}\in\Hang^4$ if and only if $\vc{\delta}=\vc{\th}-\1\in\Hang^4$.

\begin{definition}
Let $\vc{c}\in\RR^4$ be a strictly half-integral vector and let $\vc{m}$ be an even
integral vertex of $\cube_{\vc{c}}$. The {\it{half truncated cube}} centered at $\vc{c}$
associated to the vertex $\vc{m}$ is
\[
\hcube_{\vc{c}}(\vc{m}):=\{\vc{p}\in\tcube_{\vc{c}}\,|\,d_1(\vc{m},\vc{p})\leq 2\}\, .
\]
\end{definition}

\begin{example}\label{halfcube}
Let $\vc{c_0}=(\frac{1}{2}, \frac{1}{2}, \frac{1}{2}, \frac{1}{2})$ and
let $\1=(1,1,1,1)$ and $\0=(0,0,0,0)$ be even vertices of $\cube_{\vc{c_0}}$.
The truncated cube $\tcube_{\vc{c_0}}$ is the union of $\hcube_{\vc{c_0}}(\1)$ and
$\hcube_{\vc{c_0}}(\0)$; moreover, the two half truncated cubes only overlap along a face.
%
%
%
%
\end{example}

\begin{center}
\begin{figurehere}
\psfrag{Q}{$Q$}
\psfrag{t1+t2-1+t}{$\th_1+\th_2-1+t$}
\psfrag{t1}{$\th_1$}
\psfrag{y=y'}{$y=y'$}
\psfrag{t2}{$\th_2$}
\psfrag{t3}{$\th_3$}
\psfrag{t4}{$\th_4$}
\psfrag{x1}{$x_1$}
\psfrag{x2}{$x_2$}
\psfrag{x3}{$x_3$}
\psfrag{x4}{$x_4$}
\psfrag{D'2}{$\Delta'_t$}
\includegraphics[width=0.3\textwidth]{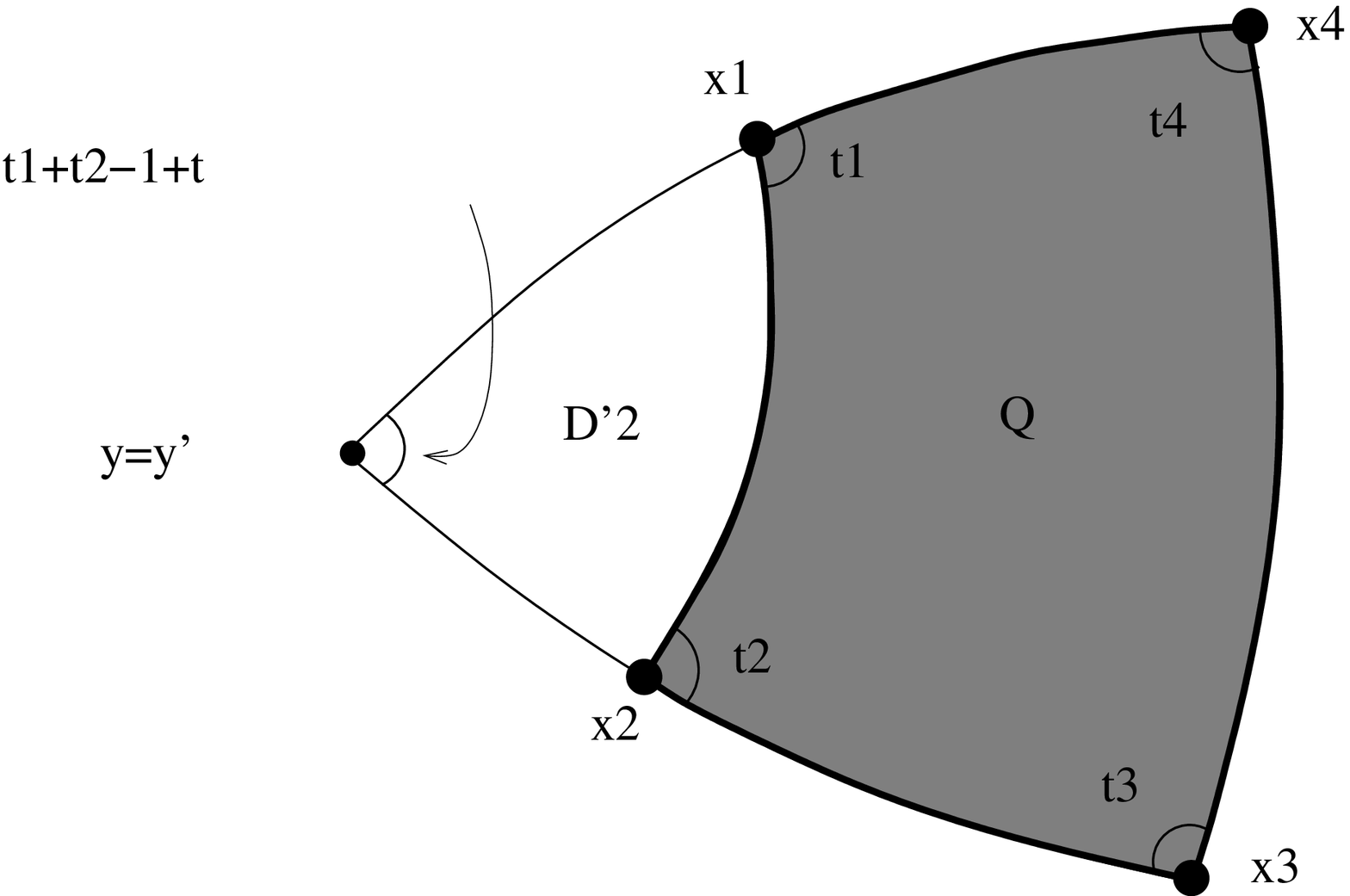}
\caption{{\small Construction of a convex quadrilateral.}}\label{fig:convex-quadrilateral}
\end{figurehere}
\end{center}

\begin{lemma}[Convex quadrilaterals]\label{convexquad}
For every $\vc{\th}=(\th_1,\dots,\th_4)$ in the interior of $\hcube_{\vc{c_0}}(\1)$
there exists a convex quadrilateral $Q$
with angles $\pi\cdot\vc{\th}$.
\end{lemma}

\begin{proof}
After possibly reversing the order and cyclically permuting
$\th_1,\dots,\th_4$, we can assume that $\th_1\ge \th_2,\th_3,\th_4$ and $\th_2\ge \th_4$. In this case one can check that the triple $(\th_1+\th_2-1,\th_3,\th_4)$ satisfies strictly both constraints (\ref{gauss-bonnet}) and (\ref{theinquality}). Hence, for $t\geq 0$ small enough, there exists a continuous family of spherical triangles $t\mapsto \Delta_t$ 
with vertices $y,x_3,x_4$
and
angles $\pi(\th_1+\th_2-1+t,\th_3,\th_4)$. For $t>0$ small enough there exists as well a continuous family $t\mapsto \Delta'_t$ of spherical triangles with vertices $y',x_2,x_1$
and
angles $\pi(\th_1+\th_2-1+t,1-\th_2,1-\th_1)$. 
Notice that the diameter  of $\Delta'_t\rar 0$ as $t\rar 0$.
Thus, for $t>0$ small enough, it is possible to inscribe $\Delta'_t$ inside
$\Delta_t$ in such a way that $y'$ coincides with $y$ and
that $y'x_2$ and $y'x_1$ are contained inside $yx_3$ and $yx_4$ respectively (see Figure \ref{fig:convex-quadrilateral}).
Hence, for such small $t>0$, we can obtain our
desired quadrilateral with vertices $x_1,x_2,x_3,x_4$
as the completion of $\Delta_t\setminus\Delta'_t$.
%
\end{proof}

\subsubsection{Non-convex quadrilaterals embedded in $\Sph$}

\begin{lemma}[Seven basic non-convex quadrilaterals]\label{1big3small}
Let $\vc{\th}\in\hcube_{\vc{c_0}}(\1)$ and consider the following table.
\[
\def\arraystretch{1.2}
\begin{array}{|c|c|c|c|}
\hline
i & f_i(\vc{\th}) & \vc{m_i} & \vc{c_i}\\
\hline
0 & (\th_1,\th_2,\th_3,\th_4) & (1,1,1,1) &  (\frac{1}{2},\frac{1}{2},\frac{1}{2},\frac{1}{2})\\
1 & (2-\th_1, 1-\th_2, \th_3, 1-\th_4) & (1,0,1,0) & (\frac{3}{2},\frac{1}{2},\frac{1}{2},\frac{1}{2})\\
2 & (\th_1+1, 1-\th_2, \th_3, \th_4) & (2,0,1,1) & (\frac{3}{2},\frac{1}{2},\frac{1}{2},\frac{1}{2})\\
3 & (\th_1+1, \th_2, \th_3, \th_4+1) & (2,1,1,2) & (\frac{3}{2},\frac{1}{2},\frac{1}{2},\frac{3}{2})\\
4 & (2-\th_1, 2-\th_4, 1-\th_3, 1-\th_2) & (1,1,0,0) & (\frac{3}{2},\frac{3}{2},\frac{1}{2},\frac{1}{2})\\
5 & (2-\th_1, \th_4, 2-\th_3, \th_2) &  (1,1,1,1) & (\frac{3}{2},\frac{1}{2},\frac{3}{2},\frac{1}{2})\\
6 & (1+\th_1,1-\th_2, 1+\th_3, 1-\th_4) & (2,0,2,0) &  (\frac{3}{2},\frac{1}{2},\frac{3}{2},\frac{1}{2})\\
7 & (1+\th_1, 1-\th_2, 2-\th_3, \th_4) & (2,0,1,1) & (\frac{3}{2},\frac{1}{2},\frac{3}{2},\frac{1}{2})\\
\hline
\end{array}
\]
For every convex quadrilateral $Q$ with cyclically ordered angles $\pi\cdot\vc{\th}$
and for every $1\leq i\leq 7$
there exists a quadrilateral $Q_i$ embedded in $\Sph$ with cyclically ordered angles $\pi\cdot f_i(\vc{\th})$.
Moreover, $f_i$ takes $\1$ to $\vc{m_i}$ and
$\hcube_{\vc{c_0}}(\1)$ to $\hcube_{\vc{c_i}}(\vc{m_i})$ through an affine map.
\end{lemma}
\begin{proof}
Let us assume that $Q$ is embedded in $\Sph$ and that
$y_i$ is the vertex of $Q$ with angle $\pi \th_i$. Denote by $y'_i$ the point
of $\Sph$ antipodal to $y_i$.
The four sides of $Q$ lie on four geodesics in $\Sph$ that cut $\Sph$ in six convex quadrilaterals and eight convex triangles. All the quadrilaterals in this lemma are assembled from these pieces and the vertices of these quadrilaterals are chosen among the points 
$y_i$ and $y'_j$.
\begin{center}
\begin{figurehere}
\psfrag{Q}{\small{$Q$}}
\psfrag{y1}{\footnotesize{$y_1$}}
\psfrag{y2}{\footnotesize{$y_2$}}
\psfrag{y3}{\footnotesize{$y_3$}}
\psfrag{y4}{\footnotesize{$y_4$}}
\psfrag{y'1}{\footnotesize{$y'_1$}}
\psfrag{y'2}{\footnotesize{$y'_2$}}
\psfrag{y'3}{\footnotesize{$y'_3$}}
\psfrag{y'4}{\footnotesize{$y'_4$}}
\psfrag{x1}{\footnotesize{$x_1$}}
\psfrag{x2}{\footnotesize{$x_2$}}
\psfrag{x3}{\footnotesize{$x_3$}}
\psfrag{x4}{\footnotesize{$x_4$}}
\psfrag{(1)}{$(Q_1)$}
\psfrag{(2)}{$(Q_2)$}
\psfrag{(3)}{$(Q_3)$}
\psfrag{(4)}{$(Q_4)$}
\psfrag{(5)}{$(Q_5)$}
\psfrag{(6)}{$(Q_6)$}
\psfrag{(7)}{$(Q_7)$}
\includegraphics[width=0.9\textwidth]{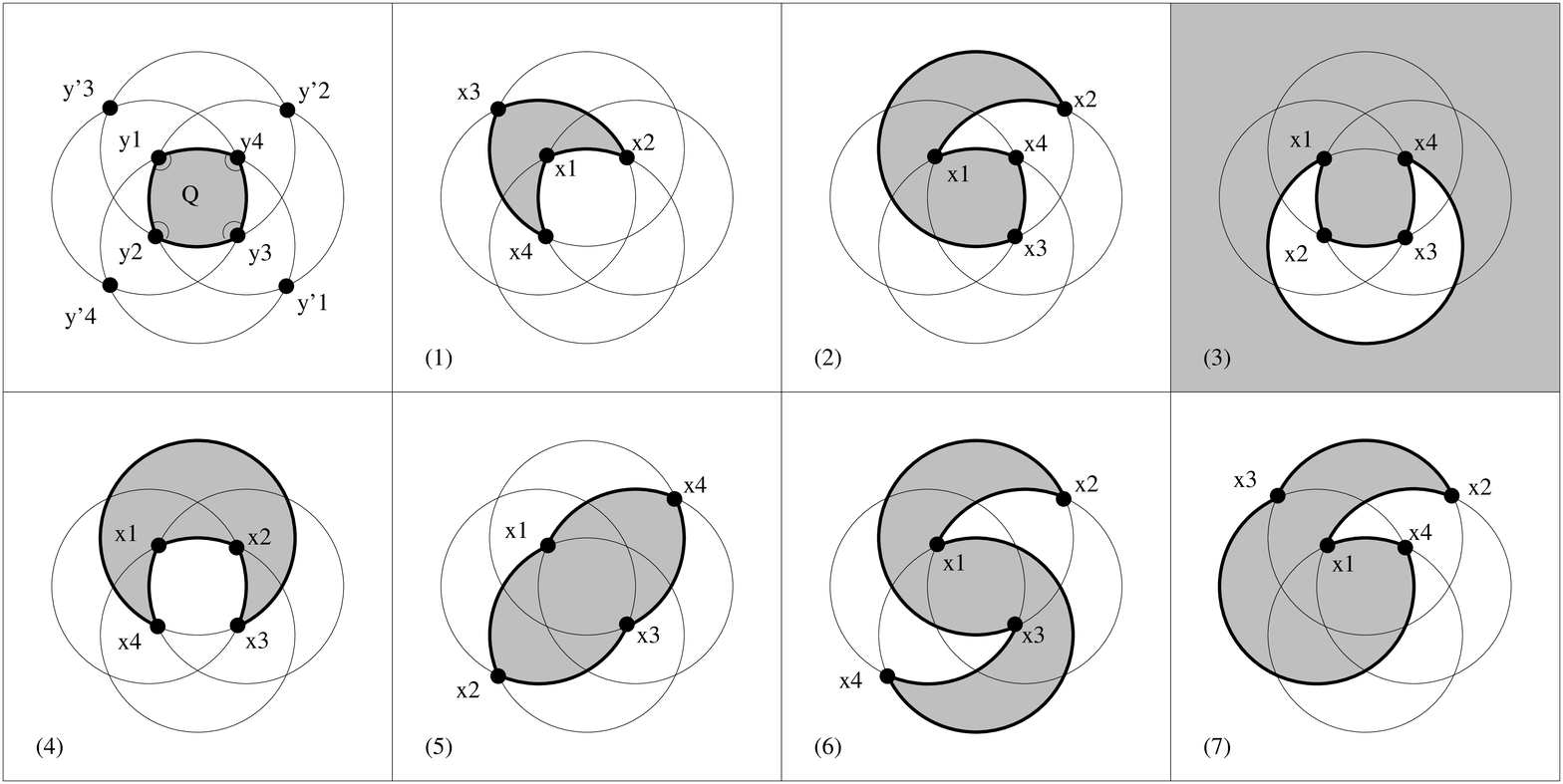}
\caption{{\small The seven basic non-convex quadrilaterals.}}\label{fig:quadrilaterals}
\end{figurehere}
\end{center}
The convex quadrilateral $Q$ and all of seven non-convex quadrilaterals $Q_i$ we wish to construct are shown in Figure \ref{fig:quadrilaterals}:
in all the cases we remove from $\Sph$ a point lying in the quadrilateral opposite to $Q$ and we draw the four great circles on which the edges of $Q$ lie.
%
%
%
%
%
%
%
%
\end{proof}




\begin{remark}\label{shortgeodesic} In quadrilaterals $(Q_1)$ and $(Q_2)$ the vertex $x_1$ is the only one with angle larger than $\pi$ and both adjacent sides (namely, $x_1x_2$ and $x_1x_4$) are shorter than $\pi$. In quadrilaterals $(Q_3),\dots,(Q_7)$ there are two opposite sides shorter than $\pi$ that join a vertex with an angle larger than $\pi$ with a vertex with an angle less than $\pi$. 
\end{remark}

We will now show that the angles of the quadrilaterals constructed in Lemma \ref{1big3small}
cover almost all points of $\tcube_{(\frac{3}{2},\frac{1}{2},\frac{1}{2},\frac{1}{2})}$
and $\tcube_{(\frac{3}{2},\frac{3}{2},\frac{1}{2},\frac{1}{2})}$.
%

\begin{corollary}[Non-convex quadrilaterals I]\label{three<1one<2}
Let $\vc{\th}$ be in the interior of  $\tcube_{(\frac{3}{2},\frac{1}{2},\frac{1}{2},\frac{1}{2})}$ but $\vc{\th}\ne (1+a,1-a,1-a,1-a)$ for all $a\in(0,1)$.
Then
for some permutation $\sigma\in\mathfrak{S}_4$
there exists a spherical quadrilateral with angles $\pi\cdot \vc{\th}_\sigma$.
\end{corollary}

\begin{proof} 
To prove this corollary we will use quadrilaterals of 
type $(Q_1)$ and $(Q_2)$.
Consider the set of points in $\tcube_{(\frac{3}{2},\frac{1}{2},\frac{1}{2},\frac{1}{2})}$ that can be represented 
by quadrilaterals of type $(Q_1)$ in Lemma \ref{1big3small}. From 
Lemma \ref{convexquad} it follows that these points are exactly 
those at distance less than $2$ from the point $\vc{m_1}=(1,0,1,0)$. 
In the same way, quadrilaterals of type $(Q_2)$ correspond to points of $\tcube_{(\frac{3}{2},\frac{1}{2},\frac{1}{2},\frac{1}{2})}$ at distance
less than $2$ from  the point $\vc{m_2}=(2,0,1,1)$. 

Now, the group of coordinate permutations preserving $\cube_{(\frac{3}{2},\frac{1}{2},\frac{1}{2},\frac{1}{2})}$ can be identified to $\mathfrak{S}_3\cong \mathrm{stab}(1)\subset\mathfrak{S}_4$.
The union of the orbits of the points $\vc{m_1}$ and $\vc{m_2}$ under this group consists
of the following six vertices of $\tcube_{(\frac{3}{2},\frac{1}{2},\frac{1}{2},\frac{1}{2})}$: 
\[
(1,0,1,0),\ (1,1,0,0), \ (1,0,0,1);  \quad (2,0,1,1), \ (2,1,0,1), \ (2,1,1,0).
\]
Hence six halves of $\tcube_{(\frac{3}{2},\frac{1}{2},\frac{1}{2},\frac{1}{2})}$
are covered by $\vc{\th}$ corresponding to spherical quadrilaterals (see Remark \ref{halfcube}). It is easy to see that  points of type $(1+a,1-a,1-a,1-a)$ are the only points that are not covered. These  are exactly the points in $\tcube_{(\frac{3}{2},\frac{1}{2},\frac{1}{2},\frac{1}{2})}$ that are on distance 
$2$ from the above six vertices. 
\end{proof}

\begin{corollary}[Non-convex quadrilaterals II]\label{two<1two<2}
Let $\vc{\th}$ be in the interior of $\tcube_{(\frac{3}{2},\frac{3}{2},\frac{1}{2},\frac{1}{2})}$.
Then
for some permutation $\sigma\in\mathfrak{S}_4$
there exists a spherical quadrilateral with angles $\pi\cdot \vc{\th}_\sigma$.
\end{corollary}

\begin{proof} 
The argument is similar to the one employed
in the proof of Corollary \ref{three<1one<2}
but in this case we use quadrilaterals of types $(Q_3)$, $(Q_4)$,
$(Q_5)$, $(Q_6)$ and $(Q_7)$.
By Lemma \ref{1big3small}, after taking coordinate permutations,
we see that $\vc{m_3},\dots,\vc{m_7}$ correspond to the eight even vertices
of $\tcube_{(\frac{3}{2},\frac{3}{2},\frac{1}{2},\frac{1}{2})}$, namely
\begin{gather*}
(2,2,1,1);\quad (1,1,0,0);\quad (1,1,1,1);\quad (2,2,0,0);\\
(2,1,0,1),\ (2,1,1,0),\ (1,2,0,1),\ (1,2,1,0)\, .
\end{gather*}
Thus, the construction of $(Q_3)$, $(Q_4)$,
$(Q_5)$, $(Q_6)$ and $(Q_7)$
provides quadrilaterals corresponding to points
in $\tcube_{(\frac{3}{2},\frac{3}{2},\frac{1}{2},\frac{1}{2})}$ at distance less than $2$ from all eight
even vertices.

\begin{center}
\begin{figurehere}
\psfrag{B'}{$B'$}
\psfrag{B''}{$B''$}
\psfrag{x1}{$x_1$}
\psfrag{x2}{$x_2$}
\psfrag{x3}{$x_3$}
\psfrag{x4}{$x_4$}
\psfrag{y''}{$y''$}
\psfrag{y'}{$y'$}
\psfrag{1/2}{{\scriptsize{$1/2$}}}
\psfrag{1}{{\scriptsize{$1$}}}
\includegraphics[width=0.3\textwidth]{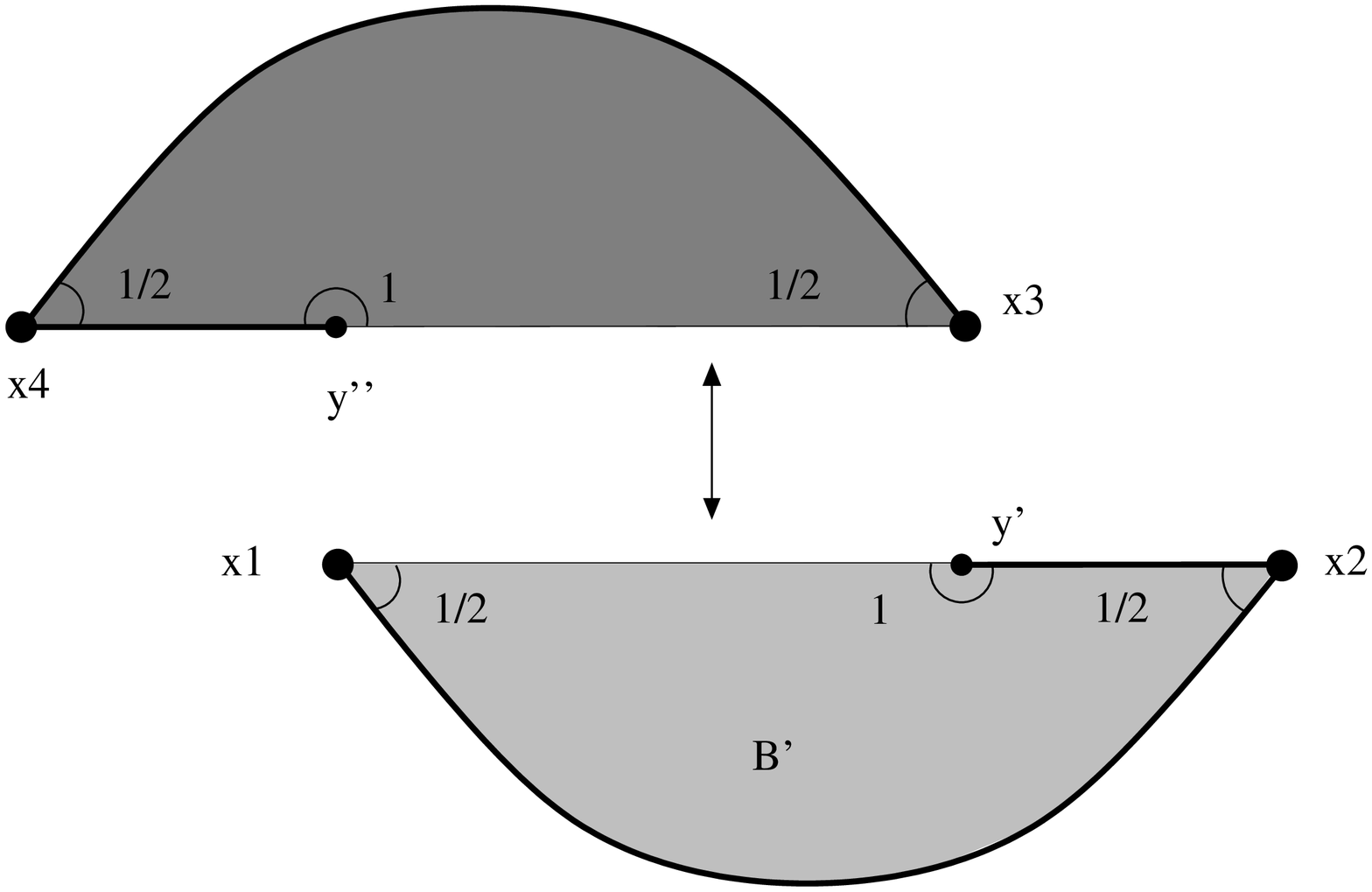}
\caption{{\small A quadrilateral with angles $\pi\cdot(\frac{3}{2},\frac{1}{2},\frac{3}{2},\frac{1}{2})$.}}\label{fig:quad-center}
\end{figurehere}
\end{center}

The only point at distance at least $2$ from all eight even vertices of
$\tcube_{(\frac{3}{2},\frac{3}{2},\frac{1}{2},\frac{1}{2})}$
is its center. 
In order to construct a quadrilateral with angles $\pi\cdot (\frac{3}{2},\frac{1}{2},\frac{3}{2},\frac{1}{2})$, consider two ordinary bigons $B'$ and $B''$ of angle $\pi/2$,
with vertices $x_1,x_2$ and $x_3,x_4$. 
Let $r\in(0,\pi)$ and pick a point $y'\in \pa B'$ at distance $r$ from $x_1$ and a point $y''\in\pa B''$ at distance $r$ from $x_3$. The wished quadrilateral is then obtained by gluing $B'$ and $B''$ via the isometric identification of $x_1 y'$ with $y'' x_3$ (see Figure \ref{fig:quad-center}).
%
%
\end{proof}

\subsubsection{Quadrilaterals immersed in $\Sph$}
In order to construct spheres with four conical points with angles
in $\tcube_{(\frac{5}{2},\frac{1}{2},\frac{1}{2},\frac{1}{2})}$,
we proceed as in the previous section.

\begin{lemma}[Three basic immersed quadrilaterals]\label{lemma:immersed}
Let $\vc{\th}\in\hcube_{\vc{c_0}}(\1)$ and consider the following table.
\[
\def\arraystretch{1.2}
\begin{array}{|c|c|c|c|}
\hline
i & f_i(\vc{\th}) & \vc{m_i} & \vc{c_i}\\
\hline
8 & (3-\th_1, 1-\th_4, \th_3, \th_2) & (2,0,1,1) & (\frac{5}{2},\frac{1}{2},\frac{1}{2},\frac{1}{2})\\
9 & (2+\th_1, 1-\th_2, \th_3, 1-\th_4) & (3,0,1,0) & (\frac{5}{2},\frac{1}{2},\frac{1}{2},\frac{1}{2})\\
10 & (2+\th_1,\th_2,\th_3,\th_4) & (3,1,1,1) &  (\frac{5}{2},\frac{1}{2},\frac{1}{2},\frac{1}{2})\\
\hline
\end{array}
\]
For every convex quadrilateral $Q$ with cyclically ordered angles $\pi\cdot\vc{\th}$
and for every $8\leq i\leq 10$
there exists a quadrilateral $Q_i$ with cyclically ordered angles $\pi\cdot f_i(\vc{\th})$.
Moreover, $f_i$ takes $\1$ to $\vc{m_i}$ and
$\hcube_{\vc{c_0}}(\1)$ to $\hcube_{\vc{c_i}}(\vc{m_i})$ through an affine map.
\end{lemma}

\begin{center}
\begin{figurehere}
\psfrag{x1}{{$x_1$}}
\psfrag{x2}{{$x_2$}}
\psfrag{x3}{{$x_3$}}
\psfrag{x4}{{$x_4$}}
\psfrag{(8)}{$(Q_8)$}
\psfrag{(9)}{$(Q_9)$}
\psfrag{(10)}{$(Q_{10})$}
\includegraphics[width=0.5\textwidth]{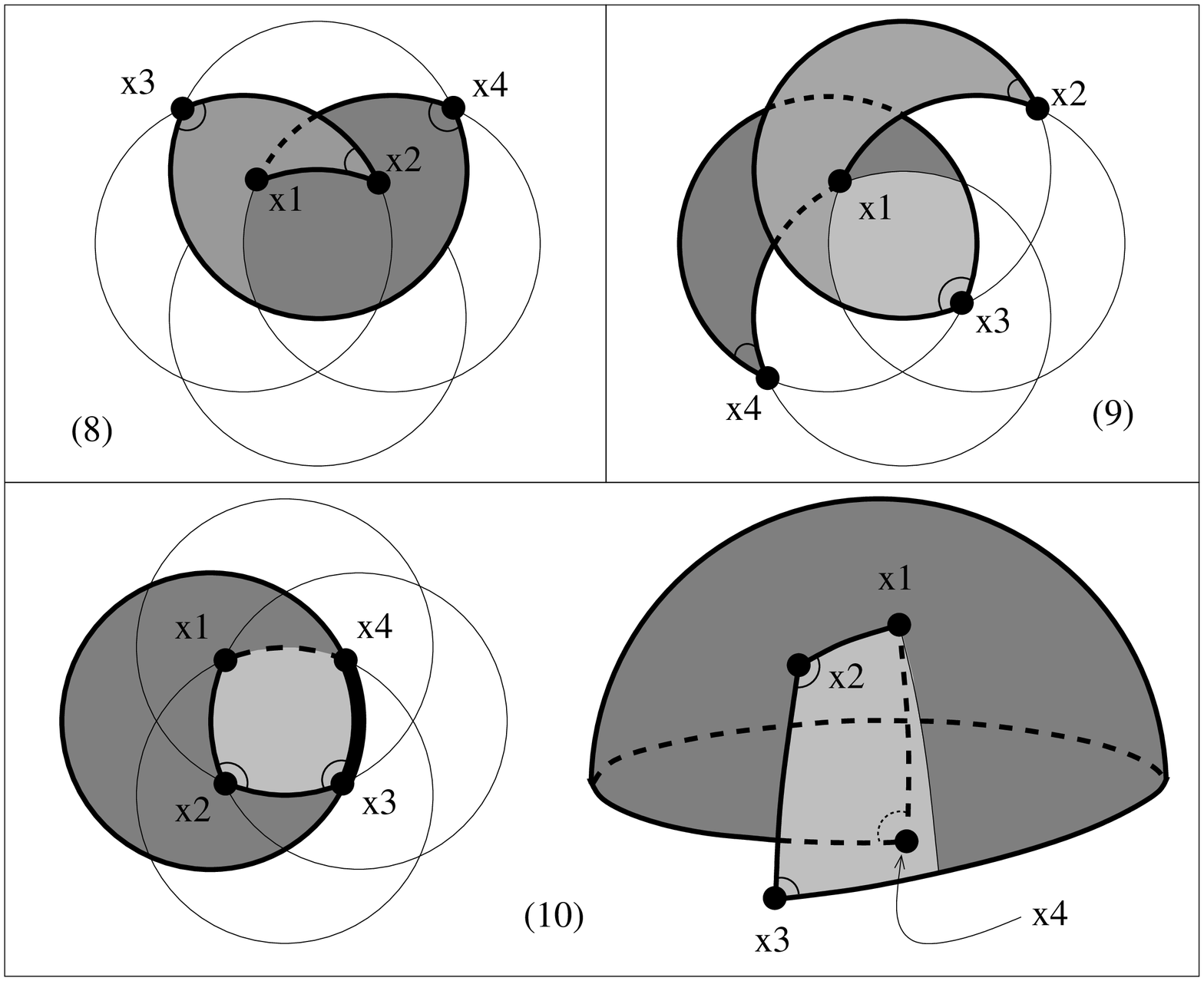}
\caption{{\small Three families of immersed quadrilaterals.}}\label{fig:immersed-quadr}
\end{figurehere}
\end{center}

\begin{proof}
The quadrilaterals are illustrated in Figure \ref{fig:immersed-quadr} as immersed in $\Sph$,
though we have included also another picture of $Q_{10}$ for clarity.

The existence of such quadrilaterals relies on Lemma \ref{1big3small}.
In fact, in order to construct $Q_8$ consider the quadrilateral
$Q_1$ and call $x'_1,x'_2,x'_3,x'_4$ its vertices and let
$B=B_{\th_2}$ be an ordinary bigon with vertices $y_4$ and $y'_4$.
Call $z$ one of the two points on the boundary of $B$ at distance
$|x'_4 x'_1|$ from $y'_4$.
The quadrilateral
$Q_8$ is obtained from $Q_1$ and $B$
by gluing the edge $x'_4 x'_1$ with $y'_4 z$,
so that $x'_4$ and $y'_4$ are identified to a smooth point
(which will not be marked as a vertex of $Q_8$): its vertices
are $x_1:=[x'_1]=[z]$, $x_2:=[x'_4]$, $x_3:=[x'_3]$ and
$x_4=:[y_4]$.

In a similar fashion, consider
$Q$ with vertices $x'_i$ and let
$B''=B_{1-\th_2}$ and $B'''=B_{1-\th_3}$ be two
ordinary bigons with vertices $x'',y''$ and $x''',y'''$ respectively.
Let $z''$ be a point on $\pa B''$ at distance
$|x'_2 x'_1|$ from $y''$ and let
$z'''$ be a point on $\pa B'''$ at distance
$|x'_3 x'_1|$ from $y'''$.
The quadrilateral $Q_9$ is then obtained
by gluing $x'_2 x'_1$ to $y''z''$ and
$x'_3x'_1$ to $y'''z'''$ and
calling $x_1:=[x'_1]=[z'']$, $x_2:=[x'']$,
$x_3:=[x'_3]$, $x_4:=[x''']$.

Finally, start again from the quadrilateral $Q$ with
vertices $x'_i$ and the triangle $T=T(1,|x'_1x'_4|,1-\th_4)$
with vertices $y_1,y_3,y_4$ of angles $\pi(2,1-\th_4,\th_4)$.
The wished $Q_{10}$ is obtained by identifying $y_1y_3$
to $x'_1x'_4$ and then calling
$x_1:=[x'_1]=[y_1]$, $x_2:=[y'_2]$, $x_3:=[x'_3]$
and $x_4:=[y_4]$.
\end{proof}

\begin{corollary}[Immersed quadrilaterals]\label{cor:immersed}
Let $\vc{\th}$ be in the interior of  $\tcube_{(\frac{5}{2},\frac{1}{2},\frac{1}{2},\frac{1}{2})}$ but $\vc{\th}\ne (2+a,a,a,a)$ for all $a\in(0,\frac{1}{2}]$.
Then for some permutation $\sigma\in\mathfrak{S}_4$
there exists a spherical quadrilateral with angles $\pi\cdot \vc{\th}_\sigma$.
Moreover, the vertex with angle larger than $2\pi$ can be joined
to any other vertex with a smooth geodesic of length strictly less than $\pi$.
\end{corollary}

\begin{proof} 
The argument is similar to the one employed
in the proof of Corollaries \ref{three<1one<2} and \ref{two<1two<2}.
This time we use quadrilaterals of types $(Q_8)$, $(Q_9)$, $(Q_{10})$.

By Lemma \ref{1big3small}, after taking coordinate permutations,
we see that $\vc{m_8},\vc{m_9},\vc{m_{10}}$ correspond to the seven even vertices
of $\tcube_{(\frac{5}{2},\frac{1}{2},\frac{1}{2},\frac{1}{2})}$, namely
\begin{gather*}
(3,1,1,1);\quad (2,1,1,0),\ (2,1,0,1),\ (2,0,1,1);\quad (3,1,0,0),\ (3,0,1,0),\ (3,0,0,1)\ .
\end{gather*}
Thus, the construction of $(Q_8)$, $(Q_9)$,
$(Q_{10})$ provides quadrilaterals corresponding to points
in $\tcube_{(\frac{3}{2},\frac{3}{2},\frac{1}{2},\frac{1}{2})}$ at distance less than $2$ from these seven vertices. The remaining points belong to the interval connecting the vertex $(2,0,0,0)$ with the centre of $\tcube_{(\frac{5}{2},\frac{1}{2},\frac{1}{2},\frac{1}{2})}$. The last assertion can be checked by direct inspection.
\end{proof}

\subsubsection{Sporadic families of 4-punctured spheres}

Though most 4-punctured spheres can be obtained by doubling
spherical quadrilaterals, it seems from Lemma \ref{three<1one<2} and Lemma \ref{lemma:immersed} that there are two 1-parameter families of $\vc{\th}\in\RR^4$ 
such that we are not able to construct quadrilaterals with angles
$\pi\cdot\vc{\th}$.
Thus, for such families of angles, here we present ad-hoc constructions.

\begin{center}
\begin{figurehere}
\psfrag{S}{$S_a$}
\psfrag{Si}{$\Sigma_a$}
\psfrag{x1}{$x_1$}
\psfrag{x2}{$x_2$}
\psfrag{x3}{$x_3$}
\psfrag{x4}{$x_4$}
\psfrag{y1}{$y_1$}
\psfrag{y2}{$y_2$}
\psfrag{y3}{$y_3$}
\psfrag{a}{$\alpha$}
\psfrag{b}{$\beta$}
\psfrag{p}{$p$}
\includegraphics[width=0.5\textwidth]{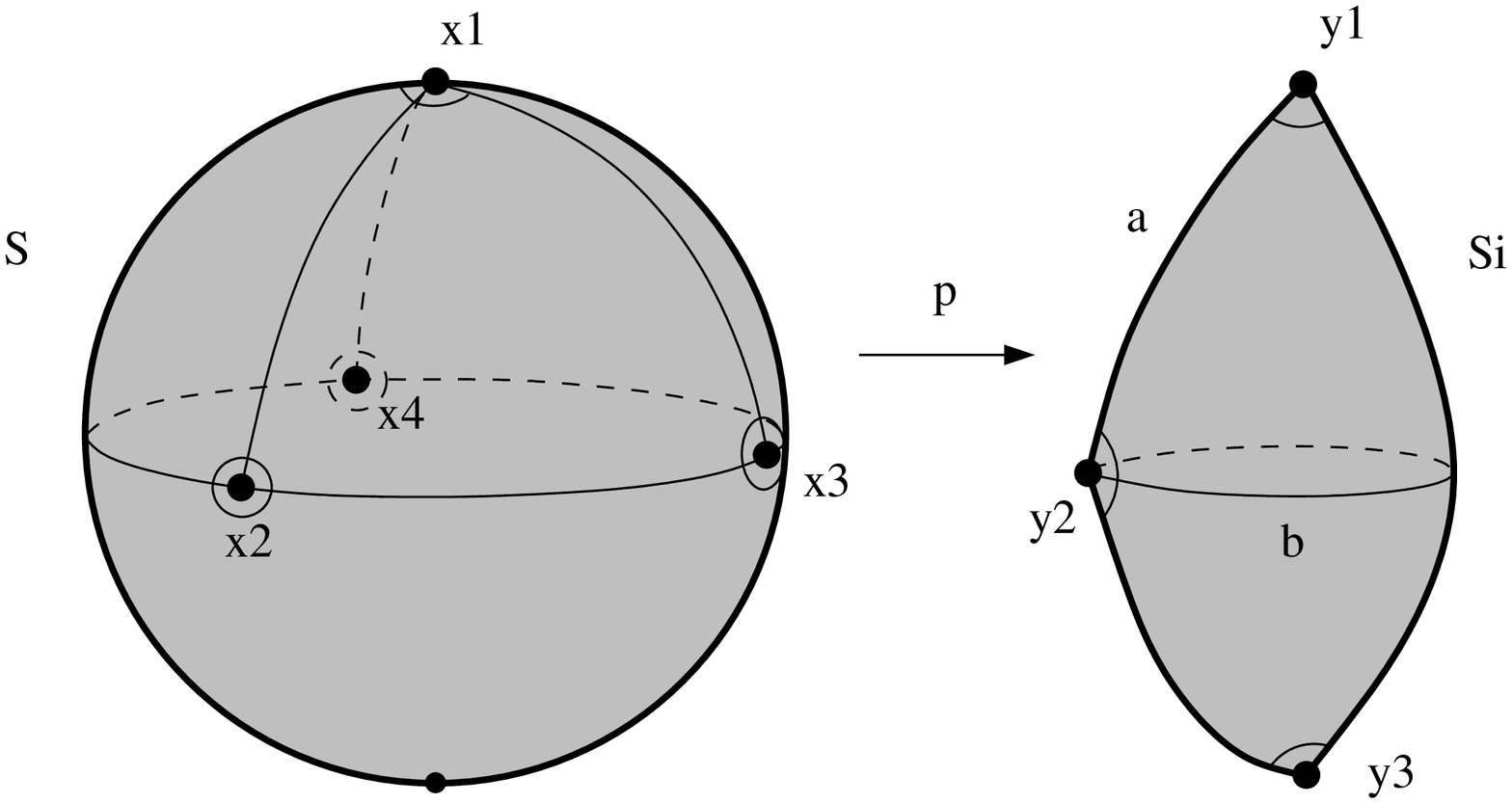}
\caption{{\small The sporadic sphere $S_a$ with the six paths $\gamma_{ij}$.}}\label{fig:cover3-1}
\end{figurehere}
\end{center}

\begin{lemma}[Sporadic $4$-punctured spheres]\label{exceptions}
\begin{itemize}
\item[(a)]
For any $a\in (0,1)$ there exists 
a 4-punctured sphere $S_a$ and a spherical metric on it with 
conical singularities $x_1,x_2,x_3,x_4$ of angles $2\pi\cdot (1+a,1-a,1-a,1-a)$.
\item[(b)]
For any $b\in (0,\frac{1}{2})$ there exists 
a 4-punctured sphere $S_b$ and a spherical metric on it with 
conical singularities $x_1,x_2,x_3,x_4$ of angles $2\pi\cdot (2+b,b,b,b)$.
\item[(c)]
There exists a 4-punctured sphere $S$ 
and a spherical metric on it with 
conical singularities $x_1,x_2,x_3,x_4$ of angles $2\pi\cdot(\frac{5}{2},\frac{1}{2},\frac{1}{2},\frac{1}{2})$.
\end{itemize}
Moreover, on the spheres constructed in all three cases there is a smooth geodesic from $x_1$ to $x_j$ of length strictly less than $\pi$ for $j=2,3,4$.
\end{lemma}

\begin{proof}
As for part (a), notice that
for every $a\in (0,1)$ there exists a sphere $\Sigma_a$ with 
three conical points $y_1,y_2,y_3$ of angles $2\pi\cdot\vc{\th}(a)$, where
$\vc{\th}(a):=(\frac{1+a}{3}, 1-a, \frac{1}{3})$.
Indeed, $\vc{\th}(0)=(\frac{2}{3}, 1, \frac{1}{3})$ and $\vc{\th}(1)=(\frac{1}{3}, 0, \frac{1}{3})$ lie on the boundary of the simplex formed by the angle vectors
$\vc{\th}\in\RR^3$ corresponding to spheres with three angles less than $2\pi$
and $\vc{\th}(a)$ lies strictly inside such a simplex for $a\in(0,1)$.

Now consider the cyclic cover $p: S_a\rar \Sigma_a$ of degree three cover
branched over $y_1$ and $y_3$. The $x_1=p^{-1}(y_1)$ is a point of 
angle $2\pi(1+a)$ and $p^{-1}(y_3)$ is a smooth point, which will not be labelled. Moreover,
$p^{-1}(y_2)$ consists of three points of angle $2\pi(1-a)$, which we label by
$x_2,x_3,x_4$ (see Figure \ref{fig:cover3-1}). Thus, $(S_a,x_1,x_2,x_3,x_4)$ is our wished spherical surface. Three geodesics that joint $x_1$ with $x_j$ are preimages on $S_a$ of the shortest geodesic in $\Sigma_a$ joining $y_1$ and $y_2$.\\

The proof of (b) is entirely analogous. 
As above, for any $b\in (0,\frac{1}{2})$ 
there exists a sphere $\Sigma_b$ with 
three conical points $y_1,y_2,y_3$ of angles $2\pi\cdot\vc{\th}(b)$, where
$\vc{\th}(b):=(\frac{2+b}{3}, b, \frac{1}{3})$.
Indeed, $\vc{\th}(0)=(\frac{2}{3}, 0, \frac{1}{3})$ and $\vc{\th}(\frac{1}{2})=(\frac{5}{6}, \frac{1}{2}, \frac{1}{3})$ lie on the boundary of the simplex formed by the angle vectors
$\vc{\th}\in\RR^3$ corresponding to spheres with three angles less than $2\pi$,
and $\vc{\th}(b)$ lies strictly inside it for $b\in(0,\frac{1}{2})$.
Now, as before take the cyclic cover of degree three $S_b\rar\Sigma_b$ branched at $y_1$ and $y_3$.\\
\begin{center}
\begin{figurehere}
\psfrag{D}{$D$}
\psfrag{D'}{$D'$}
\psfrag{y1}{$y_1$}
\psfrag{y2}{$y_2$}
\psfrag{y'3}{$y'_3$}
\psfrag{y'2}{$y'_2$}
\psfrag{y'1}{$y'_1$}
\psfrag{w'}{$w'$}
\psfrag{z'}{$z'$}
\psfrag{z}{$z$}
\includegraphics[width=0.5\textwidth]{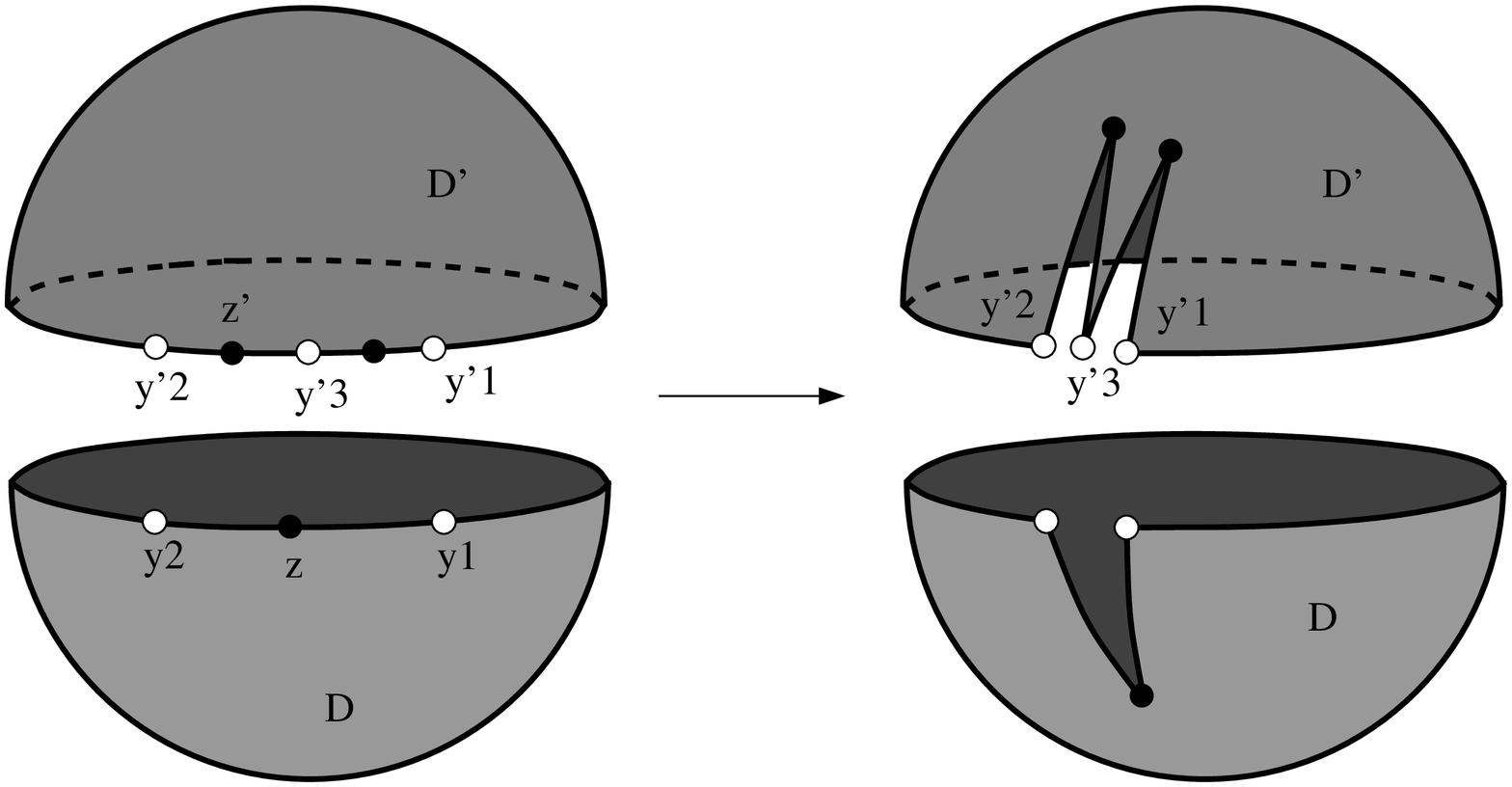}
\caption{{\small A sphere with angles $2\pi\cdot(\frac{5}{2},\frac{1}{2},\frac{1}{2},\frac{1}{2})$.}}\label{fig:quad-5-2}
\end{figurehere}
\end{center}
About (c), consider two hemispheres $D$ and $D'$ and let $\ell\in (0,\pi)$.
On $\pa D$ pick points $y_1,z,y_2$ (in this cyclic order with respect to the orientation induced on $\pa D$)
in such a way that $|y_1 y_2|=\ell$ and $z$ is the midpoint of $y_1y_2$;
on $\pa D'$ pick points $y'_2,z',y'_3,w',y'_1$ (in this cyclic order) in such a way that $|y'_2y'_1|=\ell$
and that $z'$ is the midpoint of $y'_2y'_3$ and $w'$ is the midpoint of $y'_3y'_1$.
As in Figure \ref{fig:quad-5-2}, the wished sphere $S$ is obtained by identifying $y_1z$ to $z y_2$ on $D$
and $y'_2z'$ to $z'y'_3$ and $y'_3w'$ to $w'y'_1$ on $D'$, and finally $y'_1y'_2$ to $y_1y_2$:
its marked points are $x_1:=[y_i]=[y'_i]$, $x_2=[z]$, $x_3=[w']$, $x_4=[z']$.
The last assertion is clear by construction.
\end{proof}

\subsubsection{Spheres with $\th_1, \th_2<2$, $\th_3,\th_4<1$ and with $1<\th_1<2$, $\th_2,\th_3,\th_4<1$}

Here we derive corollaries from the statements proven in the previous sections.
Denote by $\Pi^4$ the box $[1,2]\times[0,2]\times[0,1]^2\subset\RR^4$.

\begin{proposition}[4-punctured spheres with non-integral angles in $\Pi^4$]\label{sumupfour} 
For every $\vc{\th}\in\inte(\Pi^4\cap \Aang^4)$ with no integral coordinate
there exists a sphere $S$ with a spherical metric $g$ and conical singularities $x_1,x_2,x_3,x_4$ of angles $2\pi\cdot\vc{\th}$, which satisfies the following properties:
\begin{itemize}
\item[(a)]
there exist six simple 
paths 
$\{\gamma_{ij}\,|\,1\leq i<j\leq 4\}$
that have no inner points of intersection and such that
$\gamma_{ij}$ joins $x_i$ and $x_j$;
\item[(b)]
either $\gamma_{13}$ or $\gamma_{14}$ is a geodesic shorter than $\pi$;
\item[(c)]
the metric $g$ has non-coaxial holonomy.
\end{itemize}
\end{proposition}
%

\begin{proof}
We will first construct the spheres and the paths $\gamma_{ij}$ and then will prove that their holonomy is not coaxial.\\

{\it{Construction of spheres with six paths.}}\\
Notice that a $\vc{\th}\in\inte(\Pi^4\cap \Aang^4)$
with no integral coordinate must belong to the interior
either of
$\tcube_{(\frac{3}{2},\frac{1}{2},\frac{1}{2},\frac{1}{2})}$ or of $\tcube_{(\frac{3}{2},\frac{3}{2},\frac{1}{2},\frac{1}{2})}$.
Hence, by doubling quadrilaterals $Q$ constructed in Corollaries \ref{three<1one<2} and \ref{two<1two<2}, we cover all the cases apart from the exceptional family treated in Lemma \ref{exceptions}(a). 

In order to find the six paths $\gamma_{ij}$, we proceed as follows.

Consider first the case of $S$ obtained as a double quadrilateral $DQ=Q\sqcup\ol{Q}/\!\sim$.
Take four geodesic paths $\gamma_{12},\gamma_{23},\gamma_{34},\gamma_{14}$ corresponding to the edges of $Q$ (or; equivalently of $\ol{Q}$)
and choose a simple path $\gamma_{13}$ inside $Q$ and a simple path
$\gamma_{24}$ inside $\ol{Q}$. All these paths will be simple, since the quadrilateral is embedded in $\Sph$. Moreover it follows from Remark \ref{shortgeodesic} that either $\gamma_{13}$ or $\gamma_{14}$ can be chosen to be a geodesic shorter than $\pi$.

Consider now the exceptional spheres $S_a$ with $a\in(0,1)$, constructed in Lemma \ref{exceptions}(a).
%
By Lemma \ref{lemma:double-triangles}, the surface
$\Sigma_a$ can be constructed by doubling of a spherical 
triangle $T_a$ with vertices $y_1,y_2,y_3$ and angles $\pi\cdot (\frac{1+a}{3}, 1-a, \frac{1}{3})$.
Choose a point $z\in T_a$ in the interior of the side $y_1y_3$. Now consider
the following two paths on $\Sigma_a=DT_a$: the geodesic
$\alpha$ determined by the edge $y_1y_2$ of $T_a$ 
and the geodesic $\beta$ obtained 
by doubling of the geodesic segment $y_2z$ contained in $T_a$.
Clearly, $\beta$ is a simple loop on $\Sigma_a$ based at $y_2$ and it is easy to
see that $\alpha$ is shorter than $\pi$.
If $p:S_a\rar\Sigma_a$ is the triple cyclic cover branched at $y_1,y_3$,
then we define $\gamma_{12},\gamma_{13},\gamma_{14}$ to be
the preimages of $\alpha$ and $\gamma_{23},\gamma_{34},\gamma_{24}$ to
be the preimages of $\beta$ through $p$.
Since $|\alpha|<\pi$, we have that both $\gamma_{13}$ and $\gamma_{14}$ are shorter than $\pi$.

This completes the proof of claims (a) and (b).\\

{\it{Non-coaxiality.}}\\  
Since $\th_i\notin \ZZ$, non-coaxiality of the holonomy of the spherical
surfaces $S$ just constructed follows from Lemma \ref{noncoaxhol}
if we are able to find two distinct conical points $x_i,x_j$ on $S$
joined by a smooth geodesic $\gamma$ of length $\ell$ with $\ell\notin \pi \ZZ$.
By the above property (b), we can choose $\gamma$ to be either the path $\gamma_{13}$ or $\gamma_{14}$. This proves (c).
\end{proof}

Analogously, we have the following.

\begin{proposition}[4-punctured spheres with non-integral angles in $\tcube_{(\frac{5}{2},\frac{1}{2},\frac{1}{2},\frac{1}{2})}$.]\label{prop:5/2} 
For every $\vc{\th}$ in the interior of $\tcube_{(\frac{5}{2},\frac{1}{2},\frac{1}{2},\frac{1}{2})}$
there exists a sphere $S$ with a spherical metric $g$ and conical singularities $x_1,x_2,x_3,x_4$ of angles $2\pi\cdot\vc{\th}$, which satisfies the following properties:
\begin{itemize}
\item[(a)]
$x_1$ and $x_2$ are joined by a smooth geodesic of length strictly less than $\pi$;
\item[(b)]
the metric $g$ has non-coaxial holonomy.
\end{itemize}
\end{proposition}
\begin{proof}
The wished spherical surface $S$ is either obtained 
from Lemma \ref{exceptions}(b-c) or by doubling
the quadrilaterals constructed in Corollary \ref{cor:immersed}.
In either case, property (a) is satisfied.

Property (b) then follows from Lemma \ref{noncoaxhol},
since $\th_1,\th_2\notin\ZZ$.
\end{proof}

\subsubsection{Existence of spheres with 4 conical points and non-integral angles}

In this section we finally prove Theorem \ref{thm:four}.
We will construct the desired spherical surfaces
starting from those produced in Proposition  \ref{sumupfour}
and applying the gluing operations of Proposition \ref{adding4piand2pi2pi}. 
Since these surgeries do not change the holonomy, non-coaxiality
of the new metrics follows from Proposition \ref{sumupfour}.


\begin{notation}
Let $\vc{e_1},\dots,\vc{e_4}$ be the standard generators of $\ZZ^4$
and define $\vc{e_{kl}}:=\vc{e_k}+\vc{e_l}$ for $1\le k<l\le 4$.
The six elements $\vc{e_{kl}}$ generate 
the semigroup
\[
\Gamma^4:=
\left\{
\vc{p}\in\ZZ_{\geq 0}^4\,\Big|\,
\|\vc{p}\|_1\in 2\ZZ\ \text{and}
\ 2p_j\le \|\vc{p}\|_1 \ \text{for all $j=1,\dots,4$}
\right\}\, .
\]
\end{notation}

\begin{proof}[Proof of Theorem \ref{thm:four}]
Let $\vc{m}$ be a point in $\ZZ_{\geq 0}^4$ such that
$\vc{\th}\in\cube_{\vc{c}}$, where
$\vc{c}=\vc{m}+(\frac{1}{2},\frac{1}{2},\frac{1}{2},\frac{1}{2})$.
Without loss of generality, we can assume $m_1\ge m_2\ge m_3\ge m_4\geq 0$, and also $m_1\geq 1$.

We will now treat two cases separately.\\

{\it{Case (a): $m_1\leq m_2+m_3+m_4$.}}\\
Suppose first that $\|\vc{m}\|_1\in 2\ZZ$
so that $\vc{m}\in\Gamma^4$.

Since $m_1\geq 1$, we have $m_2\geq 1$ and so $\vc{m'}=\vc{m}-\vc{e_{12}}\in \Gamma^4$. As a consequence, we have a presentation
\[
\vc{m'}=m_{12}\vc{e_{12}}+\dots+m_{34}\vc{e_{34}},
\]
for suitable $m_{ij}\in\ZZ_{\geq 0}$ and so $\vc{\th}-\vc{m'}\in[1,2]\times[1,2]\times[0,1]^2\subset\Pi^4$.
Since $\|\vc{m'}\|\in 2\ZZ$, the vector $\vc{\th}-\vc{m'}\in\Aang^4$ and so
$\vc{\th}-\vc{m'}\in\inte(\Pi^4\cap\Aang^4)$ and it has no integral coordinate.

By Proposition \ref{sumupfour} there exists a sphere $S'$ with four conical points of angles $2\pi\cdot (\vc{\th}-\vc{m'})$ and six simple paths $\gamma_{ij}$ joining $x_i$ and $x_j$
that may only intersect at their endpoints.
 
The wished spherical surface $S$ is obtained by performing the surgery
described in Lemma \ref{adding4piand2pi2pi}(a) along these paths, gluing the sphere $S'$ with $m'_{ij}$ copies of $\Sph\setminus {\dev}_{\gamma_{ij}}$ along each $\gamma_{ij}$ for all $1\leq i<j\leq 4$.
This settles the case $\|\vc{m}\|_1\in 2\ZZ$.

To treat the case when $\|\vc{m}\|_1$ is odd,
it is enough to choose $\vc{m'}=\vc{m}-\vc{e_1}$ and to observe that $\vc{m'}\in \Gamma^4$ and that
$\vc{\th}-\vc{m'}\in [1,2]\times[0,1]^3\subset\Pi^4$.
Then the above argument carries on.\\
 
{\it{Case (b): $m_1>m_2+m_3+m_4$ and $\|\vc{m}\|_1$ odd.}}\\
Since $\vc{\th}-\vc{m}+\vc{e_{1}}\in[1,2]\times[0,1]\times[0,1]^2\subset\Pi^4$
and $\vc{m}-\vc{e_{1}}$ is even, we have $\vc{\th}-\vc{m}+\vc{e_{1}}\in
\inte(\Pi^4\cap \Aang^4)$ with no integral coordinate. Moreover, $(m_1-1)-m_2-m_3-m_4=2d$ with $d\in \ZZ_{\geq 0}$.
By Proposition \ref{sumupfour} there exists a sphere $S'$ with conical
points $x_1,\dots,x_4$ of angles
$2\pi\cdot (\vc{\th}-\vc{m}+\vc{e_{1}})$;
%
moreover, 
either $\gamma_{13}$ or $\gamma_{14}$ is a geodesic
of length less than $\pi$, which we will denote by $\gamma$. 

Hence, we can apply to $S'$ the surgery operation 
described in Lemma \ref{adding4piand2pi2pi}(b) along $\gamma$,
thus producing a sphere
$S''$ with angles $2\pi\cdot(\vc{\th}-\vc{m}+(2d+1)\vc{e_1})$.
%
%
 
Finally, we obtain our wished spherical surface by applying the operation of 
Lemma \ref{adding4piand2pi2pi}(a) to 
$S''$ by gluing 
$m_2$ copies of $\Sph\setminus {\dev}_{\gamma_{12}}$ along $\gamma_{12}$, 
$m_3$ copies of $\Sph\setminus {\dev}_{\gamma_{13}}$ along $\gamma_{13}$, and
$m_4$ copies of $\Sph\setminus {\dev}_{\gamma_{14}}$ along $\gamma_{14}$.\\

{\it{Case (c): $m_1>m_2+m_3+m_4$ with $m_2>0$ and $\|\vc{m}\|_1$ even.}}\\
We proceed analogously to case (b).
Since $\vc{\th}-\vc{m}+\vc{e_{12}}\in[1,2]\times[1,2]\times[0,1]^2\subset\Pi^4$
and $\vc{m},\vc{e_{12}}$ are even, we have $\vc{\th}-\vc{m}+\vc{e_{12}}\in
\inte(\Pi^4\cap \Aang^4)$ with no integral coordinate.
Moreover, $m_1-m_2-m_3-m_4=2d$ with $d\in \ZZ_{\geq 0}$. 
By Proposition \ref{sumupfour} there exists a sphere $S'$ with conical
points $x_1,\dots,x_4$ of angles
$2\pi\cdot (\vc{\th}-\vc{m}+\vc{e_{12}})$;
moreover, 
either $\gamma_{13}$ or $\gamma_{14}$ is a geodesic
of length less than $\pi$, which we will denote by $\gamma$. 

Hence, we can apply to $S'$ the surgery operation 
described in Lemma \ref{adding4piand2pi2pi}(b) along $\gamma$,
thus producing a sphere
$S''$ with angles $2\pi\cdot(\vc{\th}-\vc{m}+\vc{e_{12}}+2d\vc{e_1})$.
%
%
 
Finally, we obtain our wished spherical surface by applying the operation of 
Lemma \ref{adding4piand2pi2pi}(a) to 
$S''$ by gluing 
$(m_2-1)$ copies of $\Sph\setminus {\dev}_{\gamma_{12}}$ along $\gamma_{12}$, 
$m_3$ copies of $\Sph\setminus {\dev}_{\gamma_{13}}$ along $\gamma_{13}$, and
$m_4$ copies of $\Sph\setminus {\dev}_{\gamma_{14}}$ along $\gamma_{14}$.\\

{\it{Case (d): $m_1>m_2+m_3+m_4$ with $m_2=0$ and $\|\vc{m}\|_1$ even.}}\\
Clearly, we must have $m_1=2+2d$ with $d\in\ZZ_{\geq 0}$ and $m_2=m_3=m_4=0$
Hence, $\vc{\th}-2d\vc{e_1}$ belongs to
the interior of $\tcube_{(\frac{5}{2},\frac{1}{2},\frac{1}{2},\frac{1}{2})}$.

By Proposition \ref{prop:5/2}, 
there exists a sphere $S'$ with conical
points $x_1,\dots,x_4$ of angles
$2\pi\cdot (\vc{\th}-2d\vc{e_{1}})$;
moreover, 
$x_1$ and $x_2$ are joined by a smooth geodesic $\gamma$ of length strictly
less than $\pi$.

Hence, we can apply to $S'$ the surgery operation 
described in Lemma \ref{adding4piand2pi2pi}(b) along $\gamma$,
thus producing a sphere
$S$ with angles $2\pi\cdot\vc{\th}$.\\

In all cases, the spherical surface $S'$ has non-coaxial holonomy
by Proposition \ref{sumupfour} in cases (a-b-c) and
by Proposition \ref{prop:5/2} in case (d).
Thus, the surface $S$ constructed
performing gluing operations as in Lemma \ref{adding4piand2pi2pi}
is non-coaxial too.
\end{proof}

\subsection{Splitting conical points}\label{sec:induction}
%

The aim of this section is to complete the proof of Theorem \ref{existence}, by showing the following.
%
%

\begin{theorem}[Existence of spherical metrics for $n\geq 5$]\label{thm:five}
Assume $n\geq 4$ and let $\th_1,\dots,\th_n$ be real numbers that 
both the positivity constraints (\ref{gauss-bonnet}) and
the holonomy constraints
(\ref{theinquality}) strictly. If $n=4$, then also assume that one $\th_i$ is integral.\\
Then there exists a sphere $S$ endowed with a spherical metric 
with $n$ conical singularities of angles $2\pi\th_1,\dots,2\pi\th_n$
and non-coaxial holonomy. Moreover, such a metric is deformable.
\end{theorem}

%
%

Clearly, this immediately leads to our main result.

\begin{proof}[Proof of Theorem \ref{existence}]
The statement for $n=3$ follows from Theorem \ref{membrane}, since
$3$-punctured spheres are obtained by doubling spherical triangles.
The statement for $n=4$ when all $\th_1,\dots,\th_4$ are not integers
has already been proven
in Theorem \ref{thm:four}
and for $n\geq 5$ is the content
of Theorem \ref{thm:five}.

Notice that, if $n=4$ and $\vc{\th}$ satisfies positivity and strict holonomy constraints, then at most one $\th_i$ can be integer.
This case is also taken care by Theorem \ref{thm:five} and so the proof is complete.
\end{proof}


Theorem \ref{thm:five} is based on an inductive argument,
whose key step can be formulated as follows.

\begin{lemma}[Inductive step]\label{lemma:induction}
Let $\vc{\th}\in\RR^n_+$ and 
$M:\RR^n\rar\RR^{n-1}$ be a merging operation of
type $M_{(i+j)}$ or $M_{(i-j)}$, where $i,j\in\{1,\dots,n\}$ are two distinct indices. If $M=M_{(i-j)}$, then assume that $\th_j\notin\ZZ$.\\
Suppose that there exists a sphere $S'$ endowed with a non-coaxial, angle-deformable spherical metric $g'$ with $n-1$ conical singularities $x'_1,\dots,x'_{n-1}$ of angles $2\pi\cdot\vc{\th'}$, where 
$\vc{\delta'}:=M(\vc{\delta})$ and $\vc{\delta},\vc{\delta'}$
are the defects associated to $\vc{\th},\vc{\th'}$.
Then there exists a sphere $S$ endowed with a non-coaxial, angle-deformable spherical metric $g$ with $n$ conical singularities of angles $2\pi\cdot\vc{\th}$.
\end{lemma}
\begin{proof}
Since the metric $g$ is angle-deformable, there exists a neighbourhood
$\Nang'\subset\RR^{n-1}$ of $\vc{\th'}$ and a continuous family of metrics
$\Nang'\ni\vc{\nu'}\mapsto g'_{\vc{\nu'}}$ on $S'$ such that
$g'_{\vc{\th'}}=g'$ and $g'_{\vc{\nu'}}$ has singularities of angles
$2\pi\cdot\vc{\nu'}$.
Up to shrinking $\Nang'$, we can assume that all $g_{\vc{\nu'}}$ are $\varepsilon$-wide at $x'_{n-1}$ and with non-coaxial holonomy.

{\it{Case $M=M_{(i+j)}$.}}\\
By Proposition \ref{sumtriangle}, there exist an $|\eta|<\varepsilon/2$ and
an $(y_1,y_2)$-angle-deformable spherical triangle $(T,g'')$ with
vertices $y_1,y_2,y_3$, angles $\pi(\th_i,\th_j,\th_i+\th_j-1+\eta)$,
which is $\pi(1-\varepsilon/2)$-wide at $y_3$.
Thus, there exists a neighbourhood $\Nang''$ of $(\th_i,\th_j)\in\RR^2$,
a function $\theta_3:\Nang''\rar\RR$ with $\theta_3(\th_i,\th_j)=\th_i+\th_j-1+\eta$ and a continuous family
$\Nang''\ni\vc{\nu''}\mapsto g''_{\vc{\nu''}}$ of spherical metrics on $T$
such that $g''_{(\th_i,\th_j)}=g''$ and $g''_{\vc{\nu''}}$ has conical
angles $\pi(\nu''_1,\nu''_2,\theta_3(\vc{\nu''}))$.

By continuity, there exists a neighbourhood $\Nang$ of $\vc{\th}\in\RR^n$ such that $\vc{\nu''}(\vc{\nu}):=(\nu_i,\nu_j)\in \Nang''$ and
$\vc{\nu'}(\vc{\nu}):=(\nu_1,\dots,\widehat{\nu}_i,\dots,\widehat{\nu}_j,\dots,\nu_n,\theta_3(\nu_i,\nu_j))\in \Nang'$
for all $\vc{\nu}\in \Nang$.

For every such $\vc{\nu}\in \Nang$, consider the surface $(S,g_{\vc{\nu}})$ obtained by gluing $(S',g'_{\vc{\nu'}(\vc{\nu})})$
and the double of $(T,g''_{\vc{\nu''}(\vc{\nu})})$ at the conical points
$x'_{n-1}\in S'$ and $[y_3]\in DT$ according to Lemma 
\ref{lemma:gluing-conical}. This construction provides
a continuous family $\Nang\ni \vc{\nu}\mapsto g_{\vc{\nu}}$ of
spherical metrics on $S$ with conical points of angles $2\pi\cdot\vc{\nu}$.
Moreover, the holonomy of $g_{\vc{\nu}}$ is non-coaxial, since it
contains that of $g'_{\vc{\nu'}(\vc{\nu})}$, which is non-coaxial.

{\it{Case $M=M_{(i-j)}$.}}\\
Since $\delta_j\notin\ZZ$, we can apply
Proposition \ref{differencetriangle} to obtain
an $|\eta|<\varepsilon/2$ and
an $(y_1,y_2)$-angle-deformable spherical triangle $(T,g'')$ with
vertices $y_1,y_2,y_3$, angles $\pi(\th_i,\th_j,\th_i-\th_j-1+\eta)$,
which is $\pi(1-\varepsilon/2)$-wide at $y_3$.
The proof then works as in the previous case.
\end{proof}

Finally, the argument is completed as follows.

\begin{proof}[Proof of Theorem \ref{thm:five}]
Let $\vc{\delta}=(\th_1-1,\dots,\th_n-1)$ as usual.\\

{\it{Case $n=4$ and $\th_i\in\ZZ$.}}\\
Since $\vc{\delta}\in \inte(\Aang^4)$, for any $j\neq i$
the operation $M=M_{(i+j)}$ satisfies $\vc{\delta'}=M(\vc{\delta})\in\inte(\Aang^3)$
by Lemma \ref{lemma:int-merging}.
By Theorem \ref{membrane}, there exists
a non-coaxial angle-deformable spherical triangle with angles $2\pi(\delta'_1+1,\delta'_2+1,\delta'_3+1)$ and so
we can apply Lemma \ref{lemma:induction}, thus obtaining
the wished non-coaxial angle-deformable spherical surface
of genus $0$ with angles $2\pi(\th_1,\dots,\th_4)$.

Together with Theorem \ref{thm:four}, this settles the case $n=4$.\\

{\it{Case $n\geq 5$: induction.}}\\
Assume now that the statement holds for 
$(n-1)$-punctured spheres: 
we will prove it for $n$-punctured spheres. 

Since $n\geq 5$, by Theorem \ref{5points} there exists
a merging operation $M$ such that $\vc{\delta'}:=M(\vc{\delta})$
belongs to $\inte(\Aang^{n-1})$.
By inductive hypothesis, there exists a surface $S'$ of genus $0$ with a non-coaxial angle-deformable spherical metric and $n-1$ conical singularities of angles $2\pi(\delta'_1+1,\dots,\delta'_{n-1}+1)$.
The conclusion now follows
by Lemma \ref{lemma:induction}.
\end{proof}

%
%
%
%
%
%
%
%
%

\setcounter{secnumdepth}{-1}

\section{List of symbols}

{\small
\begin{multicols*}{2}
\Trick
\begin{supertabular}{lp{5cm}}
$\vc{e_i}$ & $i$-th vector of the standard basis of $\RR^n$\\
$\vc{\1}$ & vector $\vc{e_1}+\dots+\vc{e_n}\in\RR^n$\\
$\vc{\th}$ & angle vector $(\th_1,\th_2,\dots,\th_n)\in\RR^n$\\
$\ol{\vc{\th}}$ & reduced angle vector $\vc{\ol{\th}}\in [0,2)^n$ with $\vc{\th}-\vc{\ol{\th}}\in 2\ZZ$\\
$\vc{\delta}$ & defect vector $\vc{\th}-\vc{\1}\in\RR^n$\\
$\ol{\vc{\delta}}$ & reduced defect vector $\ol{\vc{\th}}-\vc{\1}\in[-1,1)^n$\\
$\Nang$ & small neighbourhood of $\vc{\th}$ in $\RR^n$\\
$\vc{\nu}$ & angle vector in $\Nang$\\
$d_1(\cdot,\cdot)$ & standard $\ell^1$ distance in $\RR^n$\\
$\|\cdot\|_1$ & standard $\ell^1$ norm in $\RR^n$\\
$\ZZ^n_o$ & subset of odd-integral vectors, i.e. $\vc{m}\in\ZZ^n\subset\RR^n$ with $\|\vc{m}\|$ odd\\
$\Hang^n$ & locus of $\vc{\delta}\in\RR^n$ such that $d_1(\vc{\delta},\ZZ^n_o)\geq 1$\\
$\Pang^n$ & locus of $\vc{\delta}\in (-1,+\infty)^n$ such that $
\sum_i \delta_i>-2$\\
$\Aang^n$ & intersection of $\Hang^n$ and $\Pang^n$\\
$M_{(i+j)}$ & algebraic positive merging operation\\
$M_{(i-j)}$ & algebraic negative merging operation\\
$\cube^n$ & unit cube with integral vertices in $\RR^n$\\
$\vc{c}$ & center of a unit cube in $\RR^n$\\
$\cube_{\vc{c}}$ & unit cube with center $\vc{c}\in\RR^n$\\
$\tcube^n$ & truncated cube $\cube^n\cap \Hang^n$\\
$\tcube_{\vc{c}}$ & truncated cube $\cube_{\vc{c}}\cap \Hang^n$\\
$\hcube_{\vc{c}}(\vc{m})$ & half truncated cube with center $\vc{c}$ and vertex $\vc{m}$\\
$\vc{\delta}^\pi$ & radial projection of $\vc{\delta}\in\tcube_{\vc{c}}$ onto $\partial\tcube_{\vc{c}}$ (for $\vc{\delta}\neq\vc{c}$)\\
$\Sph,\Spt$ & unit spheres endowed with the standard metric\\
$T^1\Sigma$ & unit tangent bundle to $\Sigma$\\
$\dev$ & developing map of a simply connected surface\\
$\dev_\gamma$ & developing map of a path $\gamma$\\
$\dot{S}$ & complement of the conical points $x_1,\dots,x_n$ in $S$\\
$\rho$ & holonomy representation in $\SO(3,\RR)$\\
$\hat{\rho}$ & standard lift of the holonomy representation to $\SU(2)$\\
$c_p$ & constant loop based at the point $p$\\
$\gamma_j$ & loop that simply winds about the $j$-th marked point\\
$U_j$ & matrix in $\SU(2)$ representing the holonomy along $\gamma_j$\\
$v_j$ & vertex of a broken geodesic on $\Spt$\\
$s_j$ & side of a broken geodesic on $\Spt$\\
$\ell_j$ & length of the side $s_j$ of a broken geodesic on $\Spt$\\
$|x_i x_{i+1}|$ & length of the edge between $x_i$ and $x_{i+1}$ in a spherical polygon\\
$DS$ & surface obtained by doubling the surface with boundary $S$\\
$B_\alpha(r)$ & standard open $r$-neighbourhood of a vertex of angle $\pi\alpha$ in a spherical polygon\\
$\ol{B}_\alpha(r)$ & standard closed $r$-neighbourhood of a vertex of angle $\pi\alpha$ in a spherical polygon\\
$S_\alpha(r)$ & standard open $r$-neighbourhood of a point of angle $2\pi\alpha$ in a spherical surface\\
$\ol{S}_\alpha(r)$ & standard closed $r$-neighbourhood of a point of angle $2\pi\alpha$ in a spherical surface\\
$B_\alpha$ & ordinary spherical bigon with angles $\pi\alpha$\\
$S_\alpha$ & double of $B_\alpha$\\
$B(d,\ell)$ & exceptional spherical bigon with angles $\pi d$ at distance $\ell$\\
$T(d,\ell,\alpha)$ & spherical triangle with sides $\ell,\ell,2\pi d$ and angles $\pi\alpha,\pi(1-\alpha),2\pi d$\\
$U_y(r)$ & complement in a spherical surface of the neighbourhood\\
& of the cone point $y$ of angle $2\pi\alpha$ isometric to $B_\alpha(r)$\\
$S\#_r S'$ & surface obtained by surgery at conical points\\
$S_\gamma\#_{\gamma'}S'$ & surface obtained by surgery along paths\\
\end{supertabular}
\end{multicols*}
}



\bibliographystyle{amsplain}
\bibliography{draft-biblio}

\end{document}